\documentclass[11pt,reqno,twoside]{amsart}
\usepackage{graphicx}
\usepackage{amssymb}
\usepackage{epstopdf}
\usepackage[asymmetric,top=2.5cm,bottom=2.8cm,left=2.38cm,right=2.38cm]{geometry}
\usepackage{hyperref}
\usepackage{booktabs} 
\usepackage{array} 
\usepackage{paralist} 
\usepackage{verbatim} 
\usepackage{subfig} 
\usepackage{tabularx}
\usepackage{amsmath,amsfonts,amsthm,mathrsfs,amssymb,cite}
\usepackage[usenames]{color}
\usepackage{bm}
\usepackage{color}

\usepackage{xcolor}
\hypersetup{
  colorlinks,
  linkcolor={blue!80!black},
  urlcolor={blue!80!black},
  citecolor={blue!80!black}
}

\newtheorem{thm}{Theorem}[section]

\newtheorem{lem}{Lemma}[section]
\newtheorem{prop}{Proposition}[section]
\theoremstyle{definition}
\newtheorem{defn}{Definition}[section]
\theoremstyle{remark}

\newtheorem{rem}{Remark}[section]
\numberwithin{equation}{section}

\newcommand{\bmf}[1]{{\mathbf{#1}}}
\def\bsi{{\mathrm{i}}}

\def\Oh{{\mathcal  O}}

\newcommand{\rmd}{\mathrm{d}}

\title[Generalized Holmgren's principle to Lam\'e operator and applications]{On generalized Holmgren's principle to the Lam\'e operator with applications to inverse elastic problems}

\author{Huaian Diao}
\address{School of Mathematics and Statistics, Northeast Normal University,
Changchun, Jilin 130024, China.}
\email{hadiao@nenu.edu.cn}

\author{Hongyu Liu}
\address{Department of Mathematics, City University of Hong Kong, Kowloon, Hong Kong SAR, China.}
\email{hongyu.liuip@gmail.com; hongyliu@cityu.edu.hk}

\author{Li Wang}
\address{School of Mathematics and Statistics, Northeast Normal University,
Changchun, Jilin 130024, China.}
\email{1322771004@qq.com}

\date{} 
\begin{document}
\maketitle

\begin{abstract}
Consider the Lam\'e operator $\mathcal{L}(\bmf u) :=\mu \Delta \bmf u+(\lambda+\mu) \nabla(\nabla \cdot \bmf u )$ that arises in the theory of linear elasticity. This paper studies the geometric properties of the (generalized) Lam\'e eigenfunction $\bmf u$, namely $-\mathcal{L}(\bmf u)=\kappa \bmf u$ with $\kappa\in\mathbb{R}_+$ and $\bmf u\in L^2(\Omega)^2$, $\Omega\subset\mathbb{R}^2$. We introduce the so-called homogeneous line segments of $\bmf{u}$ in $\Omega$, on which $\bmf{u}$, its traction or their combination via an impedance parameter is vanishing. We give a comprehensive study on characterizing the presence of one or two such line segments and its implication to the uniqueness of $\bmf{u}$. The results can be regarded as generalizing the classical Holmgren's uniqueness principle for the Lam\'e operator in two aspects. We establish the results by analyzing the development of analytic microlocal singularities of $\bmf{u}$ with the presence of the aforesaid line segments. Finally, we apply the results to the inverse elastic problems in establishing two novel unique identifiability results. It is shown that a generalized impedance obstacle as well as its boundary impedance can be determined by using at most four far-field patterns. 
Unique determination by a minimal number of far-field patterns is a longstanding problem in inverse elastic scattering theory.



\end{abstract}

\section{Introduction}

Consider the following partial differential operator (PDO) acting on a $\mathbb{C}^2$-valued function $\mathbf{u}=({u}_\ell({\mathbf x}))_{\ell=1}^2$, $\mathbf{x}=(x_\ell)_{\ell=1}^2\in\mathbb{R}^2$:
\begin{equation}\label{eq:pdo1}
\mathcal{L}(\bmf u) :=\mu \Delta \bmf u+(\lambda+\mu) \nabla(\nabla \cdot \bmf u ).
\end{equation}
$\mathcal{L}$ is known as the Lam\'e operator that arises in the theory of linear elasticity. $\lambda,\mu$ are the Lam\'e constants satisfying the following strong convexity condition
\begin{equation}\label{eq:convex}
\mu>0 \text { and } \lambda+\mu>0.
\end{equation}
Let $\Omega\subset\mathbb{R}^2$ be an open set. $\bmf u=(u_\ell)_{\ell=1}^2\in L^2(\Omega)^2$ is said to be a (generalised) Lam\'e eigenfunction if 
\begin{equation}\label{eq:lame}
-\mathcal{L} (\bmf{u})=\kappa \bmf{u}\ \text { in }\ \Omega,\ \ \kappa \in\mathbb{R}_+.
\end{equation}
It is noted that there is no boundary condition prescribed on $\partial\Omega$ for $\bmf{u}$. 
\eqref{eq:lame} arises in the study of the time-harmonic wave scattering. Let $\Gamma_h\Subset\Omega$ be an open connected line segment, where $h\in\mathbb{R}_+$ signifies the length of the line segment. Let
\begin{equation}\label{eq:nutau}
	\bmf{\nu }=(\nu_1,\nu_2)^\top \mbox{ and }
 \boldsymbol{\tau}=(-\nu_2,\nu_1)^\top 
\end{equation}
respectively, signify the unit normal and tangential vectors to $\Gamma_h$. The traction $T_\nu\bmf{u}$ on $\Gamma_h$ is defined by
\begin{equation}\label{eq:Tu}
 T_{\bmf{\nu }}\bmf{u}=2\mu\partial_{\bmf{\nu }}\bmf{u}+\lambda\bmf{\nu} \left(\nabla \cdot \bmf{u}\right)+\mu\boldsymbol{\tau}(\partial_{2}u_{1}-\partial_{1}u_{2}), \end{equation}
where 
$$
\nabla\bmf u:=\begin{bmatrix}
	\partial_1 u_1 & \partial_2 u_1 \cr
	\partial_1 u_2 & \partial_2 u_2 
\end{bmatrix}, 
$$
$\partial_\nu\bmf u:=\nabla\bmf u \cdot  \nu$,  $\partial_j u_i:=\partial u_i/\partial x_j$. It is noted that $\bmf{u}|_{\Gamma_h}$ and $T_\nu\bmf{u}|_{\Gamma_h}$ consist of the Cauchy data on $\Gamma_h$ to $\bmf{u}$ in \eqref{eq:lame} . 

We recall the classical Holmgren's theorem for an elliptic PDO $\mathcal{P}$ with real-analytic coefficients (cf. \cite{TF}). If $\mathcal{P}\bmf{u}$ is real analytic in a connected open neighbourhood of $\Omega$, then $\bmf{u}$ is also real-analytic. The Holmgren's theorem applied to $\bmf{u}$ in \eqref{eq:lame}, we immediately see that $\bmf{u}$ is real-analytic in $\Omega$. Suppose that 
\begin{equation}\label{eq:cond1l}
\bmf{u}=\mathbf{0}\quad\mbox{and}\quad T_\nu\bmf{u}=\mathbf{0}\quad\mbox{on}\ \ \Gamma_h, 
\end{equation}
then by the Cauchy-Kowalevski theorem, one readily has that $\bmf{u}\equiv 0$ in $\Omega$. This is known as the Holmgren's uniqueness principle. It also holds when $\Gamma_h$ is replaced to be an analytic curve. In this paper, we shall generalize the Holmgren's principle with the Cauchy data on a line segment to the Lam\'e operator $\mathcal{L}+\kappa$ in two aspects. First, we note that in \eqref{eq:cond1}, both Cauchy data are required to vanish on the line segment $\Gamma_h$. We ask whether this is the minimal/optimal requirement to ensure the uniqueness of $\bmf{u}$. Can the Holmgren's principle still hold, say if one only requires that 
\[
\bmf{u}(\bmf{x}_0)=\mathbf{0}\quad\mbox{and}\quad T_\nu \bmf{u}|_{\Gamma_h}=\mathbf{0},
\]
where $\bmf{x}_0\in\Gamma_h$ is a single point? Clearly, in general, this cannot be true for a generic PDO. However, it is one of the interesting discoveries of the present paper that one of the two homogeneous conditions in \eqref{eq:cond1l} can indeed be replaced by a certain point-value condition. Second, we view \eqref{eq:cond1l} as the existence of two line segments $\Gamma_h^\pm$ such that: (i) $\bmf{u}|_{\Gamma_h^-}=\bmf{0}$ and $T_\nu\bmf{u}|_{\Gamma_h^+}=\mathbf{0}$; (ii) $\angle (\Gamma_h^-, \Gamma_h^+)=\pi$. Hence, a natural generalization is to consider the case that the two line segments are not of a straight intersection, namely, $\angle(\Gamma_h^-, \Gamma_h^+)\neq \pi$. In such a case, we can also establish a certain uniqueness principle for the solution to \eqref{eq:lame}. It is interesting to point out that for the latter generalization, the Cauchy data of $\bmf{u}$ are no longer prescribed on an analytic curve. Furthermore, we would like to emphasize that for both cases mentioned above, we also include the more general Robin-type condition into our study, namely $\big(\bmf{u}+\eta T_\nu\bmf{u}\big)|_{\Gamma_h}=\mathbf{0}$, which is known as an impedance condition with $\eta$ called an impedance parameter. We refer to the above discoveries as the generalized Holmgren's principle to the Lam\'e operator. The implication of the generalized principle to the uniqueness of a solution to the elastic problem \eqref{eq:lame} is obvious. Moreover, our study is clearly related to the geometric structures of the (generalized) Lam\'e eigenfunctions in \eqref{eq:lame}, which is also a central topic in the spectral theory of PDOs; see \cite{CDLZ,CDLZ2} and the references therein for more related discussions.

According to our discussion above, the results obtained are clearly of independent interest for their own sake in the PDE theory of elasticity and the spectral theory of the Lam\'e operator. Moreover, as an interesting practical application of our theoretical findings, we apply the results to the inverse scattering problem of determining an elastic obstacle as well as its possible surface impedance parameter by a minimal/optimal number of far-field measurements. This is a challenging problem with a strong applied background. In its abstract formulation, the problem can be roughly described as a nonlinear operator equation,
\begin{equation}\label{eq:ipa1}
\mathcal{F}(\Omega, \eta)=\mathcal{M}(\hat x; \bmf{u}^i_j),\quad \hat x\in\mathbb{S}^{1}, j=1,2,\ldots, N,
\end{equation}
where the scattering map $\mathcal{F}$ is defined by a certain PDE system in the exterior of a domain $\Omega$. $\eta$ is a boundary impedance parameter on $\partial\Omega$. Through solving the aforementioned PDE system, the scattering map $\mathcal{F}$ sends $\Omega$ and $\eta$ to a real-analytic function $\mathcal{M}$ on the unit sphere, which signifies the observation data. This correspondence also depends on $\bmf{u}_j^i$, $j=1,2,\ldots, N$, known as the incident fields, that account for the number of measurements in the practical scenario. We shall give more relevant details about \eqref{eq:ipa1} in Section~\ref{sect:6}. We are mainly concerned with the unique identifiability issue for \eqref{eq:ipa1}. That is, we aim to establish the unique one-to-one correspondence between the target object $(\Omega, \eta)$ and the measurement data $\mathcal{M}$, particularly with the minimal/optimal number of measurements. Geometrically speaking, a single measurement, namely $\mathcal{M}(\hat x)$, $\hat x\in\mathbb{S}^1$, corresponding to a single incident field $\bmf{u}^i$ (or at most a few), may serve as a global parametrization for $\partial\Omega$. However, there is very limited progress in the literature on this challenging geometrical problem. The more recent progress is concerned with the case that $\Omega$ is of general polygonal shape \cite{ElschnerYama2010,LiuXiao}. The mathematical machinery therein is mainly based on certain reflection and path arguments, which cannot deal with the more challenging case that $\eta$ is not identically $0$ or $\infty$. Using the generalized Holmgren's principle established in this paper, we can provide a different and unified approach in tackling with the geometrical inverse problem\eqref{eq:ipa1} in the case that $\Omega$ is of general polygonal shape with at most a few measurements. More importantly, our method can deal with the more challenging case that $\eta$ is finite and not identically zero. We derive a comprehensive study for this geometrical inverse problem. It is mentioned in passing that unique determination by a minimal number of far-field patterns is a longstanding problem in the inverse scattering theory. We refer to \cite{AR,CDLZ,CDLZ2,CY,CK,CK18,LPRX,LRX,Liu-Zou,Liu-Zou3} and the references therein for related studies for the inverse acoustic and electromagnetic wave scattering problems. In addition to the application to the inverse problem, we believe that the generalized Holmgren's principle may find more interesting applications in different contexts. 

Finally, we would like to briefly discuss about the technicality of our study. We shall be mainly based on analyzing the microlocal singularities of the solution $\bmf{u}$ to \eqref{eq:lame} due to the presence of the homogeneous line segments discussed earlier. Clearly, the singularities are developed across the aforementioned line segments and are of analytic type. In the case that there are two intersecting line segments with a non-straight intersecting angle, it is not surprising that the singularities are developed at the intersecting point. However, we shall show that the singularities can even be developed across a single line segment, which are really subtle and tricky to capture. In this paper, we mainly focus on the two-dimensional case. As can be seen that even in the two dimensions, the analyses are highly technical and lengthy with tedious calculations. We shall present the extensions to the three dimensions as well as to the case with Cauchy data on an analytic curve instead of a straight line segment in forthcoming articles. 

The result of the paper is organized as follows. Sections 2--4 are devoted to establishing the generalized Holmgren's principle in different scenarios. Section 5 presents the unique identifiability results for the inverse elastic obstacle problem \eqref{eq:ipa1}.

\section{Auxiliary results}\label{sect:2}

We first introduce two important definitions. 

\begin{defn}\label{def:1}
 Let $\bmf u=(u_\ell)_{\ell=1}^2$ be a generalized Lam\'e eigenfunction to \eqref{eq:lame} associated with an eigenvalue $\kappa\in\mathbb{R}_+$. An open and connected line segment $\Gamma_{h}\Subset\Omega$ is called {\it a rigid line} of $\bmf{u}$ if $\bmf{u}|_{\Gamma_h}=\bmf{0}$; {\it a traction-free line} if $T_{\bmf{\nu}}\bmf{u}|_{\Gamma_h}=\bmf{0}$; and an {\it impedance line} if
 \begin{equation}\label{eq:im line}
(T_{\bmf{\nu}}\bmf{u}+\eta \bmf{u})\big|_{\Gamma_h}=\bmf{0},
 \end{equation}
 where $\eta\in\mathbb{C}$ is constant and referred to as an impedance parameter. Set ${\mathcal R}_\Omega^{\kappa}  $, ${\mathcal T}_\Omega^{\kappa}  $  and  ${\mathcal I}_\Omega^{\kappa}  $ to respectively denote the sets of rigid, traction-free and impedance lines in $\Omega$ of $\bmf{u}$.

 \end{defn}
 
\begin{defn}\label{def:generalized line} Recall that the unit normal vector $\nu $ and the tangential vector $\boldsymbol{\tau}$ to $\Gamma_h$ are defined in \eqref{eq:nutau}, respectively. 
Define
\begin{subequations}
\begin{align}
	{\mathcal S}\left(   {\mathcal R}_\Omega^{\kappa} \right):=& \{\Gamma_h \in  {\mathcal R}_\Omega^{\kappa}  ~|~ \exists \bmf{x}_0 \in \Gamma_h \mbox{ such that }  \boldsymbol{\tau} \cdot \partial_{\nu} \bmf{u} |_{\bmf{x}={\bmf{x}}_0 }=0  \},\label{eq:def1} \\
	{\mathcal S}\left(   {\mathcal T}_\Omega^{\kappa} \right):=& \{\Gamma_h \in  {\mathcal T}_\Omega^{\kappa}  ~|~ \exists \bmf{x}_0 \in \Gamma_h \mbox{ such that }  \bmf{u}(\bmf{x}_0)=\bmf{0} \mbox{ and }\boldsymbol{\tau} \cdot \partial_{\nu} \bmf{u} |_{\bmf{x}={\bmf{x}}_0 }=0  \},\label{eq:def2}\\
	{\mathcal S}\left(   {\mathcal I}_\Omega^{\kappa} \right):=& \{\Gamma_h \in  {\mathcal I}_\Omega^{\kappa}  ~|~ \exists \bmf{x}_0 \in \Gamma_h \mbox{ such that }  \bmf{u}(\bmf{x}_0)=\bmf{0} \mbox{ and }\boldsymbol{\tau} \cdot \partial_{\nu} \bmf{u} |_{\bmf{x}={\bmf{x}}_0 }=0 \},\label{eq:def3}
\end{align}
\end{subequations}
where ${\mathcal S}\left(   {\mathcal R}_\Omega^{\kappa} \right)$, ${\mathcal S}\left(   {\mathcal T}_\Omega^{\kappa} \right)$ and ${\mathcal S}\left(   {\mathcal I}_\Omega^{\kappa} \right)$ are named as the sets of the \emph {singular rigid, singular traction-free } and \emph{singular impedance} lines of $\bmf{u}$ respectively.
Let ${\mathcal S}( \Omega )=	{\mathcal S}\left(   {\mathcal R}_\Omega^{\kappa} \right)\cup {\mathcal S}\left(   {\mathcal T}_\Omega^{\kappa} \right) \cup {\mathcal S}\left(   {\mathcal I}_\Omega^{\kappa} \right) $ be a set of singular lines of $\bmf{u}$ in $\Omega$.
\end{defn}

It is noted that compared to the homogeneous lines introduced in Definition~\ref{def:1}, the singular lines in Definition~\ref{def:generalized line} are further required to satisfy a number of conditions on a specific point. In what follows, we shall prove that if a (generalised) Lam\'e eigenfunction $\bmf{u}$ to \eqref{eq:lame} possesses a singular line in $\Omega$, then $\bmf{u}$ is identically zero. We prove this by quantitatively characerizing $\bmf{u}$ in the phase space across the lines. This is the reason that we call them (microlocally) singular lines. Furthermore, we show that the generic intersections of the homogeneous lines of Definition~\ref{def:1} shall also generate microlocal singularities, which prevent the occurrence of such intersections unless $\bmf{u}$ is trivially zero. In this article, we provide a comprehensive characterization of all those cases. To our best knowledge, those results are new to the literature.

\begin{figure}
	\centering
	\includegraphics[width=0.3\textwidth]{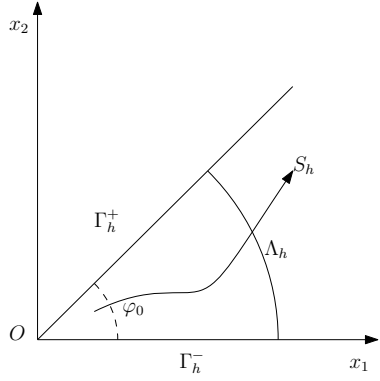}
	\caption{Schematic of the geometry of two intersecting lines with an angle $\varphi_0$ with $0< \varphi_0\leq \pi $.}
	\label{fig1}
\end{figure}

Next we introduce the geometric setup of our study. Consider two line segments respectively defined by (see Fig.~\ref{fig1} for a schematic illustration)
\begin{equation}\label{eq:gamma_pm}
	\begin{split}
		\Gamma_h^+&=\{\bmf{x} \in \mathbb R^2~|~\bmf{x}=r\cdot (\cos \varphi_0, \sin \varphi_0 )^\top ,\quad 0\leq r\leq h,\quad 0<\varphi_0\leq 2\pi   \}, \\
		\Gamma_h^-&=\{\bmf{x} \in \mathbb R^2~|~\bmf{x}=r\cdot (1, 0 )^\top,\quad 0\leq r\leq h   \},\ \ h\in\mathbb{R}_+. 
	\end{split}
	\end{equation}
Clearly, the intersecting angle between $\Gamma_h^+$ and $\Gamma_h^-$ is
\begin{equation}\label{eq:angle1}
\angle(\Gamma_h^+,\Gamma_h^{-})=\varphi_0, \quad 0< \varphi_0 \leq 2\pi.
\end{equation}
It is noted that if $\varphi_0=\pi$ or $2\pi$, $\Gamma_h^+$ and $\Gamma_h^-$ are actually lying on a same line. In such a case, the intersection between $\Gamma_h^+$ and $\Gamma_h^-$ is said to be degenerate. In our subsequent study, $\Gamma_h^\pm$ shall be the homogeneous lines in Definition~\ref{def:1} or the singular lines in Definition~\ref{def:generalized line}. In fact, for any two of such lines that are intersecting in $\Omega$ (or one line in the degenerate case), since the PDO $\mathcal{L}$ defined in \eqref{eq:lame} is invariant under rigid motions, one can always have two lines as introduced in \eqref{eq:gamma_pm} after a straightforward coordinate transformation such that the homogeneous conditions in Definitions~\ref{def:1} and \ref{def:generalized line} are still satisfied on $\Gamma_h^\pm$. We assume that $h\in\mathbb{R}_+$ is sufficiently small such that $\Gamma_h^\pm$ are contained entirely in $\Omega$. Moreover, if $\Gamma_h^\pm$ are impedance lines, we assume that the impedance parameters on $\Gamma_h^\pm$ are respectively two constants $\eta_1$ and $\eta_2$.  As also noted before that $\bmf{u}$ is analytic in $\Omega$, it is sufficient for us to consider the case that $0<\varphi_0\leq \pi$. In fact, if $\pi  < \varphi_0 \leq 2\pi$, we see that $\Gamma_h^+$ belongs to the half-plane of $x_2<0$ (see Fig.~\ref{fig1}). Let $\widetilde{\Gamma}_h^+$ be the extended line segment of length $h$ in the half-plane of $x_2>0$. By the analytic continuation, we know that $\widetilde{\Gamma}_h^+$ is of the same type of $\Gamma_h^+$, namely $\bmf{u}$ satisfies the same homogeneous condition on $\widetilde{\Gamma}_h^+$ as that on $\Gamma_h^+$. Hence, instead of studying the intersection of $\Gamma_h^+$ and $\Gamma_h^-$, one can study the intersection of $\widetilde{\Gamma}_h^+$ and $\Gamma_h^-$. Clearly, $\angle(\widetilde{\Gamma}_h^+, \Gamma_h^-)\in (0,\pi]$. 

Let $B_h$ be the central disk of radius $h \in \mathbb{R}_+$. Let $\Gamma^\pm$
signify the infinite extension of $\Gamma_{h}^\pm$ in the half-space $x_2\geq 0$.  Consider the open sector
\begin{equation}\label{eq:K}
	\mathcal{K} =\left\{\bmf{x}=(x_1,x_2)  \in \mathbb{R}^{2}~ |~ \bmf{x}\neq \bmf{0},\quad  0<\arg \left(x_{1}+\mathrm{i} x_{2}\right)<\varphi_{0}\right\},\quad \mathrm{i}:=\sqrt{-1},
\end{equation}
which is formed by the two half-lines $\Gamma^-$ and $\Gamma^+$. In the sequel, we set
\begin{equation}\label{eq:sh}
	S_h=\mathcal{K}\cap B_h,
\end{equation}
where $\partial S_h=\Gamma_h^+\cup\Gamma_h^-\cup\Lambda_h$ and
\begin{equation}\label{eq:Lambda_h}
\begin{aligned}
& \Lambda_h=\mathcal{K}\cap \partial B_h.\\
\end{aligned}
\end{equation}
Clearly in $S_h$, the unit outward normal vectors to $\Gamma_h^+$ and $\Gamma_h^-$ are respectively
\begin{equation}\label{eq:nu}
	\bmf{\nu }\big |_{\Gamma_h^+}=(-\sin\varphi_0,\cos\varphi_0),\quad \bmf{\nu }\big |_{\Gamma_h^-}=(0,-1).
\end{equation}

Throughout the rest of the paper, we set
\begin{equation}\label{eq:kpks}
	k_{{p}}=\sqrt{ \frac{\kappa }{\lambda+2 \mu} } \text { and } k_{{s}}=\sqrt{ \frac{\kappa }{\mu}}, 
\end{equation}
which are known as the compressional and shear wave numbers, respectively. We next present a few lemmas that will be needed in our subsequent analysis. The following lemma from \cite{DR95,SP} states the Fourier expansion in terms of the radial wave functions of the solution $\bmf{u}$ to \eqref{eq:lame} around the origin.
 \begin{lem}\cite{DR95,SP}\label{lem:u exp}
Recall that  $J_m(t)$ is the first-kind Bessel function of order $m\in \mathbb{N} \cup \{0\} $ and $\bmf{x}=r(\cos \varphi, \sin \varphi )^\top$. $\bmf{u}(\bmf{x})$ to \eqref{eq:lame} has the following radial wave expansion at the origin,
\begin{equation}\label{eq:radial}
 \begin{aligned}
 \mathbf{u}(\mathbf{x})
 =& \sum_{m=0}^{\infty}\left\{a_{m}\left\{k_{p} J_{m}^{\prime}\left(k_{p} r\right) \mathrm{e}^{\bsi m \varphi} \mathbf{\hat{r}}+\frac{\bsi m}{r} J_{m}\left(k_{p} r\right) \mathrm{e}^{\bsi m \varphi}  \bm{\hat{\varphi}}\right\}\right.\\ &+b_{m}\left\{\frac{\bsi m}{r} J_{m}\left(k_{s} r\right) \mathrm{e}^{\bsi m \varphi} \mathbf{\hat{r}}-k_{s} J_{m}^{\prime}\left(k_{s} r\right) \mathrm{e}^{\bsi m \varphi} \bm{\hat{\varphi}}\right\}\bigg\},
 \end{aligned}
\end{equation}
where $a_m$ and $b_m$ are constants, $\bm{\hat{\varphi}}=\left(\begin{array}{c}{-\sin{\varphi}}\\ {\cos{\varphi}}
\end{array}\right)$,
$\bmf{\hat{r}}=\left(\begin{array}{c}{\cos{\varphi}}\\ {\sin{\varphi}}
\end{array}\right)$ and the prime denotes the differentiation with respect
to $k_a r$, $a=p$ or $a=s$. Note that \eqref{eq:radial} converges uniformly on compact subsets of $\mathbb R^2$.
\end{lem}

By the analyticity of $\bmf{u}$ in the interior domain of $\Omega$ and the analytic continuation principle, we have the following proposition.
\begin{prop}\label{prop:21}	
Suppose $\mathbf{0}\in \Omega$ and $\bmf{u}$ has the expansion \eqref{eq:radial} around the origin such  that $a_m=b_m=0$ for $\forall m \in \mathbb \cup \{0\}$. Then
$$
\bmf{u} \equiv \bmf{0} \mbox{ in } \Omega .
$$
\end{prop}

The recursive relationship of the first-kind Bessel function and its derivative can be found in \cite{Abr}. In fact we have
\begin{lem}\cite{Abr}\label{lem:J exp}
Recall that $J_m(t)$ is the first-kind Bessel function of the order  $m \in \mathbb{N}\cup \{0\}$. Then
\begin{equation}\label{eq:lem recursive}
	J_m'\left(t\right)=\frac{J_{m-1}\left(t\right)-J_{m+1}\left(t\right)}{2},\quad J_m\left(t\right)=\frac{t\left(J_{m-1}\left(t\right)+J_{m+1}\left(t\right)\right).
}{2m}
\end{equation}
Moreover, we have
\begin{equation}\label{eq:J-1}
	J_{-m}(t)=(-1)^m J_m(t).
\end{equation}
\end{lem}

\begin{rem}
	Using  Lemma \ref{lem:J exp}, one can derive that
\begin{equation}\label{eq:J2}
\begin{aligned}
& \frac{J_{m-1}(k_p r)}{r}=\frac{k_p}{2(m-1)}\left(J_{m-2}(k_p r)+J_m(k_p r)\right), \frac{J_{m+1}(k_p r)}{r}=\frac{k_p}{2(m+1)}\left(J_m(k_p r)+J_{m+2}(k_p r)\right),\\
& \frac{J_{m-1}(k_s r)}{r}=\frac{k_s}{2(m-1)}\left(J_{m-2}(k_s r)+J_m(k_s r)\right), \frac{J_{m+1}(k_s r)}{r}=\frac{k_s}{2(m+1)}\left(J_m(k_s r)+J_{m+2}(k_s r)\right).
\end{aligned}
\end{equation}
\end{rem}

\begin{lem}\label{lem:uu exp}
The radial wave expression of $\mathbf{u}(\bmf{x})$ to \eqref{eq:lame} at the origin can be written as
\begin{equation}\label{eq:u}
 \begin{aligned}
 \mathbf{u}(\mathbf{x})=  \sum_{m=0}^{\infty} & \left\{ \frac{k_p}{2} a_{m} \mathrm{e}^{\bsi m \varphi} \left\{J_{m-1}\left(k_{p} r\right)\mathrm{e}^{-\bsi \varphi}\mathbf{e}_1
 -J_{m+1}\left(k_{p}r\right)\mathrm{e}^{\bsi \varphi}\mathbf{e}_2 \right\}\right.\\
 & + \frac{\bsi k_s}{2} b_{m} \mathrm{e}^{\bsi m \varphi} \left\{J_{m-1}\left(k_{s} r\right)\mathrm{e}^{-\bsi \varphi}\mathbf{e}_1
  +J_{m+1}\left(k_{s} r\right)\mathrm{e}^{\bsi \varphi}
  \mathbf{e}_2
   \right\} \bigg\} .
 \end{aligned}
\end{equation}
where and also throughout the rest of the paper, $\mathbf{e}_1:=(1,\bsi)^\top \mbox{ and }\mathbf{e}_2:=(1,-\bsi)^\top$.
\end{lem}

\begin{proof}
Using \eqref{eq:lem recursive} and Euler's formula, it can be deduced that
\begin{equation}\label{eq:lem 110}
\begin{aligned}
& k_{p} J_{m}^{\prime}\left(k_{p} r\right) \mathrm{e}^{\bsi m \varphi} \mathbf{\hat{r}}+\frac{\bsi m}{r} J_{m}\left(k_{p} r\right) \mathrm{e}^{\bsi m \varphi}  \bm{\hat{\varphi}}\\
= & k_p \frac{J_{m-1}\left(k_p r\right)-J_{m+1}\left(k_p r\right)}{2} \mathrm{e}^{\bsi m \varphi} \mathbf{\hat{r}}+ \frac{\bsi k_p }{2}\left(J_{m-1}\left(k_p r\right)+J_{m+1}\left(k_p r\right)\right)\mathrm{e}^{\bsi m \varphi}  \bm{\hat{\varphi}}\\
= & \frac{k_p}{2} \mathrm{e}^{\bsi m \varphi} \left\{J_{m-1}\left(k_p r\right)\mathrm{e}^{-\bsi\varphi}\mathbf{e}_1- J_{m+1}\left(k_p r\right)\mathrm{e}^{\bsi\varphi}\mathbf{e}_2\right\}.
\end{aligned}
\end{equation}
Similarly, we have
\begin{equation}\label{eq:lem 111}
\begin{aligned}
& \frac{\bsi m}{r} J_{m}\left(k_{s} r\right) \mathrm{e}^{\bsi m \varphi} \mathbf{\hat{r}}-k_{s} J_{m}^{\prime}\left(k_{s} r\right) \mathrm{e}^{\bsi m \varphi}\bm{\hat{\varphi}}  \\
= & \frac{\bsi k_s}{2}\left(J_{m-1}\left(k_s r\right)+J_{m+1}\left(k_s r\right)\right)\mathrm{e}^{\bsi m \varphi}  \mathbf{\hat{r}}-k_s \frac{J_{m-1}\left(k_s r\right)-J_{m+1}\left(k_s r\right)}{2} \mathrm{e}^{\bsi m \varphi} \bm{\hat{\varphi}}\\
= & \frac{\bsi k_s}{2} \mathrm{e}^{\bsi m \varphi} \left\{J_{m-1}\left(k_s r\right)\mathrm{e}^{-\bsi\varphi}\mathbf{e}_1+ J_{m+1}\left(k_s r\right)\mathrm{e}^{\bsi\varphi}\mathbf{e}_2\right\}.
\end{aligned}
\end{equation}
Substituting \eqref{eq:lem 110} and \eqref{eq:lem 111} into \eqref{eq:radial}, after some algebraic calculations, we can prove \eqref{eq:u}.
\end{proof}

\begin{rem}
	In view of \eqref{eq:u}, we have
\begin{equation}\label{eq:u comp}
\begin{aligned}
u_1 \left(\bmf{x}\right)=  \sum_{m=0} ^\infty & \Big [ \frac{k_p}{2} a_m \left(\mathrm{e}^{\bsi \left(m-1\right) \varphi} J_{m-1} \left(k_p r\right) - \mathrm{e}^{\bsi \left(m+1\right) \varphi} J_{m+1} \left(k_p r\right) \right) \\
& + \frac{\bsi k_s}{2} b_m \left(\mathrm{e}^{\bsi \left(m-1\right) \varphi} J_{m-1} \left(k_s r\right) + \mathrm{e}^{\bsi \left(m+1\right) \varphi} J_{m+1} \left(k_s r\right) \right)  \Big ],\\
u_2 \left(\bmf{x}\right)=  \sum_{m=0} ^\infty & \Big [ \frac{\bsi k_p}{2} a_m \left( \mathrm{e}^{\bsi \left(m-1\right) \varphi} J_{m-1} \left(k_p r\right) + \mathrm{e}^{\bsi \left(m+1\right) \varphi} J_{m+1} \left(k_p r\right) \right)  \\
& + \frac{ k_s}{2} b_m \left(-\mathrm{e}^{\bsi \left(m-1\right) \varphi} J_{m-1} \left(k_s r\right) + \mathrm{e}^{\bsi \left(m+1\right) \varphi} J_{m+1} \left(k_s r\right)  \right) \Big ].
\end{aligned}
\end{equation}
\end{rem}

Using \eqref{eq:radial}, we can obtain the corresponding Fourier representation of  the boundary traction operator $T_{\bmf{\nu }}\mathbf u \Big |_{\Gamma^\pm_h }$ defined in \eqref{eq:Tu} as follows. 


\begin{lem}\label{lem:Tuu exp}
Let $\bmf{u}(\bmf{x})$ be a Lam\'e eigenfunction to \eqref{eq:lame} with the Fourier expansion \eqref{eq:u} and $\Gamma^\pm_h$ is defined in \eqref{eq:gamma_pm}. Then $ T_\nu \mathbf{u}\Big |_{\Gamma^+_h}$ possesses the following radial wave expansion at the origin
\begin{equation}\label{eq:Tu1}
 \begin{aligned}
 T_\nu \mathbf{u}\Big |_{\Gamma^+_h} & =  \sum_{m=0}^{\infty}  \left\{\frac{\bsi k_p^2}{2} a_{m} \Big[ \mathrm{e}^{\bsi \left(m-2\right) \varphi} \mathrm{e}^{\bsi \varphi_0} \mu  J_{m-2}(k_p r) \mathbf{e}_1 +  \mathrm{e}^{\bsi m \varphi} \mathrm{e}^{ - \bsi \varphi_0} (\lambda+\mu) J_m(k_p r) \mathbf{e}_1\right.\\
 &  - \mathrm{e}^{\bsi (m+2) \varphi} \mathrm{e}^{-\bsi \varphi_0} \mu J_{m+2}\left(k_{p} r\right) \mathbf{e}_2  - \mathrm{e}^{\bsi m \varphi} \mathrm{e}^{\bsi \varphi_0} \left(\lambda+\mu\right) J_m\left(k_{p} r\right) \mathbf{e}_2\Big] \\
 &- \frac{k_s^2}{2} b_{m}\Big[  \mathrm{e}^{\bsi \left(m-2\right) \varphi} \mathrm{e}^{\bsi \varphi_0} \mu J_{m-2}\left(k_{s} r\right) \mathbf{e}_1 -  \mathrm{e}^{\bsi \left(m+2\right) \varphi} \mathrm{e}^{-\bsi \varphi_0} \mu J_{m+2}\left(k_{s} r\right) \mathbf{e}_2\Big]
 \bigg\}.
 \end{aligned}
\end{equation}
Similarly, the radial wave expansion of $ T_\nu \mathbf{u}\Big |_{\Gamma^-_h}$ at the origin is given by
\begin{equation}\label{eq:Tu2}
 \begin{split}
 & T_\nu \mathbf{u}\Big |_{\Gamma^-_h}=  \sum_{m=0}^{\infty}  \left\{-\frac{\bsi k_p^2}{2} a_{m} \mu  J_{m-2}(k_p r) \mathbf{e}_1 - \frac{\bsi k_p^2}{2} a_{m}   (\lambda+\mu) J_m(k_p r) \mathbf{e}_1\right.\\
 & + \frac{k_s^2}{2} b_{m}  \mu J_{m-2}\left(k_{s} r\right) \mathbf{e}_1 + \frac{\bsi k_p^2}{2} a_{m} \left(\lambda+\mu\right) J_m\left(k_{p} r\right) \mathbf{e}_2 + \frac{\bsi k_p^2}{2} a_{m}   \mu J_{m+2}\left(k_{p} r\right) \mathbf{e}_2 + \frac{k_s^2}{2} b_{m}  \mu J_{m+2}\left(k_{s} r\right) \mathbf{e}_2
 \bigg\}.
 \end{split}
\end{equation}
\end{lem}

The proof of Lemma~\ref{lem:Tuu exp} involves rather tedious calculations and it is postponed to be given in Appendix. Combing Lemmas \ref{lem:uu exp} and \ref{lem:Tuu exp}, for impedance boundary conditions defined on $\Gamma_h^\pm$ with the boundary parameters being constant on  $\Gamma_h^\pm$, we have
\begin{lem}\label{Tu+ exp}
Let  $\Gamma_h^-$ and $\Gamma_h^+$ be two impedance lines of $\bmf{u}$ with constant boundary parameters  $\eta_1$ and $\eta_2$ respectively. We have
\begin{equation}\label{eq:Tu3}
 \begin{aligned}
& \left(T_\nu \mathbf{u}+\eta_2 \boldsymbol{u}\right )\Big |_{\Gamma^+_h}=  \sum_{m=0}^{\infty}  \bigg \{\frac{\bsi k_p^2} {2} a_{m} \Big[ \mathrm{e}^{\bsi \left(m-2\right) \varphi} \mathrm{e}^{\bsi \varphi_0} \mu  J_{m-2}(k_p r) \mathbf{e}_1 +  \mathrm{e}^{\bsi m \varphi} \mathrm{e}^{ - \bsi \varphi_0} (\lambda+\mu) J_m(k_p r) \mathbf{e}_1 \\
 & -\mathrm{e}^{\bsi m \varphi} \mathrm{e}^{\bsi \varphi_0} \left(\lambda+\mu\right) J_m\left(k_{p} r\right) \mathbf{e}_2  -\mathrm{e}^{\bsi (m+2) \varphi} \mathrm{e}^{-\bsi \varphi_0} \mu J_{m+2}\left(k_{p} r\right) \mathbf{e}_2 \Big] \\
 &  - \frac{k_s^2}{2} b_{m}\Big[  \mathrm{e}^{\bsi \left(m-2\right) \varphi} \mathrm{e}^{\bsi \varphi_0} \mu J_{m-2}\left(k_{s} r\right) \mathbf{e}_1+  \mathrm{e}^{\bsi \left(m+2\right) \varphi} \mathrm{e}^{-\bsi \varphi_0} \mu J_{m+2}\left(k_{s} r\right) \mathbf{e}_2\Big] +\frac{\eta_2 k_p}{2} a_{m} \mathrm{e}^{\bsi m \varphi}\times \\ 
& \bigg[J_{m-1}\left(k_{p} r\right)\mathrm{e}^{-\bsi  \varphi}\mathbf{e}_1 -J_{m+1}\left(k_{p}r\right)\mathrm{e}^{\bsi \varphi}\mathbf{e}_2 \bigg] + \frac{\bsi \eta_2 k_s}{2} b_{m} \mathrm{e}^{\bsi m \varphi} \left[J_{m-1}\left(k_{s} r\right)\mathrm{e}^{-\bsi \varphi}\mathbf{e}_1
  +J_{m+1}\left(k_{s} r\right)\mathrm{e}^{\bsi \varphi} \mathbf{e}_2
 \right]
 \bigg\},
 \end{aligned}
\end{equation}
and
 \begin{equation}\label{eq:Tu4}
 \begin{aligned}
& \left( T_\nu \mathbf{u}+\eta_1 \boldsymbol{u}\right) \Big|_{\Gamma_h^-}=\sum_{m=0}^{\infty}  \biggl \{-\frac{\bsi {k}_p^2}{2} a_{m} \Big[ \mu   J_{m-2}(k_p r) \mathbf{e}_1  +   (\lambda+\mu) J_m(k_p r) \mathbf{e}_1\\
& -  \left(\lambda+\mu\right) J_m\left(k_{p} r\right) \mathbf{e}_2 -  \mu J_{m+2}\left(k_{p} r\right) \mathbf{e}_2\Big]   + \frac{k_s^2}{2} b_{m}\Big[   \mu J_{m-2}\left(k_{s} r\right) \mathbf{e}_1 +  \mathrm{e}^{-\bsi \varphi_0} \mu J_{m+2}\left(k_{s} r\right) \mathbf{e}_2 \Big] \\
 & +\frac{\eta_1 k_p}{2} a_{m} \left[J_{m-1}\left(k_{p} r\right)\mathbf{e}_1
 -J_{m+1}\left(k_{p}r\right) \mathbf{e}_2 \right] + \frac{\bsi \eta_1 k_s}{2} b_{m} \left[J_{m-1}\left(k_{s} r\right)\mathbf{e}_1
  +J_{m+1}\left(k_{s} r\right)\mathbf{e}_2
 \right]
 \bigg\}.
 \end{aligned}
\end{equation}

\end{lem}

\begin{lem}{\cite{CDLZ}}\label{lem:co exp}
Suppose that for $0<h \ll 1$ and $t \in (0,h)$, $\sum_{n=0}^{\infty}\alpha_{n}J_{n}(t)=0,$
where $J_{n}\left(t\right)$ is the $n$-th Bessel function of the first kind. Then $\alpha_{n}=0,n=0,1,2,...$
\end{lem}


Next we set
\begin{equation}\label{eq:v}
\bmf{v}(\bmf{x} )=\left(\begin{array}{c}{\exp (-s \sqrt{r} \exp(\bsi \varphi/2))} \\ \bsi \cdot {\operatorname{exp}(-s \sqrt{r} \exp(\bsi \varphi/2) )}\end{array}\right):=
\left(\begin{array}{c}{v_1(\bmf{x})} \\ {v_2(\bmf{x})}\end{array}\right)=v_1(\bmf{x}) \bmf{e}_1,
\end{equation}
where $\bmf{x}=r\cdot (\cos  \varphi, \sin  \varphi ) $,  $s \in \mathbb{R}_{+}$,  $-\pi<\varphi\leqslant \pi$ and $\bmf{e}_1$ is defined in Lemma~\ref{lem:uu exp}.
$\bmf{v}$ is known as the Complex Geometrical Optics (CGO) solution for the Lam\'e operator and it was first introduced in \cite{EBL}. We have

\begin{lem}\cite[Lemma 2.1]{EBL}
	Let  $\Omega \subset {  \mathbb R}^2 $ be such that $\Omega \cap\left(\mathbb{R}_{-} \cup\{\bmf{0}\}\right)=\emptyset$ and $\bmf{v}$ be defined in \eqref{eq:v}. Then there holds $\mathcal{L}  \bmf{v} =0 \mbox{ in } \Omega$.
\end{lem}

%
%

By direct calculations, one can derive the following lemma. 

\begin{lem}\label{eq:lem22 t}
	Let $\bmf{v}$ be defined in \eqref{eq:v}. Denote
	\begin{equation}\label{eq:zeta}
		\zeta{(\varphi)}=-\mathrm{e}^{\bsi \varphi/2 }.
	\end{equation}
	For any given  curve $\Gamma \Subset \mathbb{R}^2$ with a unit normal vector $\nu =(\nu_1,\nu_2 ) $, if $\bmf{v}$ is complex analytic in a neighbourhood of  $\Gamma$, then
	$$
	T_\nu \bmf{v}\big |_{\Gamma} (\bmf{x})= \mu (\nu_1+\bsi \nu_2 ) \frac{s \exp ( sr^{1/2} \zeta{(\varphi)}  )}{r^{1/2} \zeta{(\varphi)}  } \bmf{e}_1, \quad \bmf{x} =r(\cos \varphi, \sin \varphi ) \in \Gamma,
	$$
	where the boundary traction operator  $T_\nu \bmf{v}\big |_{\Gamma} $  and $\bmf{e}_1$ are defined in \eqref{eq:Tu} and Lemma~\ref{lem:uu exp}, respectively.
\end{lem}

Using Lemma \ref{eq:lem22 t} and the fact that the CGO solution $\bmf{v}$ is complex analytic on a  neighborhood  of  $\Gamma_h^\pm \backslash\{ \bmf{0}\} $, and noting that the unit normal vectors on $\Gamma_h^\pm$ are given by \eqref{eq:nu}, we have
\begin{lem}\label{lem:lemTv}
Let 	$\bmf{v}$ be defined in \eqref{eq:v}. Recall that $\Gamma_h^\pm$ are given by \eqref{eq:gamma_pm}, where their corresponding unit normal vectors $\nu\big|_{\Gamma_h^+ }$ and $\nu\big|_{\Gamma_h^- }$ are defined in \eqref{eq:nu}. Then we have
\begin{equation*}
	\begin{split}
		T_{\nu} \bmf{v}\big|_{ \Gamma_h^+ \backslash\{\bmf{0}\} }=\bsi s \mu \zeta{(\varphi_0)}\frac{\exp(sr^{1/2} \zeta{(\varphi_0))}}{r^{1/2}}\mathbf{e}_1,\ \
			T_{\nu} \bmf{v}\big|_{ \Gamma_h^- \backslash\{\bmf{0}\} }=\bsi s \mu \frac{\exp(-sr^{1/2} )}{r^{1/2}} \mathbf{e}_1,\\
	\end{split}
\end{equation*}
where $\zeta(\varphi_0)$ is defined in \eqref{eq:zeta}.
\end{lem}

%

Next, we derive the expansions of $ T_\nu \mathbf{u}\cdot \bmf{v} $, $ T_\nu \mathbf{v}\cdot \bmf{u} $ and $ \mathbf{u}\cdot \bmf{v} $ on $\Gamma_h^\pm$ around the origin, where $\bmf{u}$ is given by \eqref{eq:u} and $\bmf{v}$ is the CGO solution defined in \eqref{eq:v}.  These expansions will be used to analyze the vanishing property of $\bmf{u}$ at the intersecting point $\bmf{0}$ of $\Gamma_h^\pm$. 

\begin{lem}\cite[Proposition 2.1.7]{krantz}\label{lem:kra}
	If the power series $\sum_{\mu } a_{\mu } \bmf{x}^\mu $ converges at a point $\bmf{x}_0$, then it converges uniformly and absolutely on compact subsets of $U(\bmf{x}_0)$, where
	$$
	U(\bmf{x}_0)=\{(r_1 x_{0,1},\ldots, r_n x_{0,n}):-1<r_j<1,j=1,\ldots,n\}, \,    \bmf{x}_0=(x_{0,1},\ldots,  x_{0,n}) \in {\mathbb R}^n.
	$$
\end{lem}

\begin{lem}\label{lem:211 expan}
Let $\bmf{u}$ be a Lam\'e eigenfunction to \eqref{eq:lame} and the CGO solution $\bmf{v}$ be defined in \eqref{eq:v}. Recall that the Lam\'e eigenfunction   $\mathbf{u}$ to \eqref{eq:lame}  has the radial wave expansion \eqref{eq:u} at the origin. Then the following expansions
\begin{equation}\label{eq:Tudotv}
	\begin{split}
		 T_\nu \mathbf{u}\Big |_{\Gamma^+_h} \cdot \bmf{v} \Big |_{\Gamma^+_h} & =
  - {\mathrm e}^{s \sqrt{ r} \zeta(\varphi_0 ) }\bigg\{\bsi k_p^2 (\lambda+\mu) \mathrm e^{\bsi \varphi_0} a_0 + \frac{\bsi}{2}  k_p^3 (\lambda+\mu) \mathrm e^{2 \bsi \varphi_0} a_1 r\\
&\quad +\frac{1}{8}(\bsi k_p^4 (\lambda+\mu) \mathrm e^{3 \bsi \varphi_0} a_2 - \bsi k_p^4 (2 \lambda+\mu) \mathrm e^{ \bsi \varphi_0} a_0+ k_s^4 \mu \mathrm e^{ \bsi \varphi_0} b_0)r^2 + R_{1,\Gamma_h^+} \bigg\},\\
  T_\nu \mathbf{u}\Big |_{\Gamma^-_h} \cdot \bmf{v} \Big |_{\Gamma^-_h} & =
{\mathrm e}^{-s \sqrt{ r}}\bigg\{\bsi k_p^2 (\lambda+\mu)  a_0 + \frac{\bsi}{2}  k_p^3 (\lambda+\mu)  a_1 r\\
&\quad +\frac{1}{8}(\bsi k_p^4 (\lambda+\mu)  a_2 - \bsi k_p^4 (2 \lambda+\mu) a_0+ k_s^4 \mu  b_0)r^2 + R_{1,\Gamma_h^-} \bigg\},
	\end{split}
\end{equation}
converge uniformly and absolutely in $r\in (0,h]$,
where
\begin{equation}\label{eq:RTuv+}
\begin{aligned}
& R_{1,\Gamma_h^+} = r^3 \Bigg\{ \bsi k_p^2 (\lambda+\mu)\left[a_0 \mathrm e^{ \bsi \varphi_0} \sum_{k=2}^{\infty} \frac{(-1)^k k_p^{2k}}{2^{2k}k! k!}r^{2k-3} \right.\\
&+\sum_{m=1}^2\sum_{k=1}^{\infty} a_m \mathrm e^{ \bsi (m+1) \varphi_0} \frac{(-1)^k k_p^{2k+m}}{2^{2k+m}k! (k+m)!}r^{2k+m-3}+\sum_{m=3}^{\infty}\sum_{k=0}^{\infty} a_m \mathrm e^{ \bsi (m+1) \varphi_0} \frac{(-1)^k k_p^{2k+m}}{2^{2k+m}k! (k+m)!}r^{2k+m-3} \bigg]\\
& + \bsi k_p^2 \mu \left[a_0 \mathrm e^{ \bsi \varphi_0} \sum_{k=1}^{\infty}\frac{(-1)^k k_p^{2k+2}}{2^{2k+2}k! (k+2)!}r^{2k-1}\right. + \sum_{m=1}^{\infty}\sum_{k=0}^{\infty} a_m \mathrm e^{ \bsi (m+1) \varphi_0} \frac{(-1)^k k_p^{2k+m+2}}{2^{2k+m+2}k! (k+m+2)!}r^{2k+m-1}\bigg]\\
& + k_s^2 \mu \left[b_0 \mathrm e^{ \bsi \varphi_0} \sum_{k=1}^{\infty}\frac{(-1)^k k_S^{2k+2}}{2^{2k+2}k! (k+2)!}r^{2k-1}\right.+ \sum_{m=1}^{\infty}\sum_{k=0}^{\infty} b_m \mathrm e^{ \bsi (m+1) \varphi_0} \frac{(-1)^k k_s^{2k+m+2}}{2^{2k+m+2}k! (k+m+2)!}r^{2k+m-1}\bigg]\Bigg\},
\end{aligned}
\end{equation}
and
\begin{equation}\label{eq:RTuv-}
\begin{aligned}
& R_{1,\Gamma_h^-} = r^3 \Bigg\{ \bsi k_p^2 (\lambda+\mu)\left[a_0  \sum_{k=2}^{\infty} \frac{(-1)^k k_p^{2k}}{2^{2k}k! k!}r^{2k-3} \right.+\sum_{m=1}^2\sum_{k=1}^{\infty} a_m  \frac{(-1)^k k_p^{2k+m}}{2^{2k+m}k! (k+m)!}r^{2k+m-3}\\
&+\sum_{m=3}^{\infty}\sum_{k=0}^{\infty} a_m  \frac{(-1)^k k_p^{2k+m}}{2^{2k+m}k! (k+m)!}r^{2k+m-3} \bigg]\\
& + \bsi k_p^2 \mu \left[a_0  \sum_{k=1}^{\infty}\frac{(-1)^k k_p^{2k+2}}{2^{2k+2}k! (k+2)!}r^{2k-1}\right.
 + \sum_{m=1}^{\infty}\sum_{k=0}^{\infty} a_m  \frac{(-1)^k k_p^{2k+m+2}}{2^{2k+m+2}k! (k+m+2)!}r^{2k+m-1}\bigg]\\
& + k_s^2 \mu \left    [b_0  \sum_{k=1}^{\infty}\frac{(-1)^k k_s^{2k+2}}{2^{2k+2}k! (k+2)!}r^{2k-1}\right. + \sum_{m=1}^{\infty}\sum_{k=0}^{\infty} b_m  \frac{(-1)^k k_s^{2k+m+2}}{2^{2k+m+2}k! (k+m+2)!}r^{2k+m-1}\bigg]\Bigg\}.
\end{aligned}
\end{equation}
Furthermore, the following expansions
\begin{equation}\label{eq:Tvu}
	\begin{split}
	& T_{\nu} \bmf{v}\big|_{ \Gamma_h^+ \backslash\{\bmf{0}\} } \cdot \bmf{u}\big|_{ \Gamma_h^+ \backslash\{\bmf{0}\} }=	
\bsi s \mu \zeta(\varphi_0 ) {\mathrm e}^{s \sqrt{ r} \zeta(\varphi_0 ) }\bigg\{\frac{1}{2}(-k_p^2 a_0 + \bsi k_s^2 b_0) {\mathrm e}^{\bsi \varphi_0} r^{\frac{1}{2}}+ \frac{1}{8}(-k_p^3 a_1 + \bsi k_s^3 b_1)\\
&\times {\mathrm e}^{2 \bsi \varphi_0} r^{\frac{3}{2}}
  +\frac{1}{48}(-k_p^4 a_2 + \bsi k_s^4 b_2) {\mathrm e}^{3 \bsi \varphi_0} r^{\frac{5}{2}} + \frac{1}{16}(k_p^4 a_0 - \bsi k_s^4 b_0) {\mathrm e}^{\bsi \varphi_0} r^{\frac{5}{2}}  + R_{2,\Gamma_h^+}\bigg\},\\
 & T_{\nu} \bmf{v}\big|_{ \Gamma_h^- \backslash\{\bmf{0}\} } \cdot \bmf{u}\big|_{ \Gamma_h^- \backslash\{\bmf{0}\} }=
\bsi s \mu  {\mathrm e}^{- s \sqrt{ r}}\bigg\{\frac{1}{2}(-k_p^2 a_0 + \bsi k_s^2 b_0) r^{\frac{1}{2}}\\
& + \frac{1}{8}(-k_p^3 a_1 + \bsi k_s^3 b_1)  r^{\frac{3}{2}} +\frac{1}{48}(-k_p^4 a_2 + \bsi k_s^4 b_2)  r^{\frac{5}{2}} + \frac{1}{16}(k_p^4 a_0 - \bsi k_s^4 b_0) r^{\frac{5}{2}} + R_{2,\Gamma_h^-}\bigg\},\\
&  \bmf{u}\big|_{\Gamma_h^+} \cdot \bmf{v}\big|_{\Gamma_h^+} ={\mathrm e}^{s\sqrt{r} \zeta{(\varphi_0} )}  \bigg\{\frac{1}{2} (-k_p^2 a_0+ \bsi k_s^2 b_0) \mathrm{e}^{\bsi \varphi_0} r + R_{0.\Gamma_h^+}\bigg\}, \\
& \bmf{u}\big|_{\Gamma_h^-} \cdot \bmf{v}\big|_{\Gamma_h^-} =   {\mathrm e}^{-s\sqrt{r}}  \bigg\{\frac{1}{2} (-k_p^2 a_0+ \bsi k_s^2 b_0)  r + R_{0.\Gamma_h^+}\bigg\},
	\end{split}
\end{equation}
converge uniformly and absolutely in $r\in (0,h]$, where
\begin{equation}\label{eq:RTvu+}
\begin{aligned}
&R_{2,\Gamma_h^+} =  r^{\frac{7}{2}}\Bigg\{  a_0 \mathrm e^{ \bsi \varphi_0} \sum_{k=2}^{\infty} \frac{(-1)^{k+1} k_p^{2k+2}}{2^{2k+1}k! (k+1)!}r^{2k-3} \\
&+ \sum_{m=1}^{2} \sum_{k=1}^{\infty} a_m \mathrm e^{ \bsi (m+1)\varphi_0} \frac{(-1)^{k+1} k_p^{2k+m+2}}{2^{2k+m+1}k! (k+m+1)!}r^{2k+m-3} \\
& + \sum_{m=3}^{\infty} \sum_{k=0}^{\infty} a_m \mathrm e^{ \bsi (m+1)\varphi_0} \frac{(-1)^{k+1} k_p^{2k+m+2}}{2^{2k+m+1}k! (k+m+1)!}r^{2k+m-3}  +\bsi  b_0 \mathrm e^{ \bsi \varphi_0} \sum_{k=2}^{\infty} \frac{(-1)^k k_s^{2k+2}}{2^{2k+1}k! (k+1)!}r^{2k-3} \\
& + \bsi \sum_{m=1}^{2} \sum_{k=1}^{\infty} b_m \mathrm e^{ \bsi (m+1)\varphi_0} \frac{(-1)^k k_s^{2k+m+2}}{2^{2k+m+1}k! (k+m+1)!}r^{2k+m-3}\\
& + \bsi \sum_{m=3}^{\infty} \sum_{k=0}^{\infty} b_m \mathrm e^{ \bsi (m+1)\varphi_0} \frac{(-1)^k k_s^{2k+m+2}}{2^{2k+m+1}k! (k+m+1)!}r^{2k+m-3}\Bigg\},
\end{aligned}
\end{equation}
and
\begin{equation}\label{eq:RTvu-}
\begin{aligned}
& R_{2,\Gamma_h^-} =  r^{\frac{7}{2}}\Bigg\{  a_0  \sum_{k=2}^{\infty} \frac{(-1)^{k+1} k_p^{2k+2}}{2^{2k+1}k! (k+1)!}r^{2k-3}  + \sum_{m=1}^{2} \sum_{k=1}^{\infty} a_m  \frac{(-1)^{k+1} k_p^{2k+m+2}}{2^{2k+m+1}k! (k+m+1)!}r^{2k+m-3} \\
& + \sum_{m=3}^{\infty} \sum_{k=0}^{\infty} a_m  \frac{(-1)^{k+1} k_p^{2k+m+2}}{2^{2k+m+1}k! (k+m+1)!}r^{2k+m-3} +\bsi  b_0  \sum_{k=2}^{\infty} \frac{(-1)^k k_s^{2k+2}}{2^{2k+1}k! (k+1)!}r^{2k-3} \\
& + \bsi \sum_{m=1}^{2} \sum_{k=1}^{\infty} b_m  \frac{(-1)^k k_s^{2k+m+2}}{2^{2k+m+1}k! (k+m+1)!}r^{2k+m-3} + \bsi \sum_{m=3}^{\infty} \sum_{k=0}^{\infty} b_m  \frac{(-1)^k k_s^{2k+m+2}}{2^{2k+m+1}k! (k+m+1)!}r^{2k+m-3}\Bigg\},
\end{aligned}
\end{equation}
and
\begin{equation}\label{eq:Ruv+}
\begin{aligned}
& R_{0,\Gamma_h^+} = r^2 \bigg\{ a_0 \mathrm e^{ \bsi \varphi_0} \sum_{k=1}^{\infty} \frac{(-1)^{k+1} k_p^{2k+2}}{2^{2k+1}k! (k+1)!}r^{2k-1} + \sum_{m=1}^{\infty} \sum_{k=0}^{\infty} a_m \mathrm e^{ \bsi (m+1)\varphi_0} \frac{(-1)^{k+1} k_p^{2k+m+2}}{2^{2k+m+1}k! (k+m+1)!}r^{2k+m-1} \\
& + \bsi b_0 \mathrm e^{ \bsi \varphi_0} \sum_{k=1}^{\infty} \frac{(-1)^{k} k_s^{2k+2}}{2^{2k+1}k! (k+1)!}r^{2k-1}  + \bsi \sum_{m=1}^{\infty} \sum_{k=0}^{\infty} b_m \mathrm e^{ \bsi (m+1)\varphi_0} \frac{(-1)^{k} k_s^{2k+m+2}}{2^{2k+m+1}k! (k+m+1)!}r^{2k+m-1} \bigg\},
\end{aligned}
\end{equation}
and
\begin{equation}\label{eq:Ruv-}
\begin{aligned}
& R_{0,\Gamma_h^-}= r^2 \bigg\{ a_0  \sum_{k=1}^{\infty} \frac{(-1)^{k+1} k_p^{2k+2}}{2^{2k+1}k! (k+1)!}r^{2k-1} + \sum_{m=1}^{\infty} \sum_{k=0}^{\infty} a_m  \frac{(-1)^{k+1} k_p^{2k+m+2}}{2^{2k+m+1}k! (k+m+1)!}r^{2k+m-1} \\
& + \bsi b_0  \sum_{k=1}^{\infty} \frac{(-1)^{k} k_s^{2k+2}}{2^{2k+1}k! (k+1)!}r^{2k-1}  + \bsi \sum_{m=1}^{\infty} \sum_{k=0}^{\infty} b_m  \frac{(-1)^{k} k_s^{2k+m+2}}{2^{2k+m+1}k! (k+m+1)!}r^{2k+m-1} \bigg\}. 
\end{aligned}
\end{equation}
\end{lem}

\begin{proof} Recall that
 \begin{equation}\label{eq:Jm ex}
	J_{m+1}(t)=\sum_{\ell=0}^{\infty} \frac{(-1)^\ell }{2^{2\ell+m+1} \ell! (m+\ell+1 )!} t^{2\ell+m+1}.
	\end{equation}
Since
\begin{equation}\label{eq:orthogonal}
		\bmf{e}_1\cdot \bmf{e}_1=\bmf{e}_2\cdot \bmf{e}_2=0,\quad \bmf{e_1} \cdot \bmf{e}_2=2,
	\end{equation}
using \eqref{eq:Jm ex}, \eqref{eq:Tu1} and \eqref{eq:v}, we obtain \eqref{eq:Tudotv}. Using Lemma \ref{lem:lemTv} and \eqref{eq:orthogonal} and in view of \eqref{eq:u} and \eqref{eq:v}, we can derive \eqref{eq:Tvu}.

Recall that $\mathcal K$ is defined in \eqref{eq:K}. Since $\bmf{u}$, $ T_\nu \mathbf{u}$, $ \bmf{v} $ and  $ T_\nu \mathbf{v}$
 are analytic in $S_{2h}$, where $S_{2h}=\mathcal K \cap B_{2h}$ and $B_{2h}$ is a disk centered at the origin with the radius $2h$, from Lemma \ref{lem:kra}, we know that \eqref{eq:Tudotv} and \eqref{eq:Tvu} are convergent uniformly and absolutely in $r\in (0,h]$.
 	\end{proof}

\begin{lem}\label{lem:uv}
	Let $\bmf{u}$ and $\bmf{v}$ be respectively given by \eqref{eq:u} and \eqref{eq:v}. Then the following expansion
	\begin{equation}\label{eq:uv old}
 \begin{aligned}
 \mathbf{u} \cdot  \mathbf{v} ={\mathrm e}^{s \sqrt{ r} \zeta(\varphi ) }  \sum_{m=0}^{\infty} & \mathrm{e}^{\bsi (m+1) \varphi}  \left[- k_p a_{m}
 J_{m+1}\left(k_{p}r\right)
  +\bsi k_s b_{m}
  J_{m+1}\left(k_{s} r\right)
  \right ]
 \end{aligned}
\end{equation}
convergences uniformly  in $S_{2h}:=\mathcal{K}\cap B_{2h}
$,
 where $\mathcal K$ is defined in \eqref{eq:K}.  For $0\leqslant r\leqslant h $,  it holds that
 	\begin{equation}\label{eq:uv}
 \begin{aligned}
&\left| \sum_{m=0}^{\infty} \mathrm{e}^{\bsi (m+1) \varphi}  \left[- k_p a_{m}
 J_{m+1}\left(k_{p}r\right)
  +\bsi k_s b_{m}
  J_{m+1}\left(k_{s} r\right)
  \right ]
  \right|  \leqslant   \frac{r\left|k_p ^2 a_0 - \bsi k_s ^2 b_0 \right|}{2} + r^2\cdot S_1, \end{aligned}
\end{equation}
where
\begin{equation*}
	\begin{aligned}
		S_1&= \   \left| \sum_{k=1}^{\infty} \frac{(-1)^k }{2^{2k+1} k! (k+1)!} \left(- k_p ^{2k+2}  a_0 + \bsi k_s ^{2k+2} b_0 \right) h^{2k-1}
\right| \\
&\quad +\sum_{m=1}^{\infty}    \left| \sum_{k=0}^{\infty} \frac{(-1)^k }{2^{2k+m+1} k! (m+k+1 )!} \left( - k_p ^{2k+m+2}  a_m + \bsi k_s ^{2k+m+2} b_m \right) h^{2k+m-1}
\right|  .
	\end{aligned}
\end{equation*}
Furthermore, if $a_0=b_0=\ldots=a_{\ell-1}=b_{\ell-1}=0$, we can conclude that
\begin{equation}\label{eq:uv1}
 \begin{aligned}
\left| \sum_{m=\ell}^{\infty} \mathrm{e}^{\bsi (m+1) \varphi}  \left[- k_p a_{m}
 J_{m+1}\left(k_{p}r\right)
  +\bsi k_s b_{m}
  J_{m+1}\left(k_{s} r\right)
  \right ]
  \right|  & \leqslant   \frac{r^{\ell+1}\left|k_p ^{\ell+2} a_\ell - \bsi k_s ^{\ell+2} b_\ell \right|}{2^{\ell+1}(\ell+1)!} + r^{\ell+2}\cdot S_1(\ell),
   \end{aligned}
\end{equation}
where
\begin{equation*}
	\begin{aligned}
		S_1(\ell)&= \   \left| \sum_{k=1}^{\infty} \frac{(-1)^k }{2^{2k+\ell+1} k! (k+\ell+1)!} \left( - k_p ^{2k+\ell+2}  a_\ell + \bsi k_s ^{2k+\ell+2} b_\ell \right) h^{2k-1}
\right| \\
&\quad +\sum_{m=\ell+1}^{\infty}    \left| \sum_{k=0}^{\infty} \frac{(-1)^k }{2^{2k+m+1} k! (m+k+1 )!} \left( - k_p ^{2k+m+2}  a_m + \bsi k_s ^{2k+m+2} b_m \right) h^{2k+m-\ell-1}
\right|  .
	\end{aligned}
\end{equation*}
\end{lem}

\begin{proof}
	Since $\bmf{v}$ defined in \eqref{eq:v}  is analytic in $S_{2h} $ and $\bmf{u}$ has the expansion \eqref{eq:u}, by noting 	\eqref{eq:orthogonal},
	we have \eqref{eq:uv old}. Since $\Re(\zeta(\varphi ) )<0$ if $\varphi\in [0,\varphi_0 ]$, when $s$ is sufficient large we have ${\mathrm e}^{s \sqrt{ r} \zeta(\varphi ) } \leqslant 1$. Since \eqref{eq:uv old} is convergent at $\bmf{x}_0 \in \partial B_{2h} \cap {\mathcal K}$ from Lemma \ref{lem:kra} and noting \eqref{eq:Jm ex},
we can obtain \eqref{eq:uv}.
\end{proof}

\begin{lem}\label{lem:Tvu}
	Let $\bmf{u}$ be given in \eqref{eq:u} and $T_\nu \bmf{v}$ be defined in Lemma \ref{lem:lemTv}. If $a_0=b_0=\cdots=a_{\ell-1}=b_{\ell-1}=0$, then the following expansion{
\begin{equation}\label{eq:Tvu guina1new}
	\begin{split}
	& T_{\nu} \bmf{v}\big|_{ \Gamma_h^+ \backslash\{\bmf{0}\} } \cdot \bmf{u}\big|_{ \Gamma_h^+ \backslash\{\bmf{0}\} }=	\bsi s \mu \zeta(\varphi_0 ) {\mathrm e}^{s \sqrt{ r} \zeta(\varphi_0 ) }\bigg\{\frac{{\mathrm e}^{\bsi (\ell+1) \varphi_0}}{2^{\ell+1}(\ell+1)!}(-k_p^{\ell+2} a_\ell + \bsi k_s^{\ell+2} b_\ell)  r^{\ell+\frac{1}{2}}\\
&+ \frac{{\mathrm e}^{\bsi (\ell+2) \varphi_0}}{2^{\ell+2}(\ell+2)!}(-k_p^{\ell+3} a_{\ell+1} + \bsi k_s^{\ell+3} b_{\ell+1})  r^{\ell+\frac{3}{2}} + \frac{{\mathrm e}^{\bsi (\ell+3) \varphi_0}}{2^{\ell+3}(\ell+3)!}(-k_p^{\ell+4} a_{\ell+2} + \bsi k_s^{\ell+4} b_{\ell+2})  r^{\ell+\frac{5}{2}}\\
&+ \frac{{\mathrm e}^{\bsi (\ell+) \varphi_0}}{2^{\ell+3}(\ell+2)!}(k_p^{\ell+4} a_{\ell} - \bsi k_s^{\ell+4} b_{\ell})  r^{\ell+\frac{5}{2}}  + \hat{R}_{2,\Gamma_h^+}\bigg\},\\
 & T_{\nu} \bmf{v}\big|_{ \Gamma_h^- \backslash\{\bmf{0}\} } \cdot \bmf{u}\big|_{ \Gamma_h^- \backslash\{\bmf{0}\} }=	\bsi s \mu {\mathrm e}^{- s \sqrt{ r}  }\bigg\{\frac{1}{2^{\ell+1}(\ell+1)!}(-k_p^{\ell+2} a_\ell + \bsi k_s^{\ell+2} b_\ell)  r^{\ell+\frac{1}{2}}\\
& + \frac{1}{2^{\ell+2}(\ell+2)!}(-k_p^{\ell+3} a_{\ell+1} + \bsi k_s^{\ell+3} b_{\ell+1})  r^{\ell+\frac{3}{2}} + \frac{1}{2^{\ell+3}(\ell+3)!}(-k_p^{\ell+4} a_{\ell+2} + \bsi k_s^{\ell+4} b_{\ell+2})  r^{\ell+\frac{5}{2}}\\
& + \frac{1}{2^{\ell+3}(\ell+2)!}(k_p^{\ell+4} a_{\ell} - \bsi k_s^{\ell+4} b_{\ell})  r^{\ell+\frac{5}{2}}  + \hat{R}_{2,\Gamma_h^-}\bigg\},\\
\end{split}
\end{equation}
}
converge uniformly and absolutely with respect to $r\in (0,h]$, where
{ \begin{equation}\label{eq:RTvu++}
\begin{aligned}
& \hat{R}_{2,\Gamma_h^+} =  r^{\ell+\frac{7}{2}}\Bigg\{  a_\ell \mathrm e^{ \bsi (\ell+1) \varphi_0} \sum_{k=2}^{\infty} \frac{(-1)^{k+1} k_p^{2k+\ell+2}}{2^{2k+\ell+1}k! (k+\ell+1)!}r^{2k-3} + \sum_{m=\ell+1}^{\ell+2} \sum_{k=1}^{\infty} a_m \mathrm e^{ \bsi (m+1)\varphi_0}\times \\
& \frac{(-1)^{k+1} k_p^{2k+m+2}}{2^{2k+m+1}k! (k+m+1)!}r^{2k+m-\ell-3} + \sum_{m=\ell+3}^{\infty} \sum_{k=0}^{\infty} a_m \mathrm e^{ \bsi (m+1)\varphi_0} \frac{(-1)^{k+1} k_p^{2k+m+2}}{2^{2k+m+1}k! (k+m+1)!}r^{2k+m-\ell-3} \\
& +\bsi  b_0 \mathrm e^{ \bsi(\ell+1) \varphi_0} \sum_{k=2}^{\infty} \frac{(-1)^k k_s^{2k+\ell+2}}{2^{2k+\ell+1}k! (k+\ell+1)!}r^{2k-3} + \bsi \sum_{m=\ell+1}^{\ell+2} \sum_{k=1}^{\infty} b_m \mathrm e^{ \bsi (m+1)\varphi_0}\times\\
& \frac{(-1)^k k_s^{2k+m+2}}{2^{2k+m+1}k! (k+m+1)!}r^{2k+m-3} + \bsi \sum_{m=\ell+3}^{\infty} \sum_{k=0}^{\infty} b_m \mathrm e^{ \bsi (m+1)\varphi_0} \frac{(-1)^k k_s^{2k+m+2}}{2^{2k+m+1}k! (k+m+1)!}r^{2k+m-\ell-3}\Bigg\},
\end{aligned}
\end{equation}
}
and
{
\begin{equation}\label{eq:RTvu--}
\begin{aligned}
& \hat{R}_{2,\Gamma_h^-} =  r^{\ell+\frac{7}{2}}\Bigg\{  a_0  \sum_{k=2}^{\infty} \frac{(-1)^{k+1} k_p^{2k+2}}{2^{2k+1}k! (k+1)!}r^{2k-3} + \sum_{m=\ell+1}^{\ell+2} \sum_{k=1}^{\infty} a_m  \frac{(-1)^{k+1} k_p^{2k+m+2}}{2^{2k+m+1}k! (k+m+1)!}r^{2k+m-3} \\
& + \sum_{m=\ell+3}^{\infty} \sum_{k=0}^{\infty} a_m  \frac{(-1)^{k+1} k_p^{2k+m+2}}{2^{2k+m+1}k! (k+m+1)!}r^{2k+m-3} +\bsi  b_0  \sum_{k=2}^{\infty} \frac{(-1)^k k_s^{2k+2}}{2^{2k+1}k! (k+1)!}r^{2k-3} \\
& + \bsi \sum_{m=\ell+1}^{\ell+2} \sum_{k=1}^{\infty} b_m  \frac{(-1)^k k_s^{2k+m+2}}{2^{2k+m+1}k! (k+m+1)!}r^{2k+m-3} + \bsi \sum_{m=\ell+3}^{\infty} \sum_{k=0}^{\infty} b_m  \frac{(-1)^k k_s^{2k+m+2}}{2^{2k+m+1}k! (k+m+1)!}r^{2k+m-3}\Bigg\}.
\end{aligned}
\end{equation}
}
Furthermore, we have the following expansions 
\begin{equation}\label{eq:Ruv guina}
\begin{aligned}
\bmf{u}\big|_{\Gamma_h^+} \cdot \bmf{v}\big|_{\Gamma_h^+} &={\mathrm e}^{s\sqrt{r} \zeta{(\varphi_0} )}  \bigg\{\frac{{\mathrm e}^{\bsi (\ell+1) \varphi_0 }}{2^(\ell+1) (\ell+1)!} (-k_p^{\ell+2} a_\ell+ \bsi k_s^{\ell+2} b_\ell)  r^{\ell+1} + \hat{R}_{0.\Gamma_h^+}\bigg\}, \\
 \bmf{u}\big|_{\Gamma_h^-} \cdot \bmf{v}\big|_{\Gamma_h^-} &=  {\mathrm e}^{-s\sqrt{r}}  \bigg\{\frac{1}{2^(\ell+1) (\ell+1)!} (-k_p^{\ell+2} a_\ell+ \bsi k_s^{\ell+2} b_\ell)  r^{\ell+1} + \hat{R}_{0.\Gamma_h^-}\bigg\}, \\
 \end{aligned}
\end{equation}
which converge uniformly and absolutely in $r\in (0,h]$, where
\begin{equation}\label{eq:Ruv++}
\begin{aligned}
& \hat{R}_{0,\Gamma_h^+} = r^{\ell+2} \bigg\{ a_\ell \mathrm e^{ \bsi (\ell+1)\varphi_0} \sum_{k=1}^{\infty} \frac{(-1)^{k+1} k_p^{2k+\ell+2}}{2^{2k+\ell+1}k! (k+\ell+1)!}r^{2k-1} + \sum_{m=\ell+1}^{\infty} \sum_{k=0}^{\infty} a_m \mathrm e^{ \bsi (m+1)\varphi_0}\times\\
 & \frac{(-1)^{k+1} k_p^{2k+m+2}}{2^{2k+m+1}k! (k+m+1)!}r^{2k+m-\ell-1} + \bsi b_0 \mathrm e^{ \bsi (\ell+1)\varphi_0} \sum_{k=1}^{\infty} \frac{(-1)^{k} k_s^{2k+\ell+2}}{2^{2k+\ell+1}k! (k+\ell+1)!}r^{2k-1} \\
& + \bsi \sum_{m=\ell+1}^{\infty} \sum_{k=0}^{\infty} b_m \mathrm e^{ \bsi (m+1)\varphi_0} \frac{(-1)^{k} k_s^{2k+m+2}}{2^{2k+m+1}k! (k+m+1)!}r^{2k+m-\ell-1} \bigg\},
\end{aligned}
\end{equation}
and
\begin{equation}\label{eq:Ruv--}
\begin{aligned}
& \hat{R}_{0,\Gamma_h^-} =  r^{\ell+2} \bigg\{ a_\ell  \sum_{k=1}^{\infty} \frac{(-1)^{k+1} k_p^{2k+\ell+2}}{2^{2k+\ell+1}k! (k+\ell+1)!}r^{2k-1}  + \sum_{m=\ell+1}^{\infty} \sum_{k=0}^{\infty} a_m  \frac{(-1)^{k+1} k_p^{2k+m+2}}{2^{2k+m+1}k! (k+m+1)!}r^{2k+m-\ell-1}\\
& + \bsi b_0  \sum_{k=1}^{\infty} \frac{(-1)^{k} k_s^{2k+\ell+2}}{2^{2k+\ell+1}k! (k+\ell+1)!}r^{2k-1}  + \bsi \sum_{m=\ell+1}^{\infty} \sum_{k=0}^{\infty} b_m  \frac{(-1)^{k} k_s^{2k+m+2}}{2^{2k+m+1}k! (k+m+1)!}r^{2k+m-\ell-1} \bigg\}.
\end{aligned}
\end{equation}

\end{lem}

\begin{lem}
Recall that $R_{1,\Gamma_h^+}$, $R_{1,\Gamma_h^-} $, $R_{2,\Gamma_h^+} $ , $R_{2,\Gamma_h^-} $ ,$R_{0,\Gamma_h^+}$ , $ R_{0,\Gamma_h^-}$, $\hat{R}_{2,\Gamma_h^+}$ , $ \hat{R}_{2,\Gamma_h^-}$  $ \hat{R}_{0,\Gamma_h^+}$ and $ \hat{R}_{0,\Gamma_h^-}$ are defined in \eqref{eq:RTuv+}, \eqref{eq:RTuv-}, \eqref{eq:RTvu+}, \eqref{eq:RTvu-}, \eqref{eq:Ruv+}, \eqref{eq:Ruv-}, \eqref{eq:RTvu++}, \eqref{eq:RTvu--}, \eqref{eq:Ruv++} and \eqref{eq:Ruv--}  respectively. Then we have
\begin{subequations}
\begin{align}
		&\left |R_{1,\Gamma_h^+}\right| \leq r^3 S_{2},\quad \left |R_{1,\Gamma_h^-}\right| \leq r^3 S_{2},\label{eq:R1 s} \\
				&\left |R_{2,\Gamma_h^+}\right| \leq   r^{7/2} S_{3},\quad \left |R_{2,\Gamma_h^-}\right| \leq  r^{7/2} S_{3},\label{eq:R2 s}\\
         &\left |R_{0,\Gamma_h^+}\right| \leq   r^2 S_{0},\quad \left |R_{0,\Gamma_h^-}\right| \leq  r^2 S_{0},\label{eq:R3 s}\\
         &\left |\hat{R}_{2,\Gamma_h^+}\right| \leq   r^{\ell+7/2} \hat{S}_{3},\quad \left |\hat{R}_{2,\Gamma_h^-}\right| \leq  r^{\ell+7/2} \hat{S}_{3},\label{eq:R4 s}\\
         &\left |\hat{R}_{0,\Gamma_h^+}\right| \leq   r^{\ell+2} \hat{S}_{0},\quad \left |\hat{R}_{0,\Gamma_h^-}\right| \leq  r^{\ell+2} \hat{S}_{0},\label{eq:R5 s}
\end{align}
\end{subequations}
where
\begin{equation}\label{eq:RTuv+S}
\begin{aligned}
S_2 =&  k_p^2 (\lambda+\mu)\bigg [|a_0|  \sum_{k=2}^{\infty} \frac{ k_p^{2k}}{2^{2k}k! k!}h ^{2k-3} +\sum_{m=1}^2\sum_{k=1}^{\infty} |a_m|  \frac{ k_p^{2k+m}}{2^{2k+m}k! (k+m)!} h^{2k+m-3}\bigg] \\
& +  k_p^2 \mu  |a_0|   \sum_{k=1}^{\infty}\frac{  k_p^{2k+2}}{2^{2k+2}k! (k+2)!}h^{2k-1} +k_s^2 \mu |b_0|  \sum_{k=1}^{\infty}\frac{ k_s^{2k+2}}{2^{2k+2}k! (k+2)!} h ^{2k-1}
 \\
& +  k_p^2 \mu  \sum_{m=1}^{2}\sum_{k=0}^{\infty} |a_m|   \frac{ k_p^{2k+m+2}}{2^{2k+m+2}k! (k+m+2)!}h^{2k+m-1} \\
& + k_s^2 \mu  \sum_{m=1}^{2}\sum_{k=0}^{\infty} |b_m|  \frac{ k_s^{2k+m+2}}{2^{2k+m+2}k! (k+m+2)!}h^{2k+m-1}\bigg]\\
&+\sum_{m=3}^{\infty}\bigg|\bsi k_p^2 (\lambda+\mu) a_m \sum_{k=0}^{\infty}  \frac{(-1)^k k_p^{2k+m}}{2^{2k+m}k! (k+m)!}h ^{2k+m-3}\\
&+ \bsi k_p^2 a_m\mu   \sum_{k=0}^\infty \frac{(-1)^k k_p^{2k+m+2}}{2^{2k+m+2}k! (k+m+2)!}h^{2k+m-1} \\
&+  k_s^2 b_m \mu \sum_{k=0}^\infty \frac{(-1)^k k_s^{2k+m+2}}{2^{2k+m+2}k! (k+m+2)!}h^{2k+m-1} \bigg|,
\end{aligned}
\end{equation}
and
\begin{equation}\label{eq:RTvu+S}
\begin{aligned}
S_3 = & |a_0|   \sum_{k=2}^{\infty} \frac{ k_p^{2k+2}}{2^{2k+1}k! (k+1)!}h^{2k-3} +| b_0 |  \sum_{k=2}^{\infty} \frac{ k_s^{2k+2}}{2^{2k+1}k! (k+1)!}h^{2k-3} \\
& + \sum_{m=1}^{2}\sum_{k=1}^{\infty} |a_m|   \frac{ k_p^{2k+m+2}}{2^{2k+m+1}k! (k+m+1)!}h^{2k+m-3} \\
& + \sum_{m=1}^{2}\sum_{k=1}^{\infty} |b_m|   \frac{ k_s^{2k+m+2}}{2^{2k+m+1}k! (k+m+1)!}h^{2k+m-3} \\
& \quad + \sum_{m=3}^{\infty} \bigg| a_m \sum_{k=0}^{\infty}  \frac{(-1)^{k} k_p^{2k+m+2}}{2^{2k+m+1}k! (k+m+1)!}h^{2k+m-3} \\
&\quad  + \bsi b_m  \sum_{k=0}^{\infty} \frac{(-1)^k k_s^{2k+m+2}}{2^{2k+m+1}k! (k+m+1)!}h^{2k+m-3}\bigg |,
\end{aligned}
\end{equation}

and
\begin{equation}\label{eq:uvS}
\begin{aligned}
S_0 = & |a_0|   \sum_{k=1}^{\infty} \frac{ k_p^{2k+2}}{2^{2k+1}k! (k+1)!}h^{2k-1} +| b_0 |  \sum_{k=1}^{\infty} \frac{ k_s^{2k+2}}{2^{2k+1}k! (k+1)!}h^{2k-1} \\
&  + \sum_{m=1}^{\infty} \bigg| a_m \sum_{k=0}^{\infty}  \frac{(-1)^{k+1} k_p^{2k+m+2}}{2^{2k+m+1}k! (k+m+1)!}h^{2k+m-1}\\
&  + \bsi b_m \sum_{k=0}^{\infty}  \frac{(-1)^{k} k_s^{2k+m+2}}{2^{2k+m+1}k! (k+m+1)!}h^{2k+m-1} \bigg |,
\end{aligned}
\end{equation}
and
\begin{equation}\label{eq:RTvu+Shat}
\begin{aligned}
\hat{S}_3 = & |a_\ell|   \sum_{k=2}^{\infty} \frac{ k_p^{2k+\ell+2}}{2^{2k+\ell+1}k! (k+\ell+1)!}h^{2k-3} +| b_\ell |  \sum_{k=2}^{\infty} \frac{ k_s^{2k+\ell+2}}{2^{2k+\ell+1}k! (k+\ell+1)!}h^{2k-3} \\
& + \sum_{m=\ell+1}^{\ell+2}\sum_{k=1}^{\infty} |a_m|   \frac{ k_p^{2k+m+2}}{2^{2k+m+1}k! (k+m+1)!}h^{2k+m-\ell-3} \\
& + \sum_{m=\ell+1}^{\ell+2}\sum_{k=1}^{\infty} |b_m|   \frac{ k_s^{2k+m+2}}{2^{2k+m+1}k! (k+m+1)!}h^{2k+m-\ell-3} \\
& \quad + \sum_{m=\ell+3}^{\infty} \bigg| a_m \sum_{k=0}^{\infty}  \frac{(-1)^{k} k_p^{2k+m+2}}{2^{2k+m+1}k! (k+m+1)!}h^{2k+m-\ell-3} \\
&\quad  + \bsi b_m  \sum_{k=0}^{\infty} \frac{(-1)^k k_s^{2k+m+2}}{2^{2k+m+1}k! (k+m+1)!}h^{2k+m-\ell-3}\bigg |,
\end{aligned}
\end{equation}

\begin{equation}\label{eq:uvShat}
\begin{aligned}
\hat{S}_0 = & |a_\ell|   \sum_{k=1}^{\infty} \frac{ k_p^{2k+\ell+2}}{2^{2k+\ell+1}k! (k+\ell+1)!}h^{2k-1} +| b_\ell |  \sum_{k=1}^{\infty} \frac{ k_s^{2k+\ell+2}}{2^{2k+\ell+1}k! (k+\ell+1)!}h^{2k-1} \\
&  + \sum_{m=\ell+1}^{\infty} \bigg| a_m \sum_{k=0}^{\infty}  \frac{(-1)^{k+1} k_p^{2k+m+2}}{2^{2k+m+1}k! (k+m+1)!}h^{2k+m-\ell-1}\\
&  + \bsi b_m \sum_{k=0}^{\infty}  \frac{(-1)^{k} k_s^{2k+m+2}}{2^{2k+m+1}k! (k+m+1)!}h^{2k+m-\ell-1} \bigg |,
\end{aligned}
\end{equation}
\end{lem}

\begin{proof}
	Since $ T_\nu \mathbf{u}\Big |_{\Gamma^\pm_h} \cdot \bmf{v} \Big |_{\Gamma^\pm_h} $ are analytic on  ${\Gamma^\pm_h}$, from \eqref{eq:Tudotv}, by root test and Lemma \ref{lem:kra}, for $r\in (0,h)$, together with straightforward though tedious calculations, one can prove \eqref{eq:R1 s}. \eqref{eq:R2 s}, \eqref{eq:R3 s}, \eqref{eq:R4 s} and \eqref{eq:R5 s} can be proved in a similar way, and we skip the details.
\end{proof}

Recall that the open sector  $\mathcal K$ and $\Gamma_h^\pm$ are defined in \eqref{eq:K} and \eqref{eq:gamma_pm}, respectively. For $\varepsilon \in {\mathbb R}_{+}$ satisfying $\varepsilon<h$, let
\begin{equation}\label{eq:Seps}
	S_\varepsilon={\mathcal K} \cap B_\varepsilon, \quad \Gamma_{(0,\varepsilon )}^\pm  =\Gamma_h^\pm \cap B_\varepsilon,\quad \Lambda_\varepsilon=
	S_\varepsilon \cap \partial B_\varepsilon.
\end{equation}

\begin{lem}\label{lem:25 int 0}
	Let $\bmf{u}$  be a Lam\'e eigenfunction to \eqref{eq:lame} and $\bmf{v}$ be defined in \eqref{eq:v}. Recall that $\Gamma_{(0,\varepsilon )}^\pm$ and $\Lambda_\varepsilon$  are defined in \eqref{eq:Seps}. Then

\begin{subequations}
	\begin{align}
		&\lim_{\varepsilon\rightarrow 0^+}	\int_{\Gamma_{(0,\varepsilon )}^+} T_\nu \bmf{v} \cdot \bmf{u} \rmd \sigma=	\lim_{\varepsilon\rightarrow 0^+}	\int_{\Gamma_{(0,\varepsilon )}^-} T_\nu \bmf{v} \cdot \bmf{u} \rmd \sigma=0,\label{eq:27a} \\
		&\lim_{\varepsilon\rightarrow 0^+}	\int_{\Gamma_{(0,\varepsilon )}^+} T_\nu \bmf{u} \cdot \bmf{v} \rmd \sigma=	\lim_{\varepsilon\rightarrow 0^+}	\int_{\Gamma_{(0,\varepsilon )}^-} T_\nu \bmf{u} \cdot \bmf{v} \rmd \sigma=0,\label{eq:27b} \\
		&\lim_{\varepsilon\rightarrow 0^+}	\int_{\Gamma_{(0,\varepsilon )}^+}  \bmf{u} \cdot \bmf{v} \rmd \sigma=	\lim_{\varepsilon\rightarrow 0^+}	\int_{\Gamma_{(0,\varepsilon )}^-}  \bmf{u} \cdot \bmf{v} \rmd \sigma=0,\label{eq:27c}
	\end{align}
\end{subequations}

\end{lem}

\begin{proof}
	From \eqref{eq:Tvu}, it is easy to see that
	$$
	\lim_{\bmf{x} \rightarrow \bmf{0}\atop  \bmf{x} \in \Gamma_h^\pm } T_{\nu} \bmf{v}\big|_{ \Gamma_h^\pm \backslash\{\bmf{0}\} } \cdot \bmf{u}\big|_{ \Gamma_h^\pm  \backslash\{\bmf{0}\} }=0.
	$$
Therefore the function $ T_{\nu} \bmf{v}\big|_{ \Gamma_h^\pm \backslash\{\bmf{0}\} } \cdot \bmf{u}\big|_{ \Gamma_h^\pm  \backslash\{\bmf{0}\} }$ is continuous at the origin. Hence by the dominant convergent theorem, we can prove \eqref{eq:27a}. Similarly, from \eqref{eq:Tudotv} and \eqref{eq:Tvu}, we know that \eqref{eq:27b} and \eqref{eq:27c} hold via the dominant convergent theorem.
\end{proof}

The following lemma gives the estimates of the integrals with respect to the CGO solution $\bmf{v}$ \eqref{eq:v} on an open sector and an arbitrary arc, which will be used in the subsequent study.
\begin{lem}\cite[Proposition 3.1]{EBL}\label{lem:CGO exp}
Let $\bmf{v} : \mathbb{R}^{2} \rightarrow \mathbb{C}^{2}$ be defined by \eqref{eq:v} and
\begin{equation}\label{eq:K phi}
	\mathcal{K}_{\varphi_m,\varphi_M} =\left\{\bmf{x} =(x_1,x_2) \in \mathbb{R}^{2} ~|~\bmf{ x}  \neq \bmf{0} , \quad \varphi_{m}<\arg \left(x_{1}+\mathrm{i} x_{2}\right)<\varphi_{M}\right\}
\end{equation}
for given angles $-\pi<\varphi_{m}<\varphi_{M}<\pi$. Then it holds that
\begin{equation}\label{eq:int 29 n}
	\int_{{\mathcal K}_{\varphi_m,\varphi_M }}  v_1 \left(\bmf{x}\right) \rmd \bmf{x}  = 6 \bsi \left(\mathrm{e}^{-2\varphi_M \bsi}-\mathrm{e}^{-2\varphi_m \bsi}\right)s^{-4}.
\end{equation}
In addition for $\alpha$,$h>0$ and $j\in\left\{1,2\right\}$, we have the upper bounds
\begin{subequations}
\begin{align}
\int_{{\mathcal K}_{\varphi_m,\varphi_M }} \left|v_j \left(\bmf{x}\right)\right| \left|\bmf{x}\right|^\alpha \rmd \bmf{x}& \leqslant \frac{2\left(\varphi_M-\varphi_m\right)\Gamma\left(2\alpha+4\right)}{\delta_{{\mathcal K}_{\varphi_m,\varphi_M }}^{2\alpha+4}}s^{-2\alpha-4},\label{eq:volume n} \\
\int_{{\mathcal K}_{\varphi_m,\varphi_M } \backslash B_h}\left|v_{j}(\bmf{x})\right| \mathrm{d} \bmf{x} &\leqslant \frac{6\left(\varphi_{M}-\varphi _{m}\right)}{\delta_{\mathcal{K}_{\varphi_m, \varphi_M }}^{4}} s^{-4} \mathrm{e}^{-\delta_{\mathcal{K}_{\varphi_m, \varphi_M } } s \sqrt{h} / 2},\label{eq:volume2 n}
\end{align}
\end{subequations}
where $\delta_{{\mathcal K}_{\varphi_m,\varphi_M }}=\min_{\varphi_m<\varphi<\varphi_M} \cos(\varphi/2)$ is a positive constant depending ${\mathcal K}_{\varphi_m,\varphi_M  }$.   Furthermore, assume that $\bmf{u} \in H^{2}\left(\mathcal{K}_{\varphi_m, \varphi_M } \cap B_l \right)$ where $B_l=B(\bmf{0}, l)$ for some $l>0$. The it holds that
\begin{equation}\label{eq:arc n}
\begin{aligned}
 & \left|\int_{\mathcal{K}_{\varphi_m,\varphi_M} \cap \partial B_l }\left[\left({T}_{\nu} \bmf{u}\right) \cdot \bmf{v}-\left({T}_{\nu} \bmf{v}\right) \cdot \bmf{u}\right] \mathrm{d} \sigma\right| \\
 & \leq \mathcal{C}_{{\mathcal K}_{{\varphi_m,\varphi_M }},B_l,\mu,\lambda}\|\bmf{u}\|_{H^2\left( {{\mathcal K}_{\varphi_m,\varphi_M }} \cap B_l \right)}\left(1+s\right)\exp\left(-\delta_{{\mathcal K}_{\varphi_m,\varphi_M }} s \sqrt{h}\right),
\end{aligned}
\end{equation}
where $\mathcal{C}_{{\mathcal K}_{{\varphi_m,\varphi_M }},B_l,\mu,\lambda}$ is a  positive constant depending on $B_l, \lambda, \mu $ and ${\mathcal K}_{\varphi_m,\varphi_M  }$.
\end{lem}

\begin{rem}
	Recall that $\mathcal K$ is defined in \eqref{eq:K}. Then ${\mathcal K}_{\varphi_m, \varphi_M }$ defined in  \eqref{eq:K phi} degenerates to $\mathcal K$ whenever  $\varphi_M:=\varphi_0$, $\varphi_m:=0$. In this situation, the constant $\delta_{{\mathcal K}_{\varphi_m,\varphi_M }}$ given in \eqref{eq:volume n}  and \eqref{eq:volume2 n} is denoted by $\delta_{{\mathcal K}}$ in the remainder of this paper. Indeed, setting  $\varphi_M:=\varphi_0$ and $\varphi_m:=0$, from \eqref{eq:int 29 n},  \eqref{eq:volume n},  \eqref{eq:volume2 n} and \eqref{eq:arc n},  we have
	\begin{subequations}
\begin{align}
&\int_{{\mathcal K} }  {\mathrm e}^{s \sqrt{ r} \zeta(\varphi ) } r  
\rmd {r} \rmd \varphi  = 6 \bsi \left(\mathrm{e}^{-2\varphi_0 \bsi}-1\right)s^{-4},\\
&\int_{{\mathcal K} } {\mathrm e}^{s \sqrt{ r} \Re(\zeta(\varphi )) }  r^{\alpha+1}  \rmd r \rmd \varphi
\leqslant
\frac{2\varphi_0\Gamma\left(2\alpha+4\right)}{\delta_{{\mathcal K}}^{2\alpha+4}}s^{-2\alpha-4},\label{eq:volume} \\
 & \left|\int_{\Lambda_h }\left[\left({T}_{\nu} \bmf{u}\right) \cdot \bmf{v}-\left({T}_{\nu} \bmf{v}\right) \cdot \bmf{u}\right] \mathrm{d} \sigma\right|  \leqslant \mathcal{C}_{{\mathcal K },B_h,\mu,\lambda}\|\bmf{u}\|_{H^2\left( {{\mathcal K}} \cap B_h \right)}\left(1+s\right){\mathrm e}^{-\delta_{{\mathcal K} } s \sqrt{h}} , \label{eq:arc}
\end{align}
\end{subequations}
where $\delta_{{\mathcal K}}=\min_{0<\varphi<\varphi_0} \cos(\varphi/2)$ and  $\Lambda_h$ is defined in \eqref{eq:Lambda_h}.
\end{rem}

We next derive several crucial integral identities.


\begin{lem}\label{lem:part exp}
Suppose $D$ is a bounded Lipschitz domain in $\mathbb R^2$ and $\bmf{u}$, $\bmf{v}$ are $H_{loc}^2(\mathbb{R}^2)$ functions. Let $\bmf{v}$ be the CGO solution defined in \eqref{eq:v}, which satisfies   $\mathcal{L}\bmf v=\bmf{0}$.  Then
%
\begin{equation}\label{eq:CGO1}
\int_D \left(\mathcal{L} \bmf{u}\right) \cdot \bmf{v} \rmd x=\int_{\partial D} \left[\left({T}_{\nu} \bmf{u}\right) \cdot \bmf{v}-\left({T}_{\nu} \bmf{v}\right) \cdot \bmf{u}\right] \mathrm{d} \sigma.
\end{equation}
If $\bmf{u}$ is a Lam\'e eigenfunction to \eqref{eq:lame}, then the following integral identity holds
\begin{equation}\label{eq:CGO2}
	I_3=I_1^+ + I_1^-+I_2,
\end{equation}
where
\begin{subequations}
	\begin{align}
		 I_1^\pm&=\int_{\Gamma_h^\pm } \left[\left({T}_{\nu} \bmf{u}\right) \cdot \bmf{v}-\left({T}_{\nu} \bmf{v}\right) \cdot \bmf{u}\right] \mathrm{d} \sigma, \label{eq:CGO6}\\
		  I_2&=\int_{\Lambda_h} \left[\left({T}_{\nu} \bmf{u}\right) \cdot \bmf{v}-\left({T}_{\nu} \bmf{v}\right) \cdot \bmf{u}\right] \mathrm{d} \sigma, \label{eq:CGO6a}\\
		 I_3&
		 =-\kappa \int_{S_h} \bmf{u} \cdot \bmf{v} \rmd \bmf{x}. \label{eq:I3 int}
	\end{align}
	\end{subequations}
	Furthermore, we have
	\begin{equation}\label{eq:I2I4}
		\begin{split}
			\left|I_2\right|&  \leqslant \mathcal{C}_{{\mathcal K },B_h,\mu,\lambda}\|\bmf{u}\|_{H^2\left( {{\mathcal K}} \cap B_h \right)}\left(1+s\right){\mathrm e}^{-\delta_{{\mathcal K} } s \sqrt{h}}
		\end{split}
	\end{equation}
	which is exponentially decays as $s\rightarrow +\infty$. Here  $\delta_{\mathcal{K}} $ is a positive constant defined in \eqref{eq:volume}.
\end{lem}
\begin{proof}
We first  prove \eqref{eq:CGO2}. Recall that $S_\varepsilon$ is defined in \eqref{eq:Seps}.  Since $\bmf{u}$ and $\bmf{v}$ are $H^2( S_h\backslash S_\varepsilon )$, from \eqref{eq:CGO1}, we have
\begin{equation}\label{eq:CGO3}
\int_{S_h \backslash{S_\varepsilon}} \left(\mathcal{L} \bmf{u}\right) \cdot \bmf{v} \rmd x=\int_{\partial ({S_h}\backslash{S_\varepsilon})} \left[\left({T}_{\nu} \bmf{u}\right) \cdot \bmf{v}-\left({T}_{\nu} \bmf{v}\right) \cdot \bmf{u}\right] \mathrm{d} \sigma.
\end{equation}
Since $\mathbf{u}$ is a Lam\'e eigenfunction to \eqref{eq:lame}, then $\mathcal{L} (\bmf{u})=-\kappa  \bmf{u} \text { in } \Omega \subset  \mathbb{R}^{2}$. Moreover, $\bmf{u}$ and $\bmf{v}$ are $H_{loc}^2(\mathbb{R}^2)$ functions, it yields that
\begin{equation}\label{eq:CGO4}
\lim_{\varepsilon\rightarrow 0^+ } \int_{S_\varepsilon} \left(\mathcal{L} \bmf{u}\right) \cdot \bmf{v} \rmd x=-\kappa  \lim_{\varepsilon\rightarrow 0^+ }  \int_{S_\varepsilon} \bmf{u} \cdot \bmf{v} \rmd x =0.
\end{equation}

Recall that $\Lambda_\varepsilon$ is defined in \eqref{eq:Seps}. It is easy to see that
\begin{equation}\label{eq:CGO5}
\int_{\partial ({S_h}\backslash{S_\varepsilon})} \left[\left({T}_{\nu} \bmf{u}\right) \cdot \bmf{v}-\left({T}_{\nu} \bmf{v}\right) \cdot \bmf{u}\right] \mathrm{d} \sigma=I_1^+ + I_1^- + I_2 -I_\varepsilon^+- I_\varepsilon^- + I_{\Lambda_\varepsilon}.
\end{equation}
where
 \begin{subequations}
 \begin{align}
 	 I_\varepsilon^\pm&=\int_{\Gamma_{(0,\varepsilon)}^\pm } \left[\left({T}_{\nu} \bmf{u}\right) \cdot \bmf{v}-\left({T}_{\nu} \bmf{v}\right) \cdot \bmf{u}\right] \mathrm{d} \sigma,  \label{eq:CGO7} \\
 	 I_{\Lambda_\varepsilon}&=\int_{\Lambda_\varepsilon} \left[\left({T}_{\nu} \bmf{u}\right) \cdot \bmf{v}-\left({T}_{\nu} \bmf{v}\right) \cdot \bmf{u}\right] \mathrm{d} \sigma. \label{eq:CGO9}
 \end{align}
 \end{subequations}
 Here the line segements $\Gamma_{(0,\varepsilon)}^\pm $ are defined in \eqref{eq:Seps}.
From \eqref{eq:27a} and \eqref{eq:27b}, we know that
\begin{equation}\label{eq:220 eps}
	\lim_{\varepsilon \rightarrow 0^+ } I_\varepsilon^\pm =0.
\end{equation}
By setting $h=\varepsilon$ in \eqref{eq:arc}, it is readily seen that
\begin{equation}\label{eq:221 lam}
	\lim_{\varepsilon \rightarrow 0^+ } I_{\Lambda_\varepsilon}  =0.
\end{equation}
In \eqref{eq:CGO3}, letting $\varepsilon \rightarrow 0^+$, using \eqref{eq:CGO4}, \eqref{eq:220 eps} and \eqref{eq:221 lam},  we obtain that
\begin{equation*}
I_3=-\kappa \int_{S_h}  \bmf{u} \cdot \bmf{v} \rmd x
=\int_{\partial S_h} \left[\left({T}_{\nu} \bmf{u}\right) \cdot \bmf{v}-\left({T}_{\nu} \bmf{v}\right) \cdot \bmf{u}\right] \mathrm{d} \sigma,
\end{equation*}
from which we prove \eqref{eq:CGO2}.

The proof is complete. 
\end{proof}

\begin{lem}\label{lem:the order of r exp2}
For a given $\zeta(\varphi)=-\mathrm{e}^{\bsi \frac{\varphi}{2}} \in \mathbb{C}$ and $\ell=\frac{m}{2},m=0,1,2,...$. we have
\begin{equation*}
\begin{aligned}
\int_{0}^{h}r^\ell \mathrm{e}^{s \sqrt{r} \zeta(\varphi)}\rmd r= & \frac{2}{s^{2 \ell+2}}\left(\frac{\left(-1\right)^{2\ell}\left(2 \ell+1\right)!}{\zeta(\varphi)^{2l+2}}\right.\\
& + \mathrm{e}^{s \sqrt{h} \zeta(\varphi)}\sum_{j=0}^{2 \ell+1}\frac{\left(-1\right)^j\left(2 \ell+1\right)!}{\left(2 \ell+1-j\right)!\zeta(\varphi)^{j+1}}
\left(s^2 h\right)^\frac{\left(2\ell+1-j\right)}{2}\bigg),\\
\int_{0}^{h}r^\ell \mathrm{e}^{s \sqrt{r} \Re(\zeta(\varphi)) }\rmd r= & \frac{2}{s^{2 \ell+2}}\left(\frac{\left(-1\right)^{2\ell}\left(2 \ell+1\right)!}{\Re( \zeta(\varphi)) ^{2l+2}}\right.\\
& + \mathrm{e}^{s \sqrt{h} \Re( \zeta(\varphi)) }\sum_{j=0}^{2 \ell+1}\frac{\left(-1\right)^j\left(2 \ell+1\right)!}{\left(2 \ell+1-j\right)! \Re( \zeta(\varphi))^{j+1}}
\left(s^2 h\right)^\frac{\left(2\ell+1-j\right)}{2}\bigg),\\
\end{aligned}
\end{equation*}
Furthermore, if $\mathfrak{R}\left(\zeta\left(\varphi\right)\right)<0$, we have the following asymptotic expansion:
\begin{subequations}
\begin{align}
\int_{0}^{h}r^\ell \mathrm{e}^{s \sqrt{r} \zeta(\varphi)}\rmd r&=\frac{2}{s^{2 \ell+2}}\frac{\left(-1\right)^{2\ell}\left(2 \ell+1\right)!}{\zeta(\varphi)^{2 \ell+2}}+\Oh \left(s^{-\frac{1}{2}} \mathrm{e}^{\sqrt{sh} \zeta(\varphi)}\right),\label{eq:268a int} \\
\int_{0}^{h}r^\ell \mathrm{e}^{s \sqrt{r} \Re(\zeta(\varphi))}\rmd r&=\frac{2}{s^{2 \ell+2}}\frac{\left(-1\right)^{2\ell}\left(2 \ell+1\right)!}{\Re( \zeta(\varphi))^{2 \ell+2}}+\Oh \left(s^{-\frac{1}{2}} \mathrm{e}^{\sqrt{sh} \Re( \zeta(\varphi))}\right). \label{eq:268b int}
\end{align}
\end{subequations}
as s $\rightarrow +\infty$, where $\ell=\frac{m}{2},m=0,1,2,
\ldots$.
\end{lem}

\begin{proof}
By induction and direct verifications, one can show the lemma. We skip the details. 
\end{proof}

In the next lemma, we derive an upper bound for the integral $I_3$ defined in \eqref{eq:I3 int}.

\begin{lem}\label{lem:I3 exp}
Recall that the Lam\'e eigenfunction   $\mathbf{u}$ to \eqref{eq:lame}  has the radial wave expression \eqref{eq:u} at the origin and  $I_3$ is defined by \eqref{eq:I3 int}, then one has
 \begin{equation}\label{eq:I3}
  |I_3| \leqslant  \left|k_p ^2 a_0 - \bsi k_s ^2 b_0 \right|  \frac{\kappa  \varphi_0\Gamma\left(6\right)  }{\delta_{{\mathcal K}}^{6} } s^{-6} + \frac{2\kappa \varphi_0\Gamma\left(8\right) S_1 }{\delta_{{\mathcal K}}^{8} } s^{-8},
\end{equation}
as $s\rightarrow +\infty$, where $\delta_{\mathcal{K}} $ is a positive constant defined in \eqref{eq:volume}
 and  $S_1$ is defined in \eqref{eq:uv}. Recall that $S_1(\ell) $ is defined in  \eqref{eq:uv1}. Furthermore, if $a_0=b_0=\cdots=a_{\ell-1}=b_{\ell-1}=0$,
 \begin{equation}\label{eq:I3+}
  |I_3| \leqslant  \left|k_p ^{\ell+2} a_{}\ell - \bsi k_s ^{\ell+2} b_{\ell} \right|  \frac{  \kappa \varphi_0\Gamma(2 \ell+6)  }{2^{\ell}(\ell+1)!\delta_{{\mathcal K}}^{2 \ell+6} } s^{-2 \ell-6} + \frac{2  \kappa \varphi_0\Gamma(2 \ell+8) S_1(\ell) }{\delta_{{\mathcal K}}^{8} } s^{-2\ell-8},
\end{equation}
as $s\rightarrow +\infty$. 
 \end{lem}

\begin{proof}
Submitting \eqref{eq:uv old} into  \eqref{eq:I3 int},  using  \eqref{eq:volume} and \eqref{eq:uv},   we derive that
\begin{equation*}
\begin{aligned}
	|I_3| &\leqslant
		 \kappa \int_{\mathcal K}   |\bmf{u} \cdot \bmf{v}|  \rmd \bmf{x} \\
		 &\leqslant  \kappa  \frac{\left|k_p ^2 a_0 - \bsi k_s ^2 b_0 \right| }{4} \int_{\mathcal K} r^2 {\mathrm e}^{s \sqrt{ r} \Re(\zeta(\varphi )) }    \rmd r \rmd \varphi + \kappa   S_1  \int_{\mathcal K}  r^3 {\mathrm e}^{s \sqrt{ r} \Re(\zeta(\varphi )) }    \rmd r
		 %
\end{aligned}
\end{equation*}
from which we complete the proof of \eqref{eq:I3}. By virtue of \eqref{eq:uv1},  \eqref{eq:I3+} can be proved in a similar way. 
\end{proof}

\section{Generalized Holmgren's principle with the presence of singular lines}\label{sect:3}

In this section, we prove that if a generic Lam\'e eigenfunction $\bmf{u}$ to \eqref{eq:lame} possesses a singular line $\Gamma_h$ in $\Omega$ as defined in Definition~\ref{def:generalized line}, then $\bmf{u}$ is identically zero. According to our discussion made at the beginning of Section~\ref{sect:2}, without loss of generality, we can assume that $\Gamma_h$ is $\Gamma_h^-$ as defined in \eqref{eq:gamma_pm}. Moreover, we can assume that the point $\bmf{x}_0$ involved in Definition~\ref{def:generalized line} is the origin, namely $\bmf{x}_0=\bmf{0}$. It is clear that the unit normal vectors $\nu $ to $\Gamma_h^-$ is $(0,\pm 1)^\top$. In this paper we choose $(0,-1)^\top$ as the unit normal vector $\nu $ to $\Gamma_h^-$. In such a case, the following conditions involved in Definition~\ref{def:generalized line}
\begin{equation}\label{eq:31}
\bmf{u}(\bmf{x}_0 )=\bmf{0} \mbox{ and/or } \boldsymbol{\tau} \cdot \partial_{\nu} \bmf{u} |_{\bmf{x}={\bmf{x}}_0 }=0,
\end{equation}
turn into
$$
\bmf{u}(\bmf{0} )=\bmf{0} \mbox{ and/or }  \partial_2 u_1(\bmf{0} )=0. 
$$

\subsection{The case with a singular rigid line}
\begin{lem}\label{lem:rigid}
Let $\mathbf{u}$ be a generalized Lam\'e eigenfunction to \eqref{eq:lame} with the radial wave  expansion given in \eqref{eq:u} around the origin. Suppose there exists $\Gamma_h^-\in {\mathcal R}_\Omega^{\kappa}  $. Then one has
\begin{equation}\label{eq:1}
 \bigg\{  \begin{array}{l} k_p a_1 + \bsi k_s b_1=0,\\
 k_p^3a_1-\bsi k_s^3 b_1=0,
  \end{array}
\end{equation}
and
\begin{equation}\label{eq:2}
\bigg\{
    \begin{array}{l}
  k_p ^2 a_0+ \bsi k_s^2 b_0- k_p^2  a_2- \bsi k_s ^2 b_2=0,\\
  k_p ^2 a_0- \bsi k_s^2 b_0=0.
 \end{array}
\end{equation}
Moreover, it holds that
\begin{equation}\label{eq:a1b1 lem}
	a_1=b_1=0.
\end{equation}
Furthermore, suppose that
\begin{equation}\label{eq:lem cond}
	a_\ell=b_\ell=0
\end{equation}
where $\ell=0,\ldots,m-1$ and $m\in {\mathbb N}$ with $m\geq 2 $, then
$$
a_\ell=b_\ell=0, \quad \forall \ell \in \mathbb N \cup \{0\}.
$$
\end{lem}

\begin{proof}
Since $\Gamma_h^-$ is a rigid line of $\bmf{u}$, then $\bmf{u}\big|_{\Gamma_h^- }=\bmf{0}$. Therefore from \eqref{eq:u} and by noting $\varphi=0$ on $\Gamma_h^-$, we have for $0 \leqslant r\leqslant h$ that
\begin{equation}\label{eq:u=0}
 \begin{aligned}
 \mathbf{0}=  \sum_{m=0}^{\infty} & \left\{ \frac{k_p}{2} a_{m} \left\{J_{m-1}\left(k_{p} r\right) \mathbf{e}_1
 -J_{m+1}\left(k_{p}r\right)\mathrm{e}^{\bsi \varphi}\mathbf{e}_2 \right\}\right. + \frac{\bsi k_s}{2} b_{m}  \left\{J_{m-1}\left(k_{s} r\right) \mathbf{e}_1
  +J_{m+1}\left(k_{s} r\right)
  \mathbf{e}_2
   \right\} \bigg\}. 
 \end{aligned}
\end{equation} 
From \eqref{eq:J-1}, we know that $J_{-1}(k_pr)=-J_{1}(k_pr) $ and $J_{-1}(k_sr)=-J_{1}(k_sr) $.  Using  Lemma \ref{lem:co exp}, comparing the coefficients of the term $r^0$ in both sides of \eqref{eq:u=0}, we obtain \eqref{eq:1}. Similarly, from Lemma \ref{lem:co exp}, we compare the coefficients of the term $r^1$ in both sides of \eqref{eq:u=0}, and obtain that
\begin{equation}\nonumber
\left(-k_p^2  a_0-\bsi k_s^2  b_0+k_p^2  a_2+\bsi k_s^2  b_2\right)\mathbf{e}_1 -\left(k_p^2 a_0-\bsi k_s^2  b_0\right)\mathbf{e}_2=\bmf{0}.
\end{equation}
Since $\mathbf{e}_1$ and $\mathbf{e}_2$ are linearly independent, we can obtain \eqref{eq:2}. Similarly, comparing the coefficients of the term $r^2$ in both sides of \eqref{eq:u=0}, we have
\begin{equation}\nonumber
\left(-2k_p^3  a_1-2 \bsi k_s^3  b_1+k_p^3 a_3+\bsi k_s^3  b_3\right)\mathbf{e}_1 -\left(k_p^3 a_1-\bsi k_s^3  b_1\right)\mathbf{e}_2=\bmf{0},
\end{equation}
which implies that the second equation of \eqref{eq:1} holds. The determinant of the coefficient matrix of \eqref{eq:1} is
$$
\left| \begin{matrix}
	k_p & \bsi k_s\\
	k_p^3 & -\bsi k_s^3
\end{matrix} \right|=-\bsi k_pk_s (k_p^2+k_s^2) \neq 0,
$$
then we know that \eqref{eq:a1b1 lem} holds.

Substituting \eqref{eq:lem cond} into \eqref{eq:u=0}, it yields that
\begin{equation}\label{eq:u=0 new}
 \begin{aligned}
 \mathbf{0}=  \sum_{\ell=m}^{\infty} & \left\{ \frac{k_p}{2} a_{\ell } \left\{J_{\ell -1}\left(k_{p} r\right) \mathbf{e}_1
 -J_{\ell +1}\left(k_{p}r\right)\mathrm{e}^{\bsi \varphi}\mathbf{e}_2 \right\}\right.\\
 & + \frac{\bsi k_s}{2} b_{\ell }  \left\{J_{\ell -1}\left(k_{s} r\right) \mathbf{e}_1
  +J_{\ell +1}\left(k_{s} r\right)
  \mathbf{e}_2
   \right\} \bigg\},\quad 0 \leqslant r\leqslant h.
 \end{aligned}
\end{equation}
Therefore the lowest order of $r$ in the right hand side of \eqref{eq:u=0 new} is $m-1$. Hence comparing the coefficients of the term $r^{m-1}$ in both sides of \eqref{eq:u=0 new}, we obtain that
\begin{equation}\label{eq:6}
 k_p ^{m } a_m+ \bsi k_s^{m}  b_m=0.
\end{equation}
 Comparing the coefficients of the term $r^{m+1}$ in both sides of \eqref{eq:u=0 new}, which are related to $J_{m-1}(k_p r)$, $J_{m-1}(k_s r)$,  $J_{m+1}(k_p r)$ and $J_{m+1}(k_s r)$, one has
\begin{equation}\nonumber
\begin{aligned}
& \left(-\left(m+1\right) k_p^{m+2} a_m-\bsi (m+1) k_s^{m+2} b_m+k_p^{m+2}  a_{m+2}+\bsi k_s^{m+2}  b_{m+2}  \right) \mathbf{e}_1\\
& - \left(k_p^{m+2} \mathrm{e}^{2 \bsi\varphi_0} a_m-\bsi k_s^{m+2}  b_m\right)\mathbf{e}_2=\bmf{0}.
\end{aligned}
\end{equation}
Since $\mathbf{e}_1$ and $\mathbf{e}_2$ are linear independent, it yields that
\begin{equation}\label{eq:7}
\left\{
    \begin{array}{l}
-\left(m+1\right) k_p^{m+2} a_m-\bsi (m+1) k_s^{m+2} b_m+k_p^{m+2}  a_{m+2}+\bsi k_s^{m+2}  b_{m+2} =0,\\
 k_p^{m+2} a_m-\bsi k_s^{m+2} b_m=0,\\
 \end{array}
 \right.
\end{equation}
Combining \eqref{eq:6} and \eqref{eq:7}, we obtain that
\begin{equation}\notag
\left\{
\begin{array}{l}
 k_p ^m a_m+ \bsi k_s^m b_m=0,\\
k_p^{m+2} a_m-\bsi k_s^{m+2} b_m=0.
\end{array}
\right.
\end{equation}
Since
\begin{equation}\notag
\left|
\begin{array}{cc}
k_p^m    &   \bsi k_s^m\\
k_p^{m+2}  &   -\bsi k_s^{m+2}\\
\end{array}
\right|=-\bsi k_p^m k_s^m \left(k_p^2 +k_s^2\right)\neq 0,
\end{equation}
one readily has that $a_m=b_m=0$. In an inductive manner,  we can prove that $a_\ell=b_\ell=0$ for $\ell=m+1,\ldots$.

The proof is complete. 
\end{proof}

\begin{lem}\label{lem:34 three conds}
Let $\mathbf{u}=(u_\ell)_{\ell=1}^2$ be a Lam\'e eigenfunction to \eqref{eq:lame} with the radial wave expansion \eqref{eq:u} around the origin.  If $\bmf{u}(\bmf{0})=\bmf{0}$, then
\begin{equation}\label{eq:u0}
k_p a_1 + \bsi k_s b_1=0.
\end{equation}
Recall that $u_1$ has the expansion \eqref{eq:u comp}. If $\partial_2 u_1 ({\mathbf 0})=0$, then
\begin{equation}\label{eq:u1}
-2 k_s^2 b_0+\bsi k_p^2 a_2 -k_s^2 b_2=0.
\end{equation}
\end{lem}

\begin{proof}
 Since $\bmf{u}(\bmf{0})=\bmf{0}$, substituting  $r=0$ in  \eqref{eq:u} we can prove \eqref{eq:u0}.
  From \eqref{eq:u3 par} in the Appendix, it is easy to know that
  \begin{equation}\notag
  \partial_1 u_1=\frac{\partial u_1}{\partial r}
  \end{equation}
   at $\varphi=0$.
   Substituting \eqref{eq:J2} into \eqref{eq:u1 par}, we can obtain that
\begin{equation}\label{eq:u3 new}
\begin{aligned}
\frac{1}{r}\frac{\partial u_1}{\partial \varphi}\Big|_{\varphi=0}=&\sum_{m=0}^\infty  \bigg\{ \bsi k_p^2 a_m J_{m-2}(k_p r)-\bsi k_p^2 a_m J_{m+2}(k_p r)- k_s^2 b_m J_{m-2}(k_s r)\\
&\quad -2 k_s^2 b_m J_m(k_s r)-k_s^2 b_m J_{m+2}(k_s r)\bigg\}.
\end{aligned}
\end{equation}
Substituting \eqref{eq:u3 new} into \eqref{eq:u3 par} and evaluating the resulting equality of   \eqref{eq:u3 par} at $r=0$ and $\varphi=0$  since $\partial_2 u_1 (\bmf{0})=0$, one can prove \eqref{eq:u1}.
%
%
%
\end{proof}

\begin{thm}\label{Thm:31 singular rigid}
	Let $\mathbf{u}\in L^2(\Omega)^2$ be a solution to \eqref{eq:lame}. If there exits a singular rigid line $\Gamma_h\subset\Omega$ of $\bmf{u}$, then $\bmf{u}\equiv \bmf{ 0}$.
\end{thm}

\begin{proof}
	Suppose that there exits a singular rigid line $\Gamma_h$ of $\bmf{u}$ as described at the beginning of this section. Then we have
	\begin{equation}\label{eq:311 cond}
		\partial_2 \bmf{u}_1(\bmf{x}_0)=0,\quad \bmf{x}_0\in \Gamma_h.
	\end{equation}
	By virtue of \eqref{eq:1} we know that
	\begin{equation}\label{eq:a1b1=0}
			a_1=b_1=0
	\end{equation}
since
	$$
	\left|\begin{matrix}
		k_p & \bsi k_s\\
		k_p^3 &-\bsi k_s^3
	\end{matrix} \right|=-\bsi k_pk_s(k_p^2+k_s^2)\neq 0.
	$$
Substituting \eqref{eq:a1b1=0} into \eqref{eq:u=0}, comparing the coefficients of the term $r^3$ in both sides of \eqref{eq:u=0}, we have
\begin{equation}\label{eq:313 eq}
	3k_p^4 a_0-3\bsi k_s^4 b_0-k_p^4 a_2+\bsi k_s^4 b_2=0.
\end{equation}
\eqref{eq:311 cond} implies that \eqref{eq:u1} is satisfied. Combing \eqref{eq:2} \eqref{eq:u1}, \eqref{eq:313 eq}, we have
$$
\begin{bmatrix}
	k_p^2 &\bsi k_s^2 & -k_p^2 & -\bsi k_s^2\\
	k_p^2&-\bsi k_s^2 &0 &0\\
	0&-2k_s^2&\bsi k_p^2& -k_s^2\\
	3k_p^4&-3\bsi k_s^4 &-k_p^4 & \bsi k_s^4
\end{bmatrix} \begin{bmatrix}
	a_0\\b_0\\a_2\\
	b_2
\end{bmatrix}=\bmf{0},
$$
whose determinant is $-4\bsi k_p^4 k_s^4 (k_p^2+k_s^2) \neq 0$, then we know that $a_0=b_0=a_2=b_2=0$. Using Lemma \ref{lem:rigid}, we can prove $a_\ell=b_\ell=0$ for $\ell\in \mathbb N$ and $\ell\geq 3$, which induces that $\bmf{u}\equiv \bmf{0} $ in $\Omega$ by Proposition \ref{prop:21}. 

The proof is complete.
\end{proof}

\subsection{The case with a singular traction-free line}

\begin{lem}\label{lem:tranction free}
Let $\mathbf{u}$ be a Lam\'e eigenfunction to \eqref{eq:lame} with the radial wave  expansion \eqref{eq:u} around the origin.  Suppose that  $\Gamma_h^- \in {\mathcal T}_\Omega^{\kappa}  $. Then we have
\begin{equation}\label{eq:8}
\left\{
\begin{array}{l}
\bsi k_p^2 a_2- k_s^2  b_2=0,\\
a_0 =0,
\end{array}
\right.
\end{equation}
and
\begin{equation}\label{eq:10}
\left\{
    \begin{array}{l}
k_s^3  b_1+\bsi k_p^3  a_3-k_s^3   b_3=0,\\
 a_1=0.
 \end{array}
 \right.
\end{equation}
Furthermore, suppose that $a_\ell=b_\ell=0$ for $\ell=0,1$, then
\begin{equation}\label{eq:lem33}
	a_\ell=b_\ell=0,\quad \forall \ell \in \mathbb N\cup \{0\}.
\end{equation}
\end{lem}

\begin{proof}
Since $\Gamma_h^-$ is a traction-free  line of $\bmf{u}$, then $T_\nu \bmf{u}\big|_{\Gamma_h^- }=\bmf{0}$, we have from \eqref{eq:Tu2} that
\begin{equation}\label{eq:Tu2=0}
 \begin{split}
\bmf{0}=  \sum_{m=0}^{\infty} & \left\{-\frac{\bsi k_p^2}{2} a_{m} \mu   J_{m-2}(k_p r) \mathbf{e}_1 - \frac{\bsi k_p^2}{2} a_{m}  (\lambda+\mu) J_m(k_p r) \mathbf{e}_1\right.\\
 & + \frac{k_s^2}{2} b_{m}  \mu J_{m-2}\left(k_{s} r\right) \mathbf{e}_1 + \frac{\bsi k_p^2}{2} a_{m}  \left(\lambda+\mu\right) J_m\left(k_{p} r\right) \mathbf{e}_2\\
 & + \frac{\bsi k_p^2}{2} a_{m}  \mu J_{m+2}\left(k_{p} r\right) \mathbf{e}_2 + \frac{k_s^2}{2} b_{m} \mu J_{m+2}\left(k_{s} r\right) \mathbf{e}_2
 \bigg\},
 \end{split}
\end{equation}
where $r\in [0,h] $.  Using \eqref{eq:J-1} and  Lemma \ref{lem:co exp}, comparing the coefficients of the term $r^0$ in both sides of \eqref{eq:Tu2=0}, we obtain
\begin{equation}\nonumber
\left(\bsi k_p^2 \left(\lambda+\mu\right) a_0+\bsi k_p^2 \mu  a_2-k_s^2 \mu  b_2\right)\mathbf{e}_1-\bsi k_p^2 \left(\lambda+\mu\right) a_0 \mathbf{e}_2=\bmf{0}.
\end{equation}
From the fact that $\mathbf{e}_1$ and $\mathbf{e}_2$ are linearly independent, we prove \eqref{eq:8} by using \eqref{eq:convex}.  Again comparing the coefficients of the term $r^1$ in both sides of \eqref{eq:Tu2=0}, we obtain
\begin{equation}\nonumber
\left(\bsi k_p^3 \lambda a_1+k_s^3 \mu b_1+\bsi k_p^3 \mu a_3-k_s^3 \mu  b_3\right)\mathbf{e}_1-\bsi k_p^3  \left(\lambda+\mu\right)a_1 \mathbf{e}_2=\bmf{0},
\end{equation}
from which one can easily see that \eqref{eq:8} holds by using \eqref{eq:convex}.

Suppose that  $a_\ell=b_\ell=0$ for $\ell=0,1$, then we want to prove \eqref{eq:lem33}. Substituting $a_\ell=b_\ell=0$ ($\ell=0,1$) into \eqref{eq:Tu2=0} and comparing the coefficients of $r^2$ in the resulting equation \eqref{eq:Tu2=0}, which are related  to $J_0(k_p r)$, $J_0(k_s r)$ and $J_2(k_p r)$, $J_2(k_s r)$ for the indexes  $m=2$ and $m=4$ in \eqref{eq:Tu2=0}, we obtain that
\begin{equation}\nonumber
\left(\bsi k_p^4 \left(\lambda-\mu\right) a_2+2k_s^4 \mu b_2+\bsi k_p^4\mu a_4-k_s^4 \mu b_4\right)\mathbf{e}_1-\bsi k_p^4 \left(\lambda+\mu\right) a_2\mathbf{e}_2=\bmf{0}.
\end{equation}
Since $\mathbf{e}_1$ and $\mathbf{e}_2$ are linear independent, we can deduce that
\begin{equation}\label{eq:12}
\left\{
    \begin{array}{l}
 \bsi k_p^4 \left(\lambda-\mu\right) a_2+2k_s^4 \mu b_2+\bsi k_p^4\mu a_4-k_s^4 \mu b_4=0,\\
 a_2=0,\\
 \end{array}
 \right.
\end{equation}
Since $a_2=0$ from \eqref{eq:12},  we have $b_2=0$ from  \eqref{eq:8}.
Thus \eqref{eq:12} can be written as $
\bsi k_p^4 a_4-k_s^4 b_4=0$.
Repeating the above argument in an inductive manner, we can prove that
\begin{equation}\label{eq:albl=0T}
	 a_\ell=b_\ell=0
\end{equation}
for $\ell=0,1,2,\ldots, m-1$ where $m\in \mathbb N$ is fixed and $m \geq 3$. Next, we prove $a_m=b_m=0$. Substituting \eqref{eq:albl=0T} into \eqref{eq:Tu2=0}, it yields that
\begin{equation}\label{eq:Tu2=0 new}
 \begin{split}
\bmf{0}=  \sum_{\ell=m}^{\infty} & \left\{-\frac{\bsi k_p^2}{2} a_{\ell } \mu   J_{\ell -2}(k_p r) \mathbf{e}_1 - \frac{\bsi k_p^2}{2} a_{\ell }  (\lambda+\mu) J_\ell (k_p r) \mathbf{e}_1\right.\\
 & + \frac{k_s^2}{2} b_{\ell }  \mu J_{\ell -2}\left(k_{s} r\right) \mathbf{e}_1 + \frac{\bsi k_p^2}{2} a_{\ell }  \left(\lambda+\mu\right) J_\ell \left(k_{p} r\right) \mathbf{e}_2\\
 & + \frac{\bsi k_p^2}{2} a_{\ell }  \mu J_{\ell +2}\left(k_{p} r\right) \mathbf{e}_2 + \frac{k_s^2}{2} b_{
\ell } \mu J_{\ell +2}\left(k_{s} r\right) \mathbf{e}_2
 \bigg\}.
 \end{split}
\end{equation}
Therefore the lowest order of $r$ in the left sides of \eqref{eq:Tu2=0 new} is $m-2$.
Comparing the coefficients of  $r^{m-2}$ in both sides of \eqref{eq:Tu2=0 new} and using \eqref{eq:Jm ex}, we have
\begin{equation}\label{eq:16}
 \bsi k_p ^m a_m- k_s^m b_m=0,
\end{equation}
since $\mu>0$. Again comparing the coefficients of  $r^{m}$ in both sides of \eqref{eq:Tu2=0 new}, which are related to $J_{m-2}(k_p r)$, $J_{m-2}(k_s r)$ and $J_m(k_p r)$, $J_m(k_s r)$,  we can derive that
\begin{equation}\nonumber
\begin{aligned}
& \left[\bsi k_p^{m+2} \left(\lambda-\left(m-1\right)\mu\right) a_m+m k_s^{m+2} \mu b_m+\bsi k_p^{m+2} \mu a_{m+2}-k_s^{m+2} \mu b_{m+2}\right]\mathbf{e}_1\\
& - \bsi k_p^{m+2} \left(\lambda+\mu\right) a_m \mathbf{e}_2=\bmf{0}.
\end{aligned}
\end{equation}
Since $\mathbf{e}_1$ and $\mathbf{e}_2$ are linearly independent, we can deduce that
\begin{equation}\label{eq:17}
\left\{
\begin{array}{l}
\bsi k_p^{m+2} \left(\lambda-\left(m-1\right)\mu\right) a_m+m k_s^{m+2} \mu b_m+\bsi k_p^{m+2} \mu a_{m+2}-k_s^{m+2} \mu b_{m+2}=0,\\
a_m=0.
\end{array}
\right.
\end{equation}
Combining \eqref{eq:16}and \eqref{eq:17},  we derive that $a_m=b_m=0$. Using the above recursive procedure, we can prove that $a_m=b_m=0$ for $m \in \mathbb{N}\cup \{0\}$.

The proof is complete. 
\end{proof}

\begin{thm}\label{thm:traction free line}
Let $\mathbf{u}\in L^2(\Omega)^2$ be a solution to \eqref{eq:lame}. If there exits a singular traction-free line $\Gamma_h\subset\Omega$ of $\bmf{u}$, then $\bmf{u}\equiv \bmf{ 0}$.

%
\end{thm}


\begin{proof}
Suppose that there exits a singular traction-free line $\Gamma_h$ of $\bmf{u}$ as described at the beginning of this section.
Recall that $\bmf{u}$ has the expansion \eqref{eq:u} around the origin. From Lemma \ref{lem:tranction free}, we know that \eqref{eq:8} holds. By virtue of the fact that $\bmf{u}(\bmf{0})=\bmf{0} \mbox{ and }\boldsymbol{\tau} \cdot \partial_{\nu} \bmf{u} (\bmf{0})=0$, from Lemma \ref{lem:34 three conds}, we know that  \eqref{eq:u1} is satisfied. Substituting \eqref{eq:8} into \eqref{eq:u1}, we can obtain that $
b_0=0.
$
Furthermore, substituting the second equation in \eqref{eq:10} into \eqref{eq:u0}, one can derive that
$
b_1=0,
$
which implies that  \eqref{eq:10} can be rewritten as
$ \bsi k_p^3 a_3-k_s^3 b_3=0$.

By now we have proven that $a_\ell=b_\ell=0$ for $\ell=0,1$. Then from Lemma \ref{lem:tranction free}, we have that \eqref{eq:lem33} holds. Therefore, from Proposition \ref{prop:21}, we know that $\bmf{u}\equiv \bmf{0}$ in $\Omega$.

The proof is complete.
\end{proof}

\subsection{The case with a singular impedance line}

\begin{lem}\label{lem:impedance eqn}
Let $\mathbf{u}$ be a solution to \eqref{eq:lame} with the radial wave expansion \eqref{eq:u} around the origin.  Suppose that there is an impedance  line $\Gamma_h^-$ of $\bmf{u}$ with a constant impedance parameter $0\neq \eta \in \mathbb C$. Then we have
\begin{equation}\label{eq:Tu+u1}
\left\{
    \begin{array}{l}
 \eta k_p a_1+\bsi \eta k_s b_1 -\bsi k_p^2 \mu a_2+k_s^2 \mu b_2=0,\\
 a_0=0,
 \end{array}
 \right.
\end{equation}
and
\begin{equation}\label{eq:Tu+u2}
\left\{
    \begin{array}{l}
 - k_s^3 \mu b_1 + \eta k_p^2 a_2+\bsi \eta k_s^2 b_2 +\bsi k_p^3 \mu a_3-k_s^3 \mu b_3=0,\\
 a_1=0.
 \end{array}
 \right.
\end{equation}
Furthermore, if $a_\ell=b_\ell=0$ for $\ell=0,1$, then
\begin{equation}\label{eq:a1bl=0}
	a_\ell=b_\ell=0, \quad \forall \ell \in \mathbb N \cup \{0\}.
\end{equation}
\end{lem}

\begin{proof}
Since $(T_\nu \bmf{u}+\eta \bmf{u}) \big|_{\Gamma_h^-}=\bmf{0}$, using  \eqref{eq:Tu4} and noting $\varphi=0$ on $\Gamma_h^-$, we have for $0 \leqslant r\leqslant h$ that
\begin{equation}\label{eq:Tu+u=0}
 \begin{aligned}
&\bmf{0}=\sum_{m=0}^{\infty}  \biggl \{-\frac{\bsi {k}_p^2}{2} a_{m} \Big[ \mu   J_{m-2}(k_p r) \mathbf{e}_1  +   (\lambda+\mu) J_m(k_p r) \mathbf{e}_1\\
&-  \left(\lambda+\mu\right) J_m\left(k_{p} r\right) \mathbf{e}_2 -  \mu J_{m+2}\left(k_{p} r\right) \mathbf{e}_2\Big]  + \frac{k_s^2}{2} b_{m}\Big[   \mu J_{m-2}\left(k_{s} r\right) \mathbf{e}_1 +  \mu J_{m+2}\left(k_{s} r\right) \mathbf{e}_2 \Big] \\
 &+\frac{\eta k_p}{2} a_{m} \left[J_{m-1}\left(k_{p} r\right)\mathbf{e}_1
 -J_{m+1}\left(k_{p}r\right) \mathbf{e}_2 \right]+ \frac{\bsi \eta k_s}{2} b_{m} \left[J_{m-1}\left(k_{s} r\right)\mathbf{e}_1
  +J_{m+1}\left(k_{s} r\right)\mathbf{e}_2
 \right]
 \bigg\}. 
 \end{aligned}
\end{equation}
Using  Lemma \ref{lem:co exp}, comparing the coefficients of the term $r^0$ in both sides of \eqref{eq:Tu+u=0}, which are related  to $J_0(k_p r)$ and $J_0(k_s r)$  for the indexes  $m=0$, $m=1$ and $m=2$ in \eqref{eq:Tu+u=0},
we have
\begin{equation}\nonumber
\left[-\bsi k_p^2 (\lambda+\mu) a_0+\eta k_p a_1+\bsi \eta k_s b_1 -\bsi k_p^2 \mu a_2+k_s^2 \mu b_2\right]\mathbf{e}_1+\bsi k_p^2(\lambda+\mu)a_0 \mathbf{e}_2=\bmf{0}.
\end{equation}
Using the fact that $\mathbf{e}_1$ and $\mathbf{e}_2$ are linearly independent, we can obtain \eqref{eq:Tu+u1} since $k_p \neq 0$ and $\lambda+\mu >0$ from \eqref{eq:convex}.  Similarly,  comparing the coefficients of the term $r^1$ in both sides of \eqref{eq:Tu+u=0}, we can derive \eqref{eq:Tu+u2}.

Now we are in the position to prove \eqref{eq:a1bl=0} under the condition $a_\ell=b_\ell=0$ for $\ell=0,1$.
Since $a_1=b_1=0$,  \eqref{eq:Tu+u1} can be rewritten as
\begin{equation}\label{eq:Tu+u1a}
\bsi k_p^2 a_2-k_s^2 b_2=0.
\end{equation}
Substituting  $a_\ell=b_\ell=0$ into \eqref{eq:Tu+u=0}, where $\ell=0,1$, comparing the coefficients of $r^2$ in the resulting equation \eqref{eq:Tu+u=0}, which are related  to $J_0(k_p r)$, $J_0(k_s r)$ and $J_2(k_p r)$, $J_2(k_s r)$ for the indexes  $m=2$ $m=3$ and $m=4$ in \eqref{eq:Tu+u=0}, we can derive that
\begin{equation}\nonumber
\left[\bsi k_p^4 (\mu-\lambda) a_2 - 2 k_s^4 \mu b_2+\eta k_p^3 a_3+\bsi \eta k_s^3 b_3 -\bsi k_p^4 \mu a_4 + k_s^4 \mu b_4\right]\mathbf{e}_1+\bsi k_p^4(\lambda+\mu)a_2 \mathbf{e}_2=\bmf{0},
\end{equation}
Since $\mathbf{e}_1$ and $\mathbf{e}_2$ are linearly independent, we can obtain that
\begin{equation}\label{eq:Tu+u3}
\left\{
    \begin{array}{l}
 - 2 k_s^4 \mu b_2+\eta k_p^3 a_3+\bsi \eta k_s^3 b_3 -\bsi k_p^4 \mu a_4+k_s^4 \mu b_4=0,\\
 a_2=0,\\
 \end{array}
 \right.
\end{equation}
Combing \eqref{eq:Tu+u1a} with \eqref{eq:Tu+u3}, we can derive that $b_2=0$.

Therefore, it is easy to see that \eqref{eq:Tu+u2} can be rewritten as
\begin{equation}\label{eq:Tu+u2a}
\bsi k_p^3 a_3-k_s^3 b_3=0.
\end{equation}
By now we have proven that $a_\ell=b_\ell=0$ ($\ell=0,1,2$)  {if $\Gamma_h^- \Subset \Omega $ is an impedance line of $\bmf{u}$.}
Repeating the above argument in an inductive manner, suppose that we have proven
\begin{equation}\label{eq:albl=0Tu+u}
	 a_\ell=b_\ell=0
\end{equation}
for $\ell=0,1,2,\ldots, m-1$ where $m\in \mathbb N$ is fixed and $m \geq 3$. We next prove $a_m=b_m=0$. Substituting \eqref{eq:albl=0Tu+u} into \eqref{eq:Tu+u=0}, it yields that
\begin{equation}\label{eq:Tu+u=0new}
 \begin{aligned}
&\bmf{0}=\sum_{\ell=m}^{\infty}  \biggl \{-\frac{\bsi {k}_p^2}{2} a_{\ell} \Big[ \mu   J_{\ell-2}(k_p r) \mathbf{e}_1  +   (\lambda+\mu) J_\ell(k_p r) \mathbf{e}_1\\
&-  \left(\lambda+\mu\right) J_\ell\left(k_{p} r\right) \mathbf{e}_2 -  \mu J_{\ell+2}\left(k_{p} r\right) \mathbf{e}_2\Big]  
+ \frac{k_s^2}{2} b_{\ell}\Big[   \mu J_{\ell-2}\left(k_{s} r\right) \mathbf{e}_1 +  \mu J_{\ell+2}\left(k_{s} r\right) \mathbf{e}_2 \Big] \\
 &+\frac{\eta k_p}{2} a_{\ell} \left[J_{\ell-1}\left(k_{p} r\right)\mathbf{e}_1
 -J_{\ell+1}\left(k_{p}r\right) \mathbf{e}_2 \right] + \frac{\bsi \eta k_s}{2} b_{\ell} \left[J_{\ell-1}\left(k_{s} r\right)\mathbf{e}_1
  +J_{\ell+1}\left(k_{s} r\right)\mathbf{e}_2
 \right]
 \bigg\}.
 \end{aligned}
\end{equation}
Therefore the lowest order of $r$ in the left sides of \eqref{eq:Tu+u=0new} is $m-2$.
Comparing the coefficients of  $r^{m-2}$ in both sides of \eqref{eq:Tu+u=0new}, we have
\begin{equation}\label{eq:Tu+u4}
 \bsi k_p ^m a_m- k_s^m b_m=0,
\end{equation}
since $\mu>0$. Again comparing the coefficients of  $r^{m}$ in both sides of \eqref{eq:Tu+u=0new}, which are related to $J_{m-2}(k_p r)$, $J_{m-2}(k_s r)$ and $J_m(k_p r)$, $J_m(k_s r)$,  we can derive that
\begin{equation}\label{eq:339 im}
\begin{aligned}
& \big[-\bsi k_p^{m+2} \left(\lambda+(1-m)\mu\right) a_m - m k_s^{m+2} \mu b_m + \eta k_p^{m+1} a_{m+1} + \bsi \eta k_s^{m+1} \mu b_{m+1} \\
&- \bsi k_p^{m+2} \mu a_{m+2} + k_s^{m+2} \mu a_{m+2} \big]\mathbf{e}_1
+ \bsi k_p^{m+2} \left(\lambda+\mu\right) a_m \mathbf{e}_2=\bmf{0}.
\end{aligned}
\end{equation}
Since $\mathbf{e}_1$ and $\mathbf{e}_2$ are linearly independent, from \eqref{eq:339 im}, we can deduce that
\begin{equation}\label{eq:Tu+u5}
\left\{
    \begin{array}{l}
 - m k_s^{m+2} \mu b_m + \eta k_p^{m+1} a_{m+1} + \bsi \eta k_s^{m+1} \mu b_{m+1}
- \bsi k_p^{m+2} \mu a_{m+2} + k_s^{m+2} \mu a_{m+2}=0,\\
 a_m=0,\\
 \end{array}
 \right.
\end{equation}
Combining \eqref{eq:Tu+u4}and \eqref{eq:Tu+u5},  we derive that $a_m=b_m=0$. Using the above recursive procedure, we can prove that $a_m=b_m=0$ for $m \in \mathbb{N}\cup \{0\}$.

The proof is complete. 
\end{proof}

\begin{thm}\label{thm:impedance line}
Let $\mathbf{u}\in L^2(\Omega)^2$ be a solution to \eqref{eq:lame}. If there exits a singular impedance line $\Gamma_h\subset\Omega$ of $\bmf{u}$ with a constant impedance parameter $0\neq \eta \in \mathbb C$ as defined in \eqref{eq:def3}, then $\bmf{u}\equiv \bmf{ 0}$.
\end{thm}

\begin{proof}
Suppose that there exits a singular impedance line $\Gamma_h$ of $\bmf{u}$ as described at the beginning of this section with a nonzero constant impedance $\eta$. From \eqref{eq:def3}, we see that
\begin{equation}\label{eq:thm 33 cond1}
	\mathbf{u}(\bmf{0})=\bmf{0}
\mbox{ and } \partial_2  u_1(\bmf{0})=0.
\end{equation}
Moreover, $(T_\nu \bmf{u}+\eta \bmf{u}) \big|_{\Gamma_h^-}=\bmf{0}$ . Hence from Lemma \ref{lem:impedance eqn}, we know that \eqref{eq:Tu+u1} and \eqref{eq:Tu+u2} hold.

By virtue of \eqref{eq:thm 33 cond1}, we know that \eqref{eq:u0} and \eqref{eq:u1} hold. Substituting $a_1=0$ in \eqref{eq:Tu+u2} into \eqref{eq:u0} we have  $b_1=0$. Therefore from \eqref{eq:Tu+u1} we derive that \begin{equation}\label{eq:339}
	-\mathrm{i} k^2_p a_2 + k^2_2 b_2=0
\end{equation} 
by using the fact that $a_1=b_1=0$. Substituting \eqref{eq:339} into \eqref{eq:u1} we have $b_0=0$.

 Therefore,  from Lemma \ref{lem:impedance eqn}, we have $a_\ell=b_\ell=0$ for $\forall \ell \in \mathbb N \cup \{0\}$. Then $\bmf{u}\equiv \bmf{0}$ in $\Omega$ via Proposition \ref{prop:21}.

The proof is complete.
\end{proof}

\section{Generalized Holmgre's principle with the non-degenerate intersection of two homogeneous lines}

In this section, we consider the homogeneous lines introduced in Definition~\ref{def:1}. We shall show that the generic non-degenerate intersections of two of such lines also generate microlocal singularities, which prevent the occurrence of such intersections unless the Lam\'e eigenfunction $\bmf{u}$ is identically vanishing. As discussed in the beginning of Section~\ref{sect:2}, we assume throughout this section that the aforementioned two homogeneous lines are given by $\Gamma_h^\pm$ in \eqref{eq:gamma_pm} with the intersecting angle $\varphi_0\in (0, \pi)$.

\begin{lem}\label{lem:FB exp}
Let $\mathbf{u}$ be a solution to \eqref{eq:lame} with the radial wave expansion \eqref{eq:u} around the origin. Suppose that there are two rigid lines $\Gamma_h^+$ and $\Gamma_h^-$ in $\Omega$ of $\bmf{u}$ that are intersecting with each other in a non-degenerate way. Then
$$
- \bsi k_p ^4 (2 \lambda + \mu) a_0 + \mu k_s ^4 b_0=0.
$$
\end{lem}

\begin{proof}
Recall that $I_1^\pm$ are defined in \eqref{eq:CGO6}. Substituting \eqref{eq:Tudotv} into \eqref{eq:CGO6}, since $\bmf{u} \big|_{\Gamma_h^\pm}=\bmf{0}$, we can obtain that
\begin{equation}\label{eq:Tudotv1}
	\begin{split}
		 I_1^+ & =  \int_{\Gamma_h^+} \left({T}_{\nu} \bmf{u}\right) \cdot \bmf{v}\mathrm{d} \sigma  =
  - \int_0^h {\mathrm e}^{s \sqrt{ r} \zeta(\varphi_0 ) }\bigg\{\bsi k_p^2 (\lambda+\mu) \mathrm e^{\bsi \varphi_0} a_0 + \frac{\bsi}{2}  k_p^3 (\lambda+\mu) \mathrm e^{2 \bsi \varphi_0} a_1 r\\
&\quad +\frac{1}{8}\left[\bsi k_p^4 (\lambda+\mu) \mathrm e^{3 \bsi \varphi_0} a_2 - \bsi k_p^4 (2 \lambda+\mu) \mathrm e^{ \bsi \varphi_0} a_0+ k_s^4 \mu \mathrm e^{ \bsi \varphi_0} b_0\right]r^2 \bigg\}\rmd r - r_{I_1^+,1} \\
&=  -\bsi \left(\lambda+\mu\right) k_p^2 \left( 2 a_0 s^{-2} - 60 k_p^2 \mathrm{e}^{-2 \bsi \varphi_0} a_0 s^{-6} + 6 k_p a_1 s^{-4} + 30 k_p^2 a_2 s^{-6} \right)\\
 &-30 \bsi \mu k_p ^4  \mathrm{e}^{-2 \bsi \varphi_0} a_0 s^{-6} - 30\mu k_s^4   \mathrm{e}^{-2 \bsi \varphi_0} b_0 s^{-6}-r_{I_1^+,1},
	\end{split}
\end{equation}
where $r_{I_1^+,1}= - \int_0^h {\mathrm e}^{s \sqrt{ r} \zeta(\varphi_0 ) } R_{1,\Gamma_h^+} \rmd r $. Similarly, we have
\begin{equation}\label{eq:Tudotv2}
	\begin{split}
		 I_1^- & =  \int_{\Gamma_h^-} \left({T}_{\nu} \bmf{u}\right) \cdot \bmf{v}\mathrm{d} \sigma  =
  \int_0^h {\mathrm e}^{-s \sqrt{ r}}\bigg\{\bsi k_p^2 (\lambda+\mu)  a_0 + \frac{\bsi}{2}  k_p^3 (\lambda+\mu)  a_1 r\\
&\quad +\frac{1}{8}(\bsi k_p^4 (\lambda+\mu)  a_2 - \bsi k_p^4 (2 \lambda+\mu) a_0+ k_s^4 \mu  b_0)r^2 \bigg\}\rmd r + r_{I_1^-,1} \\
& =  \bsi \left(\lambda+\mu\right) k_p^2 \left( 2 a_0 s^{-2} - 60 k_p^2  a_0 s^{-6} + 6 k_p a_1 s^{-4} + 30 k_p^2 a_2 s^{-6} \right)\\
 &\quad +30\bsi \mu k_p ^4   a_0 s^{-6} + 30 \mu  k_s^4  b_0 s^{-6}+ r_{I_1^-,1},
	\end{split}
\end{equation}
where $r_{I_1^-,1}= \int_0^h {\mathrm e}^{- s \sqrt{ r}  } R_{1,\Gamma_h^-} \rmd r $. Therefore by straightforward calculations, we can derive that
\begin{equation}\label{eq:EFI}
	I_1^++I_1^-= 30 (1-\mathrm{e}^{-2\bsi \varphi_0})(-\bsi k_p ^4 (2 \lambda+\mu)a_0+ \mu k_s^4 b_0)s^{-6}+ r_{I_1^-,1}-r_{I_1^+,1}.
\end{equation}

By virtue of \eqref{eq:R1 s} and \eqref{eq:268b int}, it is easy to see that
\begin{equation}\label{eq:54 r}
\left|r_{I_1^+,1} \right|	\leq S_2 \cdot \Oh( s^{-8} ),\quad \left|r_{I_1^-,1} \right|	\leq S_2 \cdot \Oh( s^{-8} ).
\end{equation}
Substituting \eqref{eq:EFI} into \eqref{eq:CGO2}, we can obtain that
\begin{equation}
\notag
	\begin{split}
		30 (1-\mathrm{e}^{-2\bsi \varphi_0})(-\bsi k_p ^4 (2 \lambda+\mu)a_0+ \mu k_s^4 b_0)s^{-6}&=I_3-I_2+r_{I_1^-,1}-r_{I_1^+,1}, 
	\end{split}
\end{equation}
which further implies that
\begin{equation}\label{eq:55 int}
	\begin{split}
		30 (1-\mathrm{e}^{-2\bsi \varphi_0})(-\bsi k_p ^4 (2 \lambda+\mu)a_0+ \mu k_s^4 b_0)&=s^{6}\left(I_3-I_2+r_{I_1^-,1}-r_{I_1^+,1}\right).
	\end{split}
\end{equation}
From  \eqref{eq:55 int}, it yields that
\begin{equation}\label{eq:56 ineq}
	\left| 30 (1-\mathrm{e}^{-2\bsi \varphi_0})(-\bsi k_p ^4 (2 \lambda+\mu)a_0+ \mu k_s^4 b_0)  \right|\leq s^6\left(|I_3|+\left|I_2\right|+ \left| r_{I_1^-,1} \right|+\left|r_{I_1^+,1}\right|  \right).
\end{equation}
Since $\Gamma_h^-$ is a rigid line of $\bmf{u}$, then from Lemma \ref{lem:rigid}, \eqref{eq:2} holds. Substituting the second equation of \eqref{eq:2} into \eqref{eq:I3}, we can obtain that
\begin{equation}\label{eq:57 I_3}
	\left| I_3\right | \leq \frac{2\kappa \varphi_0 \Gamma(6) S_1 }{\delta_{\mathcal K }^8}s^{-8}.
\end{equation}

In \eqref{eq:56 ineq},  using \eqref{eq:I2I4}, \eqref{eq:54 r} and \eqref{eq:57 I_3}, letting $s\rightarrow +\infty$, under the condition $\varphi_0\neq \pi$, we can finish the proof of this lemma.
\end{proof}

\begin{thm}\label{thm:rigid line exp}
Let $\mathbf{u}\in L^2(\Omega)^2$ be a solution to \eqref{eq:lame}. There cannot exit two intersecting rigid lines $\Gamma_h^\pm$ of $\bmf{u}$ with the intersecting angle $\angle(\Gamma_h^+,\Gamma_h^-)=\varphi_0\neq \pi $ unless $\bmf{u}\equiv \bmf{0}$.
\end{thm}
\begin{proof} Suppose that there are two intersecting rigid lines $\Gamma_h^\pm$ of $\bmf{u}$ with the intersecting angle $\angle(\Gamma_h^+,\Gamma_h^- )=\varphi_0\neq \pi$. Recall that $\bmf{u}$ has the raidal wave expansion \eqref{eq:u}.
From \eqref{eq:2} and Lemma \ref{lem:FB exp},
it yields that
 \begin{equation}
\bigg\{
\begin{array}{l}
   -\bsi k_p^4(2\lambda+\mu)a_0+\mu k_s^4 b_0=0,\\
   k_p ^2 a_0- \bsi k_s^2 b_0=0,\\
\end{array}
\end{equation}
where $k_p$ and $k_s$ are defined in \eqref{eq:kpks}. Moreover, the eigenvalue $\kappa $  of \eqref{eq:lame} is positive. Hence by using \eqref{eq:convex} it is easy to  see that
$$
(2\lambda+\mu)k_p^2 +\mu k_s^2= \frac{3(\lambda+\mu  )\kappa }{\lambda+2\mu }>0. 
$$
Therefore
\begin{equation}\notag
\left|
\begin{array}{cc}
-\bsi k_p^4(2\lambda+\mu)    &   \mu k_s^4\\
                  k_p^2  &   -\bsi k_s^2\\
\end{array}
\right|=-k_p^2 k_s^2 \left[(2\lambda+\mu)k_p^2 +\mu k_s^2\right]<  0,
\end{equation}
 which implies that
 \begin{equation}\label{eq:a0b0=0}
 	 a_0=b_0=0.
 \end{equation}
Substituting \eqref{eq:a0b0=0} into \eqref{eq:u=0}, we compare the coefficients of $r^2$ in both sides of \eqref{eq:u=0} to obtain that
\begin{equation}\nonumber
\left(-2 k_p^3 a_1-2 \bsi k_s^3 b_1+k_p^3  a_3+\bsi k_s^3 b_3\right)\mathbf{e}_1-\left(k_p^3 a_1-\bsi k_s^3 b_1\right)\mathbf{e}_2=\bmf{0},
\end{equation}
which can be used to  further derive that
\begin{equation}\label{eq:4}
\left\{
    \begin{array}{l}
-2 k_p^3 a_1-2 \bsi k_s^3 b_1+k_p^3  a_3+\bsi k_s^3 b_3=0,\\
 k_p^3 a_1-\bsi k_s^3 b_1=0.
 \end{array}
 \right.
\end{equation}
Combining \eqref{eq:1} and \eqref{eq:4}, it yields that
\begin{equation}\label{eq:ll1}
k_p a_1 + \bsi k_s b_1=0,\ \
k_p^3 a_1-\bsi k_s^3 b_1=0.
\end{equation}
Since
\begin{equation}\notag
\left|
\begin{array}{cc}
k_p    &   \bsi k_s\\
k_p^3  &   -\bsi k_s^3\\
\end{array}
\right|=-\bsi k_p k_s \left(k_p^2 +k_s^2\right)\neq 0,
\end{equation}
which together with \eqref{eq:ll1} readily implies that
\begin{equation}\label{eq:a1b1=0}
	a_1=b_1=0.
\end{equation}
 In view of \eqref{eq:a0b0=0} and \eqref{eq:a1b1=0}, from Lemma \ref{lem:rigid} we obtain that $a_\ell=b_\ell=0$ for $\ell\in \mathbb N$, which implies that $\bmf{u} \equiv \bmf{0}$ in $\Omega$ by Proposition \ref{prop:21}.

 The proof is complete.
\end{proof}

\begin{lem}\label{lem:singular line}
Let $\mathbf{u}$ be a solution to \eqref{eq:lame} with the radial wave expansion \eqref{eq:u} around the origin. Suppose that there are two traction-free lines $\Gamma_h^+$ and $\Gamma_h^-$. If $\angle(\Gamma_h^+,\Gamma_h^-)= \varphi_0 \neq\pi$ and

\begin{equation}\label{eq:512 cond}
	 \frac{4  \varphi_0  }{3  \delta_{\mathcal K}^6 }<1,
\end{equation}
where $\delta_{\mathcal K}=\min_{ 0<\varphi<\varphi_0 }\cos(\varphi/2 )>0$ is defined in \eqref{eq:volume n}, then
$
b_0= 0$.
\end{lem}

\begin{proof}

Recall that $I_1^\pm$ are defined in \eqref{eq:CGO6}. Substituting \eqref{eq:Tvu} into \eqref{eq:CGO6}, since $T_\nu\bmf{u} \big|_{\Gamma_h^\pm}=\bmf{0}$, we can obtain that
\begin{equation}\label{eq:Tvdotu1}
	\begin{split}
		 I_1^+ & = - \int_{\Gamma_h^+} \left({T}_{\nu} \bmf{v}\right) \cdot \bmf{u}\mathrm{d} \sigma  =
  s \mu\zeta(\varphi_0) \int_0^h {\mathrm e}^{s \sqrt{ r} \zeta(\varphi_0 ) }\bigg\{\frac{1}{2}(\bsi k_p^2 a_0 + k_s^2 b_0) \mathrm e^{\bsi \varphi_0} r^{1/2}\\
  & + \frac{1}{8}(\bsi k_p^3 a_1 + k_s^3 b_1) \mathrm e^{2 \bsi \varphi_0} r^{3/2} + \frac{1}{48}(\bsi k_p^4 a_2 + k_s^4 b_2) \mathrm e^{3 \bsi \varphi_0} r^{5/2}\\
  & - \frac{1}{16}(\bsi k_p^4 a_0 + k_s^4 b_0) \mathrm e^{\bsi \varphi_0} r^{5/2}\bigg\}\rmd r - r_{I_1^+,2} \\
 & = -2 \mu(\bsi k_p^2 a_0+ k_s^2 b_0)s^{-2}-6 \mu(\bsi k_p^3 a_1+ k_s^3 b_1)s^{-4}-120 \mu(\bsi k_p^4 a_2+ k_s^4 b_2)s^{-6}\\
   & +90 \mathrm{e}^{-2 \bsi \varphi_0}\mu(\bsi k_p^4 a_0+ k_s^4 b_0)s^{-6}-r_{I_1^+,2},
	\end{split}
\end{equation}
where 
\begin{equation}\label{eq:r1+}
	r_{I_1^+,2}= \bsi s \mu\zeta(\varphi_0) \int_0^h {\mathrm e}^{s \sqrt{ r} \zeta(\varphi_0 ) } R_{2,\Gamma_h^+} \rmd r. 
\end{equation}
Similarly, we have
\begin{equation}\label{eq:Tvdotu2}
	\begin{split}
		 I_1^- & = - \int_{\Gamma_h^-} \left({T}_{\nu} \bmf{v}\right) \cdot \bmf{u}\mathrm{d} \sigma  =
  s \mu \int_0^h {\mathrm e}^{- s \sqrt{ r} }\bigg\{\frac{1}{2}(\bsi k_p^2 a_0 + k_s^2 b_0)  r^{1/2}\\
  & + \frac{1}{8}(\bsi k_p^3 a_1 + k_s^3 b_1)  r^{3/2} + \frac{1}{48}(\bsi k_p^4 a_2 + k_s^4 b_2)  r^{5/2}\\
  & - \frac{1}{16}(\bsi k_p^4 a_0 + k_s^4 b_0) r^{5/2}\bigg\}\rmd r - r_{I_1^-,2} \\
 & = 2 \mu(\bsi k_p^2 a_0+ k_s^2 b_0)s^{-2}+6 \mu(\bsi k_p^3 a_1+ k_s^3 b_1)s^{-4}+120 \mu(\bsi k_p^4 a_2+ k_s^4 b_2)s^{-6}\\
   & -90 \mu(\bsi k_p^4 a_0+ k_s^4 b_0)s^{-6} - r_{I_1^-,2},
	\end{split}
\end{equation}
where 
\begin{equation}\label{eq:r1-}
	r_{I_1^-,2} = \bsi s \mu\int_0^h {\mathrm e}^{- s \sqrt{ r}  } R_{2,\Gamma_h^-} \rmd r.
\end{equation} Therefore, from \eqref{eq:Tvdotu1}
 and \eqref{eq:Tvdotu2}, after straightforward algebraic calculations, we derive that
\begin{equation}
\notag
	I_1^++I_1^-= -90 \mu (1-\mathrm{e}^{-2\bsi \varphi_0})(\bsi k_p ^4 a_0+  k_s^4 b_0)s^{-6} - r_{I_1^-,2}-r_{I_1^+,2},
\end{equation}
which can be further reduced to
\begin{equation}\label{eq:EFI1}
	I_1^++I_1^-= -90 \mu (1-\mathrm{e}^{-2\bsi \varphi_0})  k_s^4 b_0s^{-6} - r_{I_1^-,2}-r_{I_1^+,2}
\end{equation}
via the second equality of \eqref{eq:8} since $\Gamma_h^-$ is a traction-free line.  By virtue of \eqref{eq:R2 s} and \eqref{eq:268b int}, it is easy to see that
\begin{equation}\label{eq:54 r1}
|r_{I_1^+,2}| \leq S_3 \cdot\Oh( s^{-8} ),|r_{I_1^-,2}| \leq S_3 \cdot\Oh( s^{-8} ),
\end{equation}
as $s\rightarrow +\infty$.

Using the fact that $a_0=0$ and \eqref{eq:I3}, we have as $s\rightarrow +\infty$ that
\begin{equation}\label{eq:55 I_3}
	\left| I_3 \right| \leq k_s^2 |b_0| \frac{\kappa \varphi_0 \Gamma(6) }{\delta_{\mathcal K}^6 }s^{-6}+ \frac{2 \kappa \varphi_0 \Gamma(8) }{\delta_{\mathcal K}^8 } s^{-8}.
\end{equation}
Substituting \eqref{eq:EFI1} into \eqref{eq:CGO2}, using \eqref{eq:I2I4},  \eqref{eq:54 r1} and \eqref{eq:55 I_3}, we obtain that
\begin{equation}\label{eq:518 in}
\begin{split}
	90 \mu |1-\mathrm{e}^{-2\bsi \varphi_0}|   k_s^4 |b_0| s^{-6} &\leq k_s^2 |b_0| \frac{\kappa \varphi_0 \Gamma(6) }{\delta_{\mathcal K}^6 }s^{-6}+ \frac{2 \kappa \varphi_0 \Gamma(8) }{\delta_{\mathcal K}^8 } s^{-8} + S_3\cdot \Oh(s^{-8})\\
	&+ \mathcal{C}_{{\mathcal K },B_h,\mu,\lambda}\|\bmf{u}\|_{H^2\left( {{\mathcal K}} \cap B_h \right)}\left(1+s\right){\mathrm e}^{-\delta_{{\mathcal K} } s \sqrt{h}}. 
	\end{split}
\end{equation}
From  \eqref{eq:512 cond}, recalling that \eqref{eq:kpks}, we can derive that
\begin{equation}\label{eq:512 cond con}
	 \mu     > \frac{4\kappa \varphi_0  }{3 k_s^2 \delta_{\mathcal K}^6 },
\end{equation}
Multiplying $s^6$ in both sides of \eqref{eq:518 in} and letting $s\rightarrow + \infty$, under the condition $\varphi_0\neq \pi$,  we can derive that
$$
 \mu     |b_0|\leq \frac{4\kappa \varphi_0  }{3 k_s^2 \delta_{\mathcal K}^6 } |b_0| 
$$
which implies that $b_0=0$  by virtue of \eqref{eq:512 cond con}. 

The proof is complete. 
\end{proof}

\begin{rem}\label{rem:varphi}
	Clearly when $\varphi_0 \in (0,\pi)$, it is easy to see that $\delta_{\mathcal K}=\min_{ 0<\varphi<\varphi_0 }\cos(\varphi/2 )= \cos(\varphi_0/2 ) >0$ and the function $f(\varphi_0)=\varphi_0/ \cos^6(\varphi_0/2 ) $ is a strictly  monotone increasing function. Denote
	$$
	g(\varphi_0 ):=\frac{4}{3}\cdot \frac{\varphi_0}{\cos^6(\varphi_0/2 ) }-1.
	$$
	Therefore $g(\varphi_0 )$ is a strictly  monotone increasing function. Let 
	\begin{equation}\label{eq:varphi0}
		\varphi_{\sf root}
	\end{equation}
 be the root of $g(\varphi_0 )$. In fact, using a standard root-finding algorithm, the numerical value of $\varphi_{\sf root} $ is $\varphi_{\sf root}^{\sf N}:=0.58043041944310849341051295527519$ and 
 $$
 g(\varphi_{\sf root}^{\sf N})=-5.5101297694794726936034525182293\times 10^{-40}. 
 $$
Hence if $\varphi_0 \in (0,  \varphi_{\sf root})$ we can claim that \eqref{eq:512 cond} is always fulfilled. 
\end{rem}

\begin{thm}\label{thm:two traction}
Let $\mathbf{u}\in L^2(\Omega)^2$ be a solution to \eqref{eq:lame}. Suppose there exit two intersecting lines $\Gamma_h^\pm$ of $\bmf{u}$ such that $\Gamma_h^\pm $ are two  traction-free lines  with the intersecting angle $\angle(\Gamma_h^+,\Gamma_h^-)=\varphi_0\neq \pi $, if $\bmf{u}(\bmf{0})=\bmf{0}$ and $\varphi_0 \in (0,  \varphi_{\sf root})$ where $\varphi_{\sf root}$ is defined in \eqref{eq:varphi0}, then $\bmf{u}\equiv \bmf{0}$.  
\end{thm}

\begin{proof}
	Suppose that there exit two intersecting lines $\Gamma_h^\pm$ of $\bmf{u}$ such that $\Gamma_h^\pm $ are two  traction-free lines satisfying that $\bmf{u}(\bmf{0})=\bmf{0}$ and  $0<\varphi_0<\varphi_{\sf root}$. Hence from Remark \ref{rem:varphi}, we know that \eqref{eq:512 cond} is fulfilled. Since $\Gamma_h^-$ is a traction free line, then \eqref{eq:8} and \eqref{eq:10} hold. Since \eqref{eq:512 cond} is fulfilled, from Lemma \ref{lem:singular line}, we know that $b_0=0$. Furthermore, under the condition $\bmf{u}(\bmf{0})=\bmf{0}$, we know \eqref{eq:u0} holds. Substituting $a_1=0$ of \eqref{eq:10} into \eqref{eq:u0}, one has $b_1=0$. In view of the second equality of  \eqref{eq:8}, we have shown that $a_\ell=b_\ell=0$ for $\ell=0,1$. Therefore from Lemma \ref{lem:tranction free} and Proposition \ref{prop:21}, we readily have $\bmf{u}\equiv \bmf{0}$ in $\Omega$, which finishes the proof. 
\end{proof}

\begin{thm}\label{thm:u&Tu exp}
Let $\mathbf{u}\in L^2(\Omega)^2$ be a solution to \eqref{eq:lame}. If there exist two intersecting lines $\Gamma_h^\pm$ of $\bmf{u}$ such that $\Gamma_h^-$ is a rigid line and $\Gamma_h^+$ is a traction-free line with the intersecting angle $\angle(\Gamma_h^+,\Gamma_h^-)=\varphi_0\neq \pi $, then $\bmf{u}\equiv \bmf{0}$. 
\end{thm}

\begin{proof}
Since $\Gamma_h^-$ is a rigid line and $\Gamma_h^+$ is a traction free line. Recall that $I_1^\pm$ are defined in \eqref{eq:CGO6}, using the definition of rigid and traction-free lines,  we can obtain that
\begin{equation}\label{eq:17a}
\begin{aligned}
 & I_1^-  =\int_{\Gamma_h^-} \left({T}_{\nu} \bmf{u}\right) \cdot \bmf{v}\mathrm{d} \sigma,\quad  I_1^+ = -\int_{\Gamma_h^+} \left({T}_{\nu} \bmf{v}\right) \cdot \bmf{u}\mathrm{d} \sigma.
\end{aligned}
\end{equation}
Using  \eqref{eq:Tudotv2} and \eqref{eq:Tvdotu1}, we can deduce that
\begin{equation}\label{eq:17b}
I_1^+ + I_1^-   = 2 ( \bsi \lambda k_p^2 a_0-  \mu k_s^2 b_0) s^{-2}+ R_{r,1} -r_{I_1^+,2}+ r_{I_1^-,1},
\end{equation}
where
\begin{equation}\label{eq:17e}
\begin{aligned}
R_{r,1} &  = 6 ( \bsi k_p^3 a_1- \mu k_s^3 b_1)s^{-4}+ 30 \bigg\{-\bsi k_p^4(\mu+2 \lambda-3 u \mathrm{e}^{-2\bsi \varphi_0})a_0\\
& + k_s^4 \mu (1+3\mathrm{e}^{-2\bsi \varphi_0})b_0 + \bsi k_p^4(\lambda-3\mu)a_2-4\mu k_s^4 b_2 \bigg\}s^{-6}. 
 \end{aligned}
\end{equation}
Therefore, from \eqref{eq:17e}, we have
\begin{equation}\label{eq:17f}
\left| R_{r,1} \right| \leq \Oh (s^{-4}),
\end{equation}
as $s\rightarrow +\infty$. 
Substituting \eqref{eq:17b} into \eqref{eq:CGO2}, we can obtain that
\begin{equation}\label{eq:17g}
	\begin{split}
		2 ( \bsi \lambda k_p^2 a_0-  \mu k_s^2 b_0) s^{-2}&=I_3-I_2-R_{r,1} -r_{I_1^-,1}+r_{I_1^+,2},
	\end{split}
\end{equation}
which can further yield that
\begin{equation}\label{eq:17h}
	\begin{split}
		2 ( \bsi \lambda k_p^2 a_0-  \mu k_s^2 b_0) &=s^{2} (I_3-I_2-R_{r,1} -r_{I_1^-,1}+r_{I_1^+,2})
	\end{split}
\end{equation}
From \eqref{eq:17h}, it gives that
\begin{equation}\label{eq:17i}
	\left| 2 ( \bsi \lambda k_p^2 a_0-  \mu k_s^2 b_0)  \right|\leq s^2\left(|I_3|+\left|I_2\right|+ \left|R_{r,1}\right|+\left| r_{I_1^-,1} \right|+\left|r_{I_1^+,2}\right|  \right).
\end{equation}
Since $\Gamma_h^-$ is a rigid line of $\bmf{u}$, then from Lemma \ref{lem:rigid}, we know \eqref{eq:2} holds. Substituting the second equation of \eqref{eq:2} into \eqref{eq:I3}, we can obtain that
\begin{equation}\label{eq:57 I_3b}
	\left| I_3\right | \leq \frac{2\kappa \varphi_0 \Gamma(6) S_1 }{\delta_{\mathcal K }^8}s^{-8}.
\end{equation}

In \eqref{eq:17i},  using \eqref{eq:I2I4}, \eqref{eq:54 r}, \eqref{eq:54 r1}, \eqref{eq:17f}  and \eqref{eq:57 I_3b}, letting $s\rightarrow +\infty$, under the condition $\varphi_0\neq \pi$,
 we can obtain that
\begin{equation}\label{eq:mix up 1}
\bsi \lambda k_p^2 a_0- \mu k_s^2 b_0=0.
\end{equation}
Combing \eqref{eq:2} with \eqref{eq:mix up 1}, it yields  that
\begin{equation}\nonumber
\left\{
\begin{array}{l}
k_p ^2 a_0 - \bsi k_s^2 b_0=0,\\
\bsi \lambda k_p^2 a_0- \mu k_s^2 b_0=0.
\end{array}
\right.
\end{equation}
In view of \eqref{eq:convex}, we have
\begin{equation}\notag
\left|
\begin{array}{cc}
k_p^2    &   -\bsi k_s^2\\
\bsi \lambda k_p^2  &   -\mu k_s^2
\end{array}
\right|=-(\lambda+\mu) k_p^2 k_s^2\neq 0,
\end{equation}
and therefore $
a_0=b_0=0.
$
Again using the fact that $\Gamma_h^-\in {\mathcal R}_\Omega^\kappa$, we have \eqref{eq:a1b1 lem}. From Lemma \ref{lem:rigid} and Proposition \ref{prop:21}, we have $\bmf{u}\equiv \bmf{0}$ in $\Omega$, which completes the proof. 
\end{proof}


\begin{thm}\label{thm:54}
Let $\mathbf{u}\in L^2(\Omega)^2$ be a solution to \eqref{eq:lame}. If there exist two intersecting lines $\Gamma_h^\pm$ of $\bmf{u}$ such that $\Gamma_h^-$ is a rigid line and $\Gamma_h^+$ is an impedance line with the intersecting angle $\angle(\Gamma_h^+,\Gamma_h^-)=\varphi_0\neq \pi $, where the associated impedance parameter is a nonzero constant $\eta_2\in \mathbb C$, then $\bmf{u}\equiv \bmf{0}$.

\end{thm}
\begin{proof}
Since $\Gamma_h^-$ is a rigid line and  $\Gamma_h^+$  is an impedance line, using the boundary conditions of $\bmf{u}$ on $\Gamma_h^\pm$ respectively, from the definition \eqref{eq:CGO6} of $I_1^\pm $, we have
\begin{equation}\notag
\begin{aligned}
 & I_1^-  =\int_{\Gamma_h^-} \left({T}_{\nu} \bmf{u}\right) \cdot \bmf{v}\mathrm{d} \sigma,\\
 & I_1^+ = \int_{\Gamma_h^+}\left[\left({T}_{\nu} \bmf{u}\right) \cdot \bmf{v}-\left({T}_{\nu} \bmf{v}\right) \cdot \bmf{u}\right] \mathrm{d} \sigma
         =-\int_{\Gamma_h^+}\left[\left(\eta_2 \bmf{u}\right) \cdot \bmf{v}+\left({T}_{\nu} \bmf{v}\right) \cdot \bmf{u}\right] \mathrm{d} \sigma. 
\end{aligned}
\end{equation}
By virtue of  \eqref{eq:Tudotv2} and \eqref{eq:Im+}, we know that
\begin{equation}\label{eq:IR}
\begin{aligned}
 I_1^- & =  \bsi \left(\lambda+\mu\right) k_p^2 \left( 2 a_0 s^{-2} - 60 k_p^2  a_0 s^{-6} + 6 k_p a_1 s^{-4} + 30 k_p^2 a_2 s^{-6} \right)\\
 & +30\bsi \mu k_p ^4   a_0 s^{-6} + 30 \mu  k_s^4  b_0 s^{-6}+ r_{I_1^-,1},\\
 I_1^+& = -2 \mu(\bsi k_p^2 a_0+ k_s^2 b_0)s^{-2}-6 \mu(\bsi k_p^3 a_1+ k_s^3 b_1)s^{-4}-120 \mu(\bsi k_p^4 a_2+ k_s^4 b_2)s^{-6}\\
   & +90 \mathrm{e}^{-2 \bsi \varphi_0}\mu(\bsi k_p^4 a_0+ k_s^4 b_0)s^{-6}-r_{I_1^+,2} + 6 \eta_2(k_p^2 a_0 - \bsi k_s^2 b_0) \mathrm e^{- \bsi \varphi_0} s^{-4} - r_{I_1^+,0},
\end{aligned}
\end{equation}
where $ r_{I_1^-,1}$, $r_{I_1^+,2} $, $r_{I_1^+,0}$ are defined in \eqref{eq:Tudotv2}, \eqref{eq:r1+} and \eqref{eq:Im+} respectively. 

Therefore, from \eqref{eq:IR},  after direct algebraic calculations, we can derive that
\begin{equation}\label{eq:IR1}
\begin{aligned}
I_1^- + I_1^+ & = (2 \bsi \lambda k_p^2 a_0-2 \mu k_s^2 b_0) s^{-2} + R_{r,2}- r_{I_1^+,0} + r_{I_1^-,1} -r_{I_1^+,2}
\end{aligned}
\end{equation}
where 
\begin{equation}\label{eq:r2}
\begin{aligned}
R_{r,2} & =  6 \big[\bsi \lambda k_p^3 a_1 - \mu k_s^3 b_1 + \eta_2 (k_p^2 a_0 - \bsi k_s^2 b_0) \mathrm{e}^{- \bsi \varphi_0}\big] s^{-4}\\
& + 30 \big[- \bsi (2 \lambda + \mu) k_p^4 a_0 + \mu k_s^4 b_0 + \bsi ( \lambda - 3  \mu) k_p^4 a_2 - 4 \mu k_s^4 b_2\\
& + 3 \mu (\bsi k_p^4 a_0 +  k_s^4 b_0) \mathrm{e}^{- 2 \bsi \varphi_0} \big] s^{-6}
\end{aligned}
\end{equation}

Substituting \eqref{eq:IR1} into \eqref{eq:CGO2},  one can deduce that
\begin{equation}\notag \label{eq:1a}
( 2 \bsi \lambda k_p^2 a_0-2 \mu k_s^2 b_0) s^{-2}= I_3 - I_2 - R_{r,2} + r_{I_1^+,0} - r_{I_1^-,1} + r_{I_1^+,2},
\end{equation}
which can be used to further derive that
\begin{equation}\label{eq:1b}
 2 \bsi \lambda k_p^2 a_0-2 \mu k_s^2 b_0 = s^2 (I_3 - I_2 - R_{r,2} + r_{I_1^+,0} - r_{I_1^-,1} + r_{I_1^+,2}). 
\end{equation}
From \eqref{eq:1b}, one can show that
\begin{equation}\label{eq:1c}
 \left|2 \bsi \lambda k_p^2 a_0-2 \mu k_s^2 b_0 \right|= s^2 (\left|I_3\right| +\left| I_2\right| +\left| R_{r,2}\right| + |r_{I_1^+,0}| + | r_{I_1^-,1}| + |r_{I_1^+,2}|). 
\end{equation}

Since $\Gamma_h^-$ is a rigid line of $\bmf{u}$, then from Lemma \ref{lem:rigid}, \eqref{eq:2} holds. Substituting the second equation of \eqref{eq:2} into \eqref{eq:I3}, we can obtain \eqref{eq:57 I_3b}.

In \eqref{eq:1c}, using \eqref{eq:I2I4}, \eqref{eq:54 r}, \eqref{eq:54 r1}, \eqref{eq:57 I_3b}, \eqref{eq:54 r3} and \eqref{eq:r2}, letting $s\rightarrow +\infty$, we can derive that
\begin{equation}\label{eq:mix up 2}
\bsi \lambda k_p^2 a_0- \mu k_s^2 b_0=0. 
\end{equation}
Recall that  $\Gamma_h^-$ is a rigid line of $\bmf{u}$.  Combining \eqref{eq:2} with \eqref{eq:mix up 2}, we have
\begin{equation}\nonumber
\left\{
\begin{array}{l}
k_p ^2 a_0 - \bsi k_s^2 b_0=0,\\
\bsi \lambda k_p^2 a_0- \mu k_s^2 b_0=0.
\end{array}
\right.
\end{equation}
By virtue of  \eqref{eq:convex},  since
\begin{equation}\notag
\left|
\begin{array}{cc}
k_p^2    &   -\bsi k_s^2\\
\bsi \lambda k_p^2  &   - \mu k_s^2
\end{array}
\right|=-(\lambda+\mu) k_p^2 k_s^2\neq 0,
\end{equation}
then $ a_0=b_0=0. $ Since  $\Gamma_h^-$ is a rigid line of $\bmf{u}$, from \eqref{eq:a1b1 lem}, we know that $a_1=b_1=0$. Therefore, from Lemma \ref{lem:rigid} and tje strong unique continuation principle, we have $\bmf{u}\equiv \bmf{0} $ in $\Omega$, which readily completes the proof. 
\end{proof}

\begin{thm}\label{thm:55}
Let $\mathbf{u}\in L^2(\Omega)^2$ be a solution to \eqref{eq:lame}. Suppose there exit two intersecting lines $\Gamma_h^\pm$ of $\bmf{u}$ such that $\Gamma_h^-$ is a traction-free line and $\Gamma_h^+$ is an impedance line associated with a nonzero impedance constant $\eta_2\in \mathbb C$, 
 with the property that $\angle(\Gamma_h^+,\Gamma_h^-)=\varphi_0\neq \pi $ and $\bmf{u}$ vanishes at the intersecting point, namely $\bmf{u}(\bmf{0})=\bmf{0}$, then $\bmf{u}\equiv \bmf{0}$. 
\end{thm}

\begin{proof}  Since $\Gamma_h^+ \in {\mathcal I}_\Omega^\kappa$, then $T_\nu \bmf{u}=-\eta_2 \bmf{u}$ on $\Gamma_h^+$.  Recall that  $I_1^+$ is defined in \eqref{eq:CGO6}. Therefore we have   \eqref{eq:Im+}. Recall that $\Gamma_h^- \in {\mathcal T}_\Omega^\kappa$ and $I_1^-$ is defined in \eqref{eq:CGO6}. We have  \eqref{eq:Tvdotu2}. Using \eqref{eq:Im+} and \eqref{eq:Tvdotu2}, we can derive that
%
%
%
%
%
%
\begin{equation}\label{eq:2a}
\begin{aligned}
I_1^- + I_1^+ & =6 \eta_2(k_p^2 a_0 - \bsi k_s^2 b_0) \mathrm e^{- \bsi \varphi_0} s^{-4}  -90 \mu (1-\mathrm{e}^{-2 \bsi \varphi_0})(\bsi k_p^4 a_0+ k_s^4 b_0)s^{-6}  \\
& - r_{I_1^+,2}- r_{I_1^-,2}  -r_{I_1^+,0}, 
\end{aligned}
\end{equation}
where $r_{I_1^+,2} $, $ r_{I_1^-,2} $ and $ r_{I_1^+,0}$ are defined in \eqref{eq:r1+},  \eqref{eq:r1-} and  \eqref{eq:Im+} respectively.  

Substituting \eqref{eq:2a} into \eqref{eq:CGO2},  one can deduce that
\begin{equation}\notag \label{eq:2b}
\begin{aligned}
6 \eta_2(k_p^2 a_0 - \bsi k_s^2 b_0) \mathrm e^{- \bsi \varphi_0} s^{-4} & = I_3 - I_2  + r_{I_1^+,0} + r_{I_1^+,2} + r_{I_1^-,2}
 + 90 \mu (1-\mathrm{e}^{-2 \bsi \varphi_0})(\bsi k_p^4 a_0+ k_s^4 b_0)s^{-6},
\end{aligned}
\end{equation}
which can further give that
\begin{equation}\label{eq:2c}
 \begin{aligned}
6 \eta_2(k_p^2 a_0 - \bsi k_s^2 b_0) \mathrm e^{- \bsi \varphi_0}  & = s^4\big(I_3 - I_2  + r_{I_1^+,0} + r_{I_1^+,2} + r_{I_1^-,2} + 90 \mu (1-\mathrm{e}^{-2 \bsi \varphi_0})(\bsi k_p^4 a_0+ k_s^4 b_0)s^{-6}\big). 
\end{aligned}
\end{equation}
From \eqref{eq:2c}, it yields that
\begin{equation}\label{eq:2d}\
\begin{aligned}
 \left|6 \eta_2(k_p^2 a_0 - \bsi k_s^2 b_0) \right| & \leq  s^4 \big(\left|I_3\right| +\left| I_2\right|  + |r_{I_1^+,0}| + | r_{I_1^+,2}| + |r_{I_1^-,2}|\\
 & +| 90 \mu (1-\mathrm{e}^{-2 \bsi \varphi_0})(\bsi k_p^4 a_0+ k_s^4 b_0)s^{-6} | \big),
 \end{aligned}
\end{equation}
In \eqref{eq:2d}, using \eqref{eq:I2I4}, \eqref{eq:I3}, \eqref{eq:54 r1} and \eqref{eq:54 r3} and letting $s\rightarrow +\infty$,
 we can obtain that
\begin{equation}\label{eq:18e}
k_p^2 a_0- \bsi k_s^2 b_0=0,
\end{equation}
Since $\Gamma_h^-$ is a traction-free line of $\bmf{u}$, from Lemma \ref{lem:tranction free}, we know that \eqref{eq:8} and \eqref{eq:10} hold.  Combining  \eqref{eq:8} with \eqref{eq:18e}, we can obtain that $a_0=b_0=0$.  
Since $\bmf{u}(\bmf{0)}=\bmf{0}$, we know that \eqref{eq:u0} holds.  Combining \eqref{eq:10} with \eqref{eq:u0} ,  it is easy to see that $a_1=b_1=0$. Therefore, from Lemma \ref{lem:tranction free} and Proposition \ref{prop:21}, we have
$
\bmf{u} \equiv \bmf{0} \mbox { in } \Omega,
$
which completes the proof. 
%
%
\end{proof}

\begin{lem}\label{lem:Im}
Let $\mathbf{u}\in L^2(\Omega)^2$ be a solution to \eqref{eq:lame}.  Suppose that there are two impedance lines $\Gamma_h^\pm$ intersecting with each other and satisfying 
\begin{equation*}
\begin{split}
	T_{\bmf{\nu}}\bmf{u}+\eta_1 \bmf{u}&=\bmf{0} \mbox{ on } \Gamma_h^-,\\
T_{\bmf{\nu}}\bmf{u}+\eta_2 \bmf{u}&=\bmf{0} \mbox{ on } \Gamma_h^+,
\end{split}
\end{equation*}
with $\eta_\ell \in \mathbb {C} \backslash \{0\}$.   If $\angle(\Gamma_h^+,\Gamma_h^-)=\varphi_0$ and
\begin{equation}\label{eq:lem53 cond}
	 \eta_2 \mathrm{e}^{-\bsi \varphi_0}+ \eta_1 \neq 0, 
\end{equation}
 then 
\begin{equation}\label{eq:lem53 eq1}
k_p^2 a_0-\bsi k_s^2 b_0=0. 
\end{equation}
Furthermore, if $a_0=b_0=\ldots=a_{\ell-1}=b_{\ell-1}=0$, under \eqref{eq:lem53 cond},  we have  
\begin{equation}\label{eq:lem53 eq2}
k_p^{\ell+2} a_{\ell}-\bsi k_s^{\ell+2} b_{\ell}=0.
\end{equation}
\end{lem}

\begin{proof}
 Since $\Gamma_h^+$ is an impedance line of $\bmf{u}$ with the impedance parameter $\eta_2$, then 
\begin{equation}\label{eq:528Tvu}
	T_\nu\bmf{u}=-\eta_2 \bmf{u}
\end{equation}
Recall that $I_1^\pm$ are defined in \eqref{eq:CGO6}. Using \eqref{eq:528Tvu}, substituting \eqref{eq:Tvu} into \eqref{eq:CGO6},  we can obtain that
 \begin{equation}\label{eq:Im+}
 \begin{aligned}
    I_1^+ &  = \int_{\Gamma_h^+}\left[\left({T}_{\nu} \bmf{u}\right) \cdot \bmf{v}-\left({T}_{\nu} \bmf{v}\right) \cdot \bmf{u}\right] \mathrm{d} \sigma
         =-\int_{\Gamma_h^+}\left[\left(\eta_2 \bmf{u}\right) \cdot \bmf{v}+\left({T}_{\nu} \bmf{v}\right) \cdot \bmf{u}\right] \mathrm{d} \sigma\\
        &  =\frac{\eta_2}{2} (k_p^2 a_0 - \bsi k_s^2 b_0) {\mathrm e}^{\bsi \varphi_0} \int_0^h {\mathrm e}^{s\sqrt{r} \zeta{(\varphi_0} )} r \rmd r- r_{I_1^+,0}\\
         &\quad  + s \mu\zeta(\varphi_0) \int_0^h {\mathrm e}^{s \sqrt{ r} \zeta(\varphi_0 ) }\bigg\{\frac{1}{2}(\bsi k_p^2 a_0 + k_s^2 b_0) \mathrm e^{\bsi \varphi_0} r^{1/2}\\
  & \quad + \frac{1}{8}(\bsi k_p^3 a_1 + k_s^3 b_1) \mathrm e^{2 \bsi \varphi_0} r^{3/2} + \frac{1}{48}(\bsi k_p^4 a_2 + k_s^4 b_2) \mathrm e^{3 \bsi \varphi_0} r^{5/2}\\
  &\quad  - \frac{1}{16}(\bsi k_p^4 a_0 + k_s^4 b_0) \mathrm e^{\bsi \varphi_0} r^{5/2}\bigg\}\rmd r - r_{I_1^+,2} \\
       & = -2 \mu(\bsi k_p^2 a_0+ k_s^2 b_0)s^{-2}-6 \mu(\bsi k_p^3 a_1+ k_s^3 b_1)s^{-4}-120 \mu(\bsi k_p^4 a_2+ k_s^4 b_2)s^{-6}\\
   & \quad +90 \mathrm{e}^{-2 \bsi \varphi_0}\mu(\bsi k_p^4 a_0+ k_s^4 b_0)s^{-6}-r_{I_1^+,2} + 6 \eta_2(k_p^2 a_0 - \bsi k_s^2 b_0) \mathrm e^{- \bsi \varphi_0} s^{-4} - r_{I_1^+,0}, 
 \end{aligned}
 \end{equation}
 where $r_{I_1^+,0}=\int_0^h {\mathrm e}^{s \sqrt{ r} \zeta(\varphi_0 ) } R_{0,\Gamma_h^+} \rmd r$ and $r_{I_1^+,2}$ is defined in \eqref{eq:r1+}. Similarly, we have
 \begin{equation}\label{eq:Im-}
 \begin{aligned}
    I_1^- &  = \int_{\Gamma_h^-}\left[\left({T}_{\nu} \bmf{u}\right) \cdot \bmf{v}-\left({T}_{\nu} \bmf{v}\right) \cdot \bmf{u}\right] \mathrm{d} \sigma
         =-\int_{\Gamma_h^-}\left[\left(\eta_1 \bmf{u}\right) \cdot \bmf{v}+\left({T}_{\nu} \bmf{v}\right) \cdot \bmf{u}\right] \mathrm{d} \sigma\\
        &  =\frac{\eta_1}{2} (k_p^2 a_0 - \bsi k_s^2 b_0)  \int_0^h {\mathrm e}^{ - s\sqrt{r}} r \rmd r
        + s \mu \int_0^h {\mathrm e}^{- s \sqrt{ r} }\bigg\{\frac{1}{2}(\bsi k_p^2 a_0 + k_s^2 b_0)  r^{1/2}\\
  &\quad  + \frac{1}{8}(\bsi k_p^3 a_1 + k_s^3 b_1)  r^{3/2} + \frac{1}{48}(\bsi k_p^4 a_2 + k_s^4 b_2)  r^{5/2}- \frac{1}{16}(\bsi k_p^4 a_0 + k_s^4 b_0) r^{5/2}\bigg\}\rmd r \\
  &\quad  - r_{I_1^-,0}- r_{I_1^-,2} \\
  & = 2 \mu(\bsi k_p^2 a_0+ k_s^2 b_0)s^{-2}+6 \mu(\bsi k_p^3 a_1+ k_s^3 b_1)s^{-4}+120 \mu(\bsi k_p^4 a_2+ k_s^4 b_2)s^{-6}\\
   &\quad  -90 \mu(\bsi k_p^4 a_0+ k_s^4 b_0)s^{-6} - r_{I_1^-,2} + 6 \eta_1(k_p^2 a_0 - \bsi k_s^2 b_0)  s^{-4} - r_{I_1^-,0}, 
 \end{aligned}
 \end{equation}
 where $r_{I_1^-,0}=\int_0^h {\mathrm e}^{- s \sqrt{ r}  } R_{0,\Gamma_h^-} \rmd r$ and $r_{I_1^-,2}$ is defined in \eqref{eq:r1-}.  Therefore, combining \eqref{eq:Im+} and \eqref{eq:Im-},   after direct algebraic calculations, we can derive that
\begin{equation}\label{eq:EFI2}
\begin{aligned}
	I_1^++I_1^- & =(6 \eta_2 \mathrm{e}^{-\bsi \varphi_0}+6 \eta_1) (k_p^2 a_0-\bsi k_s^2 b_0) s^{-4}- r_{I_1^-,0}-r_{I_1^+,0}\\
 & -90 \mu (1-\mathrm{e}^{-2\bsi \varphi_0})(\bsi k_p ^4 a_0+  k_s^4 b_0)s^{-6} - r_{I_1^-,2}-r_{I_1^+,2}.
 \end{aligned}
\end{equation}
By virtue of \eqref{eq:Ruv+} and \eqref{eq:R3 s}, it can be directly verified that
\begin{equation}\label{eq:54 r3}
|r_{I_1^+,0}| \leq S_0 \cdot\Oh( s^{-6} ),\ \ |r_{I_1^-,0}| \leq S_0 \cdot\Oh( s^{-6} ).
\end{equation}

 Substituting \eqref{eq:EFI2} into \eqref{eq:CGO2}, we can obtain that
\begin{equation}
\notag
	\begin{split}
		(6 \eta_2 \mathrm{e}^{-\bsi \varphi_0}+6 \eta_1) (k_p^2 a_0-\bsi k_s^2 b_0) s^{-4}&=I_3-I_2+r_{I_1^-,0}+r_{I_1^+,0}+r_{I_1^-,2}+r_{I_1^+,2}\\
                  & + 90 \mu (1-\mathrm{e}^{-2\bsi \varphi_0})(\bsi k_p ^4 a_0+  k_s^4 b_0)s^{-6}, 
	\end{split}
\end{equation}
which can further give that
\begin{equation}\label{eq:56 int}
	\begin{split}
		(6 \eta_2 \mathrm{e}^{-\bsi \varphi_0}+6 \eta_1) (k_p^2 a_0-\bsi k_s^2 b_0)&= s^4\bigg(I_3-I_2+r_{I_1^-,0}+r_{I_1^+,0}+r_{I_1^-,2}+r_{I_1^+,2}\\
                  & + 90 \mu (1-\mathrm{e}^{-2\bsi \varphi_0})(\bsi k_p ^4 a_0+  k_s^4 b_0)s^{-6}\bigg)
	\end{split}
\end{equation}
From  \eqref{eq:56 int}, it yields that
\begin{equation}\label{eq:56 ineq}
\begin{aligned}
	\left| 6 \eta_2 \mathrm{e}^{-\bsi \varphi_0}+6 \eta_1) (k_p^2 a_0-\bsi k_s^2 b_0)  \right| & \leq s^4 \bigg(|I_3|+\left|I_2\right|+ \left| r_{I_1^-,0} \right|+\left|r_{I_1^+,0}\right| + \left| r_{I_1^-,2} \right|+\left|r_{I_1^+,2}\right|\\
 & + \left| 90 \mu (1-\mathrm{e}^{-2\bsi \varphi_0})(\bsi k_p ^4 a_0+  k_s^4 b_0)s^{-6}\right| \bigg).
 \end{aligned}
\end{equation}
In \eqref{eq:56 ineq}, using \eqref{eq:I3}, \eqref{eq:arc}, \eqref{eq:54 r1} and \eqref{eq:54 r3}, and letting $s\rightarrow +\infty$,  under \eqref{eq:lem53 cond},  we can obtain that \eqref{eq:lem53 eq1}. 


Suppose that $a_0=b_0=\ldots=a_{\ell-1}=b_{\ell-1}=0$, by direct calculations and using  \eqref{eq:Tvu guina1new} and \eqref{eq:Ruv guina}, we can derive that
\begin{equation}\label{eq:EFI3}
\begin{aligned}
	I_1^++I_1^- & =\frac{(2 \ell + 3)!}{2^\ell (\ell+1)!}( \eta_2 \mathrm{e}^{-\bsi \varphi_0}+ \eta_1) (k_p^{\ell+2} a_\ell-\bsi k_s^{\ell+2} b_\ell) s^{-2 \ell-4}- \hat{r}_{I_1^-,0}-\hat{r}_{I_1^+,0}\\
 & -\frac{(2 \ell+6)!}{2^{\ell+2}(\ell+2)!} \mu (1-\mathrm{e}^{-2\bsi \varphi_0})(\bsi k_p ^{\ell+4} a_\ell +  k_s^{\ell+4} b_\ell)s^{-2 \ell-6} - \hat{r}_{I_1^-,2}-\hat{r}_{I_1^+,2}, 
 \end{aligned}
\end{equation}
where
\begin{equation*}
\begin{aligned}
& \hat{r}_{I_1^+,0}=\int_0^h \mathrm{e}^{s \sqrt{r} \zeta(\varphi_0)} \hat{R}_{0,\Gamma_h^+} \rmd r,\quad \hat{r}_{I_1^-,0}=\int_0^h \mathrm{e}^{ - s \sqrt{r}} \hat{R}_{0,\Gamma_h^-} \rmd r,\\
 & \hat{r}_{I_1^+,2}=\bsi s \mu \zeta(\varphi_0) \int_0^h \mathrm{e}^{s \sqrt{r} \zeta(\varphi_0)} \hat{R}_{2,\Gamma_h^+} \rmd r ,\quad  \hat{r}_{I_1^-,2}=\bsi s \mu  \int_0^h \mathrm{e}^{ - s \sqrt{r} } \hat{R}_{2,\Gamma_h^-} \rmd r.
\end{aligned}
\end{equation*}
Here  $\hat{R}_{2,\Gamma_h^+}$, $\hat{R}_{2,\Gamma_h^-}$, $\hat{R}_{0,\Gamma_h^+}$ and $\hat{R}_{0,\Gamma_h^-}$ are defined in \eqref{eq:RTvu++}, \eqref{eq:RTvu--}, \eqref{eq:Ruv++} and \eqref{eq:Ruv--} respectively.

From \eqref{eq:R4 s} and \eqref{eq:R5 s},  using \eqref{eq:268b int},  it is readily seen that
\begin{equation}\label{eq:54 r4}
\begin{aligned}
 & |\hat{r}_{I_1^+,0}| \leq \hat{S}_0 \cdot\Oh( s^{-2\ell-6} ),\quad  |\hat{r}_{I_1^-,0}| \leq \hat{S}_0 \cdot\Oh( s^{-2\ell-6} ),\\
  & |\hat{r}_{I_1^+,2}| \leq \hat{S}_3 \cdot\Oh( s^{-2\ell-6} ),\quad |\hat{r}_{I_1^-,2}| \leq \hat{S}_3 \cdot\Oh( s^{-2\ell-6} ), 
\end{aligned}
\end{equation}
where $\hat{S}_0$ and $\hat{S}_3 $ are defined in \eqref{eq:RTvu+Shat} and \eqref{eq:uvShat} respectively.  Substituting \eqref{eq:EFI3} into \eqref{eq:CGO2}, we can obtain that
\begin{equation}
\notag
	\begin{split}
		& \frac{(2 \ell + 3)!}{2^\ell (\ell+1)!}( \eta_2 \mathrm{e}^{-\bsi \varphi_0}+ \eta_1) (k_p^{\ell+2} a_\ell-\bsi k_s^{\ell+2} b_\ell) s^{-2 \ell-4} \\
& =I_3-I_2+\hat{r}_{I_1^-,0}+\hat{r}_{I_1^+,0}+\hat{r}_{I_1^-,2}+\hat{r}_{I_1^+,2}\\
                  & \quad + \frac{(2 \ell+6)!}{2^{\ell+2}(\ell+2)!} \mu (1-\mathrm{e}^{-2\bsi \varphi_0})(\bsi k_p ^{\ell+4} a_\ell +  k_s^{\ell+4} b_\ell)s^{-2 \ell-6}
	\end{split}
\end{equation}
which can be further reduced to
\begin{equation}\label{eq:57 int}
	\begin{split}
		& \frac{(2 \ell + 3)!}{2^\ell (\ell+1)!}( \eta_2 \mathrm{e}^{-\bsi \varphi_0}+ \eta_1) (k_p^{\ell+2} a_\ell-\bsi k_s^{\ell+2} b_\ell)  \\
& =s^{2\ell+4}\bigg(I_3-I_2+\hat{r}_{I_1^-,0}+\hat{r}_{I_1^+,0}+\hat{r}_{I_1^-,2}+\hat{r}_{I_1^+,2}\\
                  & \quad + \frac{(2 \ell+6)!}{2^{\ell+2}(\ell+2)!} \mu (1-\mathrm{e}^{-2\bsi \varphi_0})(\bsi k_p ^{\ell+4} a_\ell +  k_s^{\ell+4} b_\ell)s^{-2 \ell-6}\bigg),
	\end{split}
\end{equation}

From  \eqref{eq:57 int}, we can readily have
\begin{equation}\label{eq:57 ineq}
\begin{aligned}
	& \left| \frac{(2 \ell + 3)!}{2^\ell (\ell+1)!}( \eta_2 \mathrm{e}^{-\bsi \varphi_0}+ \eta_1) (k_p^{\ell+2} a_\ell-\bsi k_s^{\ell+2} b_\ell) \right|\\ & \leq s^{2\ell+4} \bigg(|I_3|+\left|I_2\right|+ \left| \hat{r}_{I_1^-,0} \right|+\left|\hat{r}_{I_1^+,0}\right| + \left| \hat{r}_{I_1^-,2} \right|+\left|\hat{r}_{I_1^+,2}\right|\\
 & + \left| \frac{(2 \ell+6)!}{2^{\ell+2}(\ell+2)!} \mu (1-\mathrm{e}^{-2\bsi \varphi_0})(\bsi k_p ^{\ell+4} a_\ell +  k_s^{\ell+4} b_\ell)s^{-2 \ell-6}\right| \bigg).
 \end{aligned}
\end{equation}
In \eqref{eq:57 ineq}, using \eqref{eq:I3+}, \eqref{eq:arc} and \eqref{eq:54 r4},  and letting $s\rightarrow +\infty$,  one finally sees that \eqref{eq:lem53 eq2} hold under \eqref{eq:lem53 cond}, which completes the proof.
\end{proof}

\begin{rem}\label{rem:imp}
	If the impedance parameters $\eta_1=\eta_2$, one can directly verify that \eqref{eq:lem53 cond} is equivalent to $\varphi_0\neq \pi$. 
\end{rem}

\begin{thm}\label{eq:impedance line exp}
Let $\mathbf{u}\in L^2(\Omega)^2$ be a solution to \eqref{eq:lame}. Suppose there exist two impedance lines  $\Gamma_h^\pm$ of $\bmf{u}$ such that $\angle(\Gamma_h^+,\Gamma_h^-)=\varphi_0\neq \pi$ satisfying that \eqref{eq:lem53 cond} is fulfilled and $\bmf{u}$ vanishes at the intersecting point, i.e. $\bmf{u}(\bmf{0})=\bmf{0}$, then $\bmf{u}\equiv \bmf{0}$. Here, $T_\nu \bmf{u}+\eta_2\bmf{u}=\mathbf{0}$ on $\Gamma_h^+$, $T_\nu \bmf{u}+\eta_1\bmf{u}=\mathbf{0}$ on $\Gamma_h^-$  and $\eta_\ell\in \mathbb{C}\backslash\{0\}, \ell=1,2.$.

\end{thm}

\begin{proof}
 Since that $\Gamma_h^-$ is an impedance line of $\bmf{u}$. From Lemma \ref{lem:impedance eqn}, we know that \eqref{eq:Tu+u1} and \eqref{eq:Tu+u2} hold. Under the condition \eqref{eq:lem53 cond}, from Lemma \ref{lem:Im}, we have \eqref{eq:lem53 eq1}. Combining \eqref{eq:Tu+u1} with \eqref{eq:lem53 eq1}, we can derive that $a_0=b_0=0$. Therefore, using Lemma \ref{lem:Im} again, under \eqref{eq:lem53 cond},  from \eqref{eq:lem53 eq2}, we have
\begin{equation}\label{eq:26C}
k_p^3 a_1 -\bsi k_s^3 b_1=0. 
\end{equation}
 Since $\bmf{u}(\bf{0})=\bmf{0}$, from Lemma \ref{lem:34 three conds},  \eqref{eq:u0} holds. Combining \eqref{eq:u0} and \eqref{eq:26C}, one readily has that
 \begin{equation}\notag
\left\{
\begin{array}{l}
k_p^3 a_1 -\bsi k_s^3 b_1=0,\\
k_p a_1 + \bsi k_s b_1=0. 
\end{array}
\right.
\end{equation}
Since 
\begin{equation}\notag
\left|
\begin{array}{cc}
k_p^3    &   -\bsi k_s^3\\
k_p  &   \bsi k_s\\
\end{array}
\right|=\bsi k_p k_s( k_p^2 + k_s^2)\neq 0, 
\end{equation}
we can derive that $a_1=b_1=0$. Using Lemma \ref{lem:impedance eqn} and Proposition \ref{prop:21}, we readily have $\bmf{u}\equiv \bmf{0}$.

The proof is complete. 
\end{proof}

\section{Unique identifiability for inverse elastic obstacle problems}\label{sect:6}

In this section, as an important and practical application, we apply the theoretical findings in the previous sections to the study of the unique identifiability for the inverse elastic  obstacle problem. As mentioned earlier, we shall refer to the theoretical findings in the previous sections as the generalized Holmgren's principle. The inverse problem is concerned with recovering the geometrical shape of a certain unknown object by using the elastic wave probing data. The inverse elastic obstacle problem arises from industrial applications of practical importance, e.g. in the geophysical exploration. We next introduce the mathematical setup of the inverse obstacle problem that expands the abstract formation \eqref{eq:ipa1}. 

Let $\Omega\subset\mathbb{R}^2$ be a bounded Lipschitz domain such that $\mathbb{R}^2\backslash\bar{\Omega}$ is connected. Let $\bmf{u}^i$ be an incident elastic wave field, which is a time-harmonic elastic plane wave of the form
\begin{equation}\label{eq:ui}
\bmf{u}^i:=\bmf{u}^i(\mathbf{ x};k_p,k_s, \mathbf{d})=\alpha_{p} \bmf{d} \mathrm{e}^{\mathrm{i} k_{p} \bmf{x} \cdot \bmf{d} }+\alpha_{s} \bmf{d}^{\perp} \mathrm{e}^{\mathrm{i} k_s  \bmf{x} \cdot \bmf{d} }, \quad \alpha_{p}, \alpha_{s} \in \mathbb{C}, \quad\left|\alpha_{p}\right|+\left|\alpha_{s}\right| \neq 0
\end{equation}
where  $\mathbf{d}\in\mathbb{S}^1$ denotes the incident direction, $\bmf{d}^{\perp}$ is orthogonal to $\bmf{d}$, $k_p$ and $k_s$ are compressional and shear wave numbers  defined in \eqref{eq:kpks}. Physically speaking, $\bmf{u}^i$ is the detecting wave field and $\Omega$ denotes an impenetrable obstacle which interrupts the propagation of the incident wave and generates the corresponding scattered wave field $\bmf{u}^{\mathrm{sc}}$. The scattered field $\bmf{u}^{\mathrm{sc} }$ in ${\mathbb R}^2 \backslash \Omega$ can be decomposed into
the sum of the compressional  part $\bmf{u}^{\mathrm{sc}}_p$ and the shear part $\bmf{u}^{\mathrm{sc}}_s$  as follows
\begin{equation}
	\bmf{u}^{\mathrm{s c} }=\bmf{u}_{p}^{\mathrm {s c} }+\bmf{u}_{s}^{\mathrm{s c} }, \quad \bmf{u}_{p}^{\mathrm {s c} }=-\frac{1}{k_{p}^{2}} \nabla\left( \nabla \cdot \bmf{ u}^{\mathrm {s c} }\right ), \quad\bmf{ u}_{s}^{\mathrm{s c}}=\frac{1}{k_{s}^{2}} \bf{curl} \operatorname{curl} u^{\mathrm {s c} },
\end{equation}
where
\[
 {\rm curl}\bmf{ u}=\partial_1 u_2-\partial_2 u_1, \quad {\bf
curl}{u}=(\partial_2 u, -\partial_1 u)^\top.
\]
Let $\omega=\sqrt{ \kappa}$ be the angular frequency., where $\kappa$ is the Lam\'e eigenvalue of \eqref{eq:lame}.  Define  $\bmf{u}:=\bmf{u}^i+\bmf{u}^{\mathrm{ sc} }$ to be the total wave field, then the forward scattering problem of this process can be described by the following system,
\begin{equation}\label{forward}
\begin{cases}
& {\mathcal L} \bmf{ u} + \omega^2  \bmf{u} = 0\qquad\quad \mbox{in }\ \ \mathbb{R}^2\backslash\overline{\Omega},\medskip\\
& \bmf{u} =\bmf{u}^i+\bmf{u}^{\mathrm{sc} }\hspace*{1.56cm}\mbox{in }\ \ \mathbb{R}^2,\medskip\\
& \mathscr{B}(\bmf{u})=\bmf{0}\hspace*{1.95cm}\mbox{on}\ \ \partial\Omega,\medskip\\
&\displaystyle{ \lim_{r\rightarrow\infty}r^{\frac{1}{2}}\left(\frac{\partial \bmf{u}_\beta^{\mathrm{sc} }}{\partial r}-\mathrm{i}k_\beta \bmf{u}_\beta^{\mathrm{sc} }\right) =\,0,} \quad \beta=p,s,
\end{cases}
\end{equation}
where the last equation is the Kupradze radiation condition that holds uniformly in $\hat{\mathbf{ x}}:=\mathbf{ x}/|\mathbf{ x}|\in\mathbb{S}^1$. The boundary condition $\mathscr{B}(u)$ on $\partial \Omega$ could be either of the following three conditions:
\begin{enumerate}
	\item the Dirichlet condition ($\Omega$ is a rigid obstacle): $\mathscr{B}(\bmf{u})=\bmf{u}$;
	\item the Neumann condition ($\Omega$ is a traction-free obstacle): $\mathscr{B}(\bmf{u})=T_\nu \bmf{u}$;
	\item the impedance condition ($\Omega$ is an impedance obstacle): $\mathscr{B}(\bmf{u})=T_\nu \bmf{u}+\eta \bmf{u},\ \Re(\eta)\geq 0 \mbox{ and } \Im(\eta)>0$, 
	\end{enumerate}
where $\nu$ denotes the exterior unit normal vector to $\partial\Omega$,  $\boldsymbol{\tau}= \nu^\perp$ and the boundary  traction operator $T_\nu$ is defined in \eqref{eq:Tu}. Moreover, in the impedance condition given above, $\eta\in L^\infty(\partial\Omega)$, and this is different from our study in the previous sections, where the impedance $\eta$ is always required to be a constant. We would also like to point out that the conditions $\Re(\eta)\geq 0 \mbox{ and } \Im(\eta)>0$ are the physical requirement. 


The elastic system \eqref{forward}  associated with either of the three kinds
of boundary conditions  is well understood with a unique
solution $ \bmf{u} \in H^1_{\mathrm{loc} } ({\mathbb R}^2 \backslash \Omega )$. We refer to \cite{ElschnerYama2010,Kupradze} for the related results. It is known that the compressional and shear parts $\bmf{u}_\beta^{\mathrm {sc} }$ ($\beta=p,s$)
of a radiating solution $\bmf{u}^{\mathrm{sc} }$ to the elastic system  \eqref{forward} possess the following asymptotic expansions 
\begin{equation}
	\begin{aligned}
	\bmf{u}_{p}^{\mathrm{sc}}(\bmf{x}; k_p,k_s, \mathbf{d}) &=\frac{\mathrm{e}^{\mathrm{i} k_{p} r}}{\sqrt{r}}\left\{u_{p}^{\infty}(\hat{\bmf{x}}; \bmf{d} ) \hat{\bmf{x}}+\Oh\left(\frac{1}{r}\right)\right\} \\
	\bmf{u}_{s}^{\mathrm{sc}}(\bmf{x};k_p,k_s, \mathbf{d}) &=\frac{\mathrm{e}^{\mathrm{i} k_{s} r}}{\sqrt{r}}\left\{u_{s}^{\infty}(\hat{\bmf{x}}; \bmf{d}) \hat{\bmf{x}}^{\perp}+\Oh\left(\frac{1}{r}\right)\right\}
	\end{aligned}
\end{equation}
as $r =|\bmf{x} | \rightarrow \infty$, where $u_{p}^{\infty}$ and $u_{s}^{\infty}$ are both scalar functions defined on $\mathbb S^1$. Hence, a Kupradze radiating solution has the asymptotic behavior
$$
\bmf{u}^{\mathrm{sc}}(\bmf{x}; k_p,k_s, \mathbf{d})=\frac{\mathrm{e}^{\mathrm{i} k_{p} r}}{\sqrt{r}} u_{p}^{\infty}(\hat{\bmf{x} }; \bmf{d} ) \hat{\bmf{x}}+\frac{\mathrm{e}^{\mathrm{i} k_{s} r}}{\sqrt{r}} u_{s}^{\infty}(\hat{\bmf{x} }; \bmf{d} ) \hat{\bmf{x} }^{\perp}+\Oh\left(\frac{1}{r^{3 / 2}}\right) \quad \text { as } \quad r \rightarrow \infty
$$
The far-field pattern $\bmf{u}^\infty$ of $\bmf{u}^{\mathrm{sc}} $ is defined as
$$
\bmf{u}_t^{\infty}(\hat{\bmf{x}}; \bmf{d} ) :=u_{p}^{\infty}(\hat{\bmf{x}}; \bmf{d} ) \hat{\bmf{x}}+u_{s}^{\infty}(\hat{\bmf{x}}; \bmf{d} ) \hat{\bmf{x}}^{\perp}.
$$
Obviously, the compressional and shear parts of the far-field are uniquely determined by $\bmf{u}^{\infty}$
as follows:
$$
\bmf{u}_{p}^{\infty}(\hat{\mathbf{x}}; \bmf{d} )=\bmf{u}^{\infty}(\hat{\bmf{x}};\bmf{d} ) \cdot \hat{\bmf{x}};  \quad \bmf{u}_{s}^{\infty}(\hat{\bmf{x}},\bmf{d} )=\bmf{u}^{\infty}(\hat{\bmf{x}}; \bmf{d} ) \cdot \hat{\bmf{x}}^{\perp}.
$$

The inverse elastic scattering problem corresponding to \eqref{forward} concerns the determination of the scatterer $\Omega$ (and $\eta$ as well in the impedance case) by knowledge of the far-field pattern $\bmf{u}_\beta^\infty(\hat{\mathbf{ x}},\mathbf{d},k)$, where $\beta=t,p$ or $s$. As in \eqref{eq:ipa1}, we introduce the operator $\mathcal{F}$ which sends the obstacle to the corresponding far-field pattern and is defined by the forward scattering system \eqref{forward}, the aforementioned inverse problem can be formulated as
\begin{equation}\label{inverse}
\mathcal{F}(\Omega, \eta)=\bmf{u}_\beta^\infty(\hat{\mathbf{ x}}; \mathbf{d}),\quad \beta=t,p, \mbox{ or } s.
\end{equation}

Next, we show that by using the generalized Holmgren's uniqueness principle, we can establish two novel unique identifiability results for \eqref{inverse} in determining an obstacle without knowing its a priori physical property as well as its possible surface impedance by at most four far-field patterns, namely $\bmf{u}_\beta^\infty(\hat{\mathbf{x}})$ corresponding to four distinct $\mathbf{d}$'s.  

\begin{defn}\label{def:61}
Let $Q\subset\mathbb{R}^2$ be a polygon in $\mathbb{R}^2$ such that
\begin{equation}\label{eq:edge}
\partial Q=\cup_{j=1}^\ell \Gamma_j,
\end{equation}
where each $\Gamma_j$ is an edge of $\partial Q$. $Q$ is said to be a generalized impedance obstacle associated with \eqref{forward} if there exists a Lipschitz dissection of $\Gamma_j$, $1\leq j\leq \ell$,
\[
\Gamma_j=\Gamma_D^j\cup\Gamma_N^j\cup\Gamma_I^j
\]
such that
\begin{equation}\label{eq:66}
\mathbf{ u} =\mathbf{0} \quad \text { on } \Gamma_{D}^j, \quad T_\nu \mathbf{u}=\mathbf{0}  \quad \text { on } \Gamma_{N}^j, \quad T_\nu \mathbf{u}+\eta_j   \mathbf{u}=\mathbf{0} \quad \text { on } \Gamma_{I}^j,
\end{equation}
where $\eta_j\in \mathbb{C}$ with $\Im\eta_j\geq 0$. 
\end{defn}
It is emphasized that in \eqref{eq:66}, either $\Gamma_D^j, \Gamma_N^j$ or $\Gamma_I^j$ could be an empty set, and hence a generalized impedance obstacle could be purely a rigid obstacle, a traction-free obstacle, an impedance obstacle or a mixed type. Moreover, one each edge of the polygonal obstacle, the impedance parameter can take different (complex) values. In order to simply notations, we formally write $T_\nu\bmf{u}+\eta\bmf{u}$ with $\eta\equiv \infty$ to signify $T_\nu\bmf{u}=\bmf{0}$. In doing so, \eqref{eq:66} can be unified as $T_\nu\bmf{u}+\eta\bmf{u}=\bmf{0}$ on $\partial\Omega$ with
\begin{equation}\label{eq:eta}
\eta=0\cdot\chi_{\cup_{j=1}^\ell \Gamma_D^j}+\infty\cdot\chi_{\cup_{j=1}^\ell \Gamma_N^j}+\sum_{j=1}^\ell \eta_j\cdot\chi_{\Gamma_I^j}.
\end{equation}
We write $(Q,\eta)$ to denote a generalized polygonal impedance obstacle as describe above with $\eta\in L^\infty(\partial Q)\cup\{\infty\}$.  
In what follows, $(\Omega , \eta)$ is said to be an admissible complex obstacle if 
\begin{equation}\label{eq:p1}
(\Omega, \eta)=\cup_{j=1}^p (\Omega_j, \eta_j), 
\end{equation}
where each $(\Omega_j, \eta_j)$ is a generalized polygonal impedance obstacle such that $\Omega_j, j=1, 2,\ldots, p$ are pairwise disjoint and 
\begin{equation}\label{eq:r2b}
\eta=\sum_{j=1}^p \eta_j\chi_{\partial\Omega_j},\quad \eta_j\in L^\infty(\partial\Omega_j)\cup\{\infty\}.  
\end{equation}

%

\begin{thm}\label{thm:uniqueness1}
Let $(\Omega, \eta)$ and $(\widetilde\Omega, \widetilde\eta)$ be two admissible complex obstacles. Let $\omega\in\mathbb{R}_+$ be fixed and $\mathbf{d}_\ell$, $\ell=1, 2,3, 4$ be four distinct incident directions from $\mathbb{S}^1$. Let $\bmf{u}_\beta^\infty$ and $\widetilde{\bmf{u}}^\infty_\beta$ be, respectively, the far-field patterns associated with $(\Omega, \eta)$ and $(\widetilde\Omega, \widetilde\eta)$, where $\beta=t,p, \mbox{ or } s$. If 
\begin{equation}\label{eq:cond1}
\bmf{u}_\beta^\infty(\hat{\mathbf{ x}}; \mathbf{d}_\ell )=\widetilde{\bmf{u}}_\beta^\infty(\hat{\mathbf{ x}}; \mathbf{d}_\ell), \ \ \hat{\mathbf x}\in\mathbb{S}^1, \ell=1, \ldots, 4, 
\end{equation}
then one has that
\begin{equation}\label{eq:u1n}
\Omega =\widetilde{\Omega}\mbox{ and } \eta=\widetilde \eta. 
\end{equation}
\end{thm}

Before giving the proof of Theorem \ref{thm:uniqueness1}, we first derive an auxiliary lemma as follows. 
\begin{lem}\label{lem:51}
	Let $\mathbf {\mathbf d}_{\ell}\in\mathbb{S}^1$, $\ell=1,\ldots, n$, be $n$ vectors 
	which are distinct from each other. Suppose that $\Omega$ is a bounded Lipschitz domain and $\mathbb R^2\backslash \overline{\Omega } $ is connected. Let the incident elastic wave filed $\bmf{u}^i(\mathbf{ x};k_p,k_s, \mathbf{d}_\ell)$   be defined in \eqref{eq:ui}. Furthermore, suppose that the total elastic wave filed $\bmf{u}(\bmf{x}; k_p,k_s, \mathbf{d}_\ell)$ associated with $\bmf{u}^i(\mathbf{ x};k_p,k_s, \mathbf{d}_\ell)$ 
 satisfies \eqref{forward}. Then {the following set of functions is linearly independent:}
	$$
	\{\bmf{u}(\bmf{x}; k_p,k_s, \mathbf{d}_\ell);~\mathbf x \in D , \ \ \ell=1,2,\ldots, n \},
	$$
	where $D \subset \mathbb R^2 \backslash \overline \Omega $ is an open set. 
\end{lem}

\begin{proof}

The lemma can be proved by following a similar argument to the proof of Theorem 5.1 in \cite{CK}, and we skip the detailed calculations. 

\end{proof}

\begin{proof}[Proof of Theorem~\ref{thm:uniqueness1}]
By an absurdity argument, we first prove that if \eqref{eq:cond1} holds, one must have that $\Omega=\widetilde \Omega$.  Suppose that $\Omega $ and $\widetilde{\Omega}$ are two different admissible complex obstacles such that $\Omega\neq \widetilde\Omega$ and \eqref{eq:cond1} holds. Let $\mathbf{G}$ denote the unbounded connected component of $\mathbb{R}^2\backslash\overline{(\Omega\cup\widetilde\Omega)}$. Then by a similar topological argument to that in \cite{Liu-Zou}, one can show that there exists a line segment $\Gamma_h\subset\partial\mathbf{G}\backslash\partial\Omega$ or $\Gamma_h\subset\partial\mathbf{G}\backslash\partial\widetilde\Omega$. Without loss of generality, we assume the former case. 

Let $\bmf{u}$ and $\widetilde{\bmf{u}}$ respectively denote the total wave fields to \eqref{forward} associated with $(\Omega, \eta)$ and $(\widetilde\Omega, \widetilde\eta)$. By \eqref{eq:cond1} and the Rellich theorem (cf. \cite{CK}), we know that
\begin{equation}\label{eq:aa3}
\bmf{u} (\mathbf x; k_p,k_s, \mathbf{d}_\ell)=\widetilde{\bmf{u} }(\mathbf x; k_p,k_s,\mathbf{d}_\ell),\quad \mathbf{x}\in\mathbf{G},\ \ell=1, \ldots, 4. 
\end{equation}
 By using \eqref{eq:aa3} as well as the generalized impedance boundary condition on $\partial\widetilde\Omega$, we readily have 
\begin{equation}\label{eq:aa4}
  T_\nu \bmf{u}+\widetilde\eta \bmf{u}=T_\nu \widetilde {\bmf{u}} +\widetilde\eta\widetilde{ \bmf{u}} =\bmf{0}\quad\mbox{on}\ \ \Gamma_{h}. 
\end{equation}
Consider a fixed point $\mathbf x_0 \in \Gamma_h $. There exits a sufficient small  positive number  $\varepsilon \in {\mathbb R}_+$ such that $B_{2\varepsilon}  (\mathbf x_0) \Subset \mathbf{G} $, where $B_{2\varepsilon }(\mathbf x_0)$ is a disk centered at $\mathbf x_0$ with the radius $2\varepsilon$. Let $\Gamma_\varepsilon =  B_\varepsilon (\mathbf x_0) \cap  \Gamma_h$, where $B_{\varepsilon }(\mathbf x_0)$ is a disk centered at $\mathbf{x_0}$ with the radius $\varepsilon$. It is also noted that 
$$
-{\mathcal L} \bmf{u}= \omega^2 \bmf{u} \mbox{ in } B_{2\varepsilon } (\mathbf x_0). 
$$
Recall that the unit normal vector $\nu $ and the tangential vector $\boldsymbol{\tau}$ to $\Gamma_h$ are defined in \eqref{eq:nutau}, respectively. Due to the linear dependence of four $\mathbb{C}^3$-vectors, it is easy to see that there exist four complex constants $a_\ell$  such that
\begin{equation}\notag
 \sum_{\ell=1}^4 a_{\ell } \begin{bmatrix}
	\bmf{u} (\mathbf x_0; k_p,k_s,  \mathbf{d}_\ell)\cr   \boldsymbol{\tau} \cdot \partial_{\nu} \bmf{u} |_{\bmf{x}={\bmf{x}}_0 }
\end{bmatrix}=\bmf{0}.
\end{equation}
Moreover, there exits at least one $a_\ell $ is not zero. 
Let
\begin{equation}\label{eq:u513}
 \bmf{u}(\bmf{x};k_p,k_s)=\sum_{\ell=1}^4 a_{\ell } \bmf{u} (\mathbf x; k_p,k_s, \mathbf{d}_\ell).
\end{equation}
Then we have
\begin{equation}\label{eq:611}
 \bmf{u}(\bmf{x}_0;k_p,k_s)=\bmf{0} \mbox{ and } \boldsymbol{\tau} \cdot \partial_{\nu} \bmf{u} |_{\bmf{x}={\bmf{x}}_0 }=0.
\end{equation}
Next we distinguish two separate cases. The first case is that $ \bmf{u}(\bmf{x};k_p,k_s)\equiv \bmf{0},\,  \forall \bmf{x} \in \mathbf{G}$. In view of \eqref{eq:u513}, since $a_\ell $ are not all zero and $\bmf{d}_\ell$ are distinct, we readily have a contradiction by Lemma \ref{lem:51}. For the second case,  we suppose that $ \bmf{u}(\bmf{x};k_p,k_s)\equiv\hspace*{-3mm}\backslash\  \bmf{0}$. 
  In view of \eqref{eq:aa4} and \eqref{eq:611}, recalling Definition \ref{def:generalized line}, we know that $\Gamma_\varepsilon $ is a singular  line of $\bmf{u}$, which implies that $\Gamma_\varepsilon  $ could be a singular rigid, or singular traction-free or singular impedance line of $\bmf{u}$ in Definition \ref{def:generalized line}. Therefore, by the generalized Holmgren's principle (cf. Theorems  \ref{Thm:31 singular rigid}, \ref{thm:traction free line} and \ref{thm:impedance line}), we obtain that
\begin{equation}\label{eq:613 con}
	\bmf{u}\equiv \bmf{0} \mbox{ in } B_{2\varepsilon }(\bmf{x}_c),
\end{equation}
which is obviously a contradiction.

Next, we prove that by knowing $\Omega =\widetilde{\Omega}$, one must have that $\eta \equiv  \widetilde \eta$.  Assume contrarily that $\eta \neq  \widetilde \eta$. One can easily show that there exists an open subset $\Sigma$ of $\partial \Omega=\partial \widetilde \Omega$ such that
$$
\mathbf{u}=T_\nu \mathbf{u}=\mathbf{0} \mbox { on } \Sigma.
$$
Therefore by the classical Holmgren's principle, we know that $\mathbf{u}\equiv \mathbf{0}$ in ${\mathbb R}^2 \backslash \Omega $, which readily yields a contradiction. 

The proof is complete. 
\end{proof}

\begin{rem}
	\eqref{eq:u1n} means that one can not only determine the shape of an admissible complex obstacle, but also its physical properties (in the case that $\eta=0$ or $\eta=\infty$). Furthermore, if it is of impedance type, one can determine the surface impedance parameter as well. 
\end{rem}

Finally, we show that if fewer far-field patterns are used, one can establish a local uniqueness result in determining a generic class of admissible complex obstacles. To that end, we first introduce a geometric notion of the degree of an admissible complex obstacle. Let $\Omega$ be defined in \eqref{eq:p1} that consists of finitely many pairwise disjoint polygons. Let $\Gamma, \Gamma'\subset\partial\Omega$ be two adjacent edges of $\partial\Omega$. Extending $\Gamma$ and $\Gamma'$ into straight lines in the plane $\mathbb{R}^2$, we denote them by $\widehat{\Gamma}$ and $\widehat{\Gamma'}$. Clearly, the intersection of $\widehat{\Gamma}$ and $\widehat{\Gamma'}$ forms two angles, with one belonging to $(0,\pi/2]$ and the other one belonging to $[\pi/2, \pi)$. We write $\angle_{\rm acute}(\Gamma,\Gamma')$ to signify the former one. Define 
\begin{equation}\label{de:dg1}
\mathrm{deg}(\Omega):=\max_{\Gamma, \Gamma'\in \partial\Omega}\{\angle_{\rm acute}(\Gamma, \Gamma')|\ \Gamma, \Gamma'\ \mbox{are two adjacent edges of}\ \partial\Omega\}. 
\end{equation}
Moreover, we let $\zeta$ and $\zeta'$ respectively signify the values of $\eta$ on $\Gamma$ and $\Gamma'$ around the vertex formed by those two edges. It is noted that $\zeta$ and $\zeta'$ may be $0, \infty$ or finite and nonzero. An admissible complex obstacle $(\Omega, \eta)$ is said to belong to the class $\mathcal{C}$ if 
\begin{equation}\label{eq:cond1n}
\mathrm{deg}(\Omega)<\varphi_{\sf root},
\end{equation} 
where $\varphi_{\sf root}$ is defined in \eqref{eq:varphi0}, and
\begin{equation}\label{eq:cond2n}
\zeta=\zeta'\quad \mbox{if both $\zeta$ and $\zeta'$ are finite and nonzero},
\end{equation}
for all two adjacent edges $\Gamma, \Gamma'$ of $\partial\Omega$.

\begin{thm}\label{thm:uniqueness2}
Let $(\Omega, \eta)$ and $(\widetilde\Omega, \widetilde\eta)$ be two admissible complex obstacles from the class $\mathcal{C}$ as described above.  Let $\omega \in\mathbb{R}_+$ be fixed and $\mathbf{d}_\ell$, $\ell=1, 2,3$ be three distinct incident directions from $\mathbb{S}^1$.  Let $\mathbf{G}$ denote the unbounded connected component of $\mathbb{R}^2\backslash\overline{(\Omega\cup\widetilde\Omega)}$. Let $\bmf{u}_\beta^\infty$ and $\widetilde{\bmf{u}}^\infty_\beta$ be, respectively, the far-field patterns associated with $(\Omega, \eta)$ and $(\widetilde\Omega, \widetilde\eta)$, where $\beta=t,p, \mbox{ or } s$. If 
\begin{equation}\label{eq:cond1 corner}
\bmf{u}_\beta^\infty(\hat{\mathbf{ x}},\mathbf{d}_\ell )=\widetilde{\bmf{u}}_\beta^\infty(\hat{\mathbf{ x}},\mathbf{d}_\ell), \ \ \hat{\mathbf x}\in\mathbb{S}^1, \ell=1, 2, 3, 
\end{equation}
then one has that
$$
\left(\partial \Omega \backslash \partial \overline{ \widetilde{\Omega }} \right  )\cup \left(\partial \widetilde{\Omega } \backslash \partial \overline{ \Omega } \right)
$$
cannot have a corner on  $\partial \mathbf{G}$.
\end{thm}

\begin{proof}
	We prove the theorem by contradiction. Assume \eqref{eq:cond1 corner}   holds but $
\left(\partial \Omega \backslash \partial \overline{ \widetilde{\Omega }} \right  )\cup \left(\partial \widetilde{\Omega } \backslash \partial \overline{ \Omega } \right)
$ has a corner $\mathbf x_c$ on $\partial \mathbf{G}$.  Clearly, $\mathbf x_c$ is either a vertex of $\Omega$ or a vertex of $\widetilde\Omega$. Without loss of generality, we assume the latter case. Let $h\in\mathbb{R}_+$ be sufficiently small such that $B_h(\bmf{x}_c)\Subset\mathbb{R}^2\backslash\overline \Omega $. Moreover, since $\bmf{x}_c$ is a vertex of $\widetilde\Omega$, we can assume that 
\begin{equation}\label{eq:aa2}
B_h(\mathbf x_c)\cap \partial\widetilde\Omega=\Gamma_h^\pm,
\end{equation}
where $\Gamma_h^\pm$ are the two line segments lying on the two edges of $\widetilde\Omega$ that intersect at $\mathbf x_c$.  Furthermore, on $\Gamma_h^\pm$ the boundary conditions are given by \eqref{eq:66}. 

By \eqref{eq:cond1 corner} and the Rellich theorem (cf. \cite{CK}), we know that
\begin{equation}\label{eq:aa5}
\bmf{u} (\mathbf x; k_p,k_s, \mathbf{d}_\ell)=\widetilde{\bmf{u} }(\mathbf x; k_p,k_s, \mathbf{d}_\ell),\quad x\in\mathbf{G},\ \ell=1, 2, 3. 
\end{equation}
It is clear that $\Gamma_h^\pm\subset\partial\mathbf{G}$. Hence, by using \eqref{eq:aa3} as well as the generalized boundary condition \eqref{eq:66} on $\partial\widetilde\Omega$, we readily have 
\begin{equation}\label{eq:aa4new}
\partial_\nu \bmf{u}+\widetilde\eta \bmf{u}=\partial_\nu \widetilde {\bmf{u}} +\widetilde\eta\widetilde{ \bmf{u}} =\bmf{0}\quad\mbox{on}\ \ \Gamma_h^\pm. 
\end{equation}
It is also noted that in $B_h(\mathbf x_c)$, $-{\mathcal L} \bmf{u}=\omega^2 \bmf{u}$.

Due to the linear dependence of three $\mathbb{C}^2$-vectors, we see that there exits three complex constants $a_\ell $ such that 
$$
\sum_{\ell=1}^3 a_\ell \bmf{u}(\bmf{x}_c;k_p,k_s,\bmf{d}_\ell)=\bmf{0},
$$
where there exits at least one $a_\ell $ is not zero. Set $\bmf{u}(\bmf{x};k_p,k_s)=\sum_{\ell=1}^3  a_\ell \bmf{u}(\bmf{x};k_p,k_s,\bmf{d}_\ell)$. Then we know that
\begin{equation}\label{eq:611new}
	 \bmf{u}(\bmf{x}_c;k_p,k_s)=\bmf{0}. 
\end{equation}
Similar to the proof of Theorem \ref{thm:uniqueness2}, we consider the following two cases. The first one is $ \bmf{u}(\bmf{x};k_p,k_s)\equiv \bmf{0},\, \forall \bmf{x}\in \bmf{G}$. Since there exits at one $a_\ell$ such that $a_\ell \neq 0$ and $\mathbf{d}_\ell$ are distinct, from Lemma \ref{lem:51}, we can arrive at a contradiction. The other case is that $ \bmf{u}(\bmf{x};k_p,k_s)\equiv\hspace*{-3mm}\backslash\  \bmf{0}$. By \eqref{eq:cond1n} and \eqref{eq:cond2n}, as well as the generalized Holmgren's principle in Theorems \ref{thm:rigid line exp}--\ref{eq:impedance line exp}, one can show that
$$
\mathbf{u} \equiv \mathbf{0} \mbox{ in } \mathbf{G}
$$
which yields a contradiction again. 

The proof is complete. 
\end{proof}

\begin{rem}
Following a similar argument, one can derive more unique identifiability results similar to Theorem~\ref{thm:uniqueness2}. For example, if one excludes the presence of $T_\nu \mathbf{u}=\mathbf{0}$ on any boundary portion in \eqref{eq:66} of Definition \ref{def:61}, then the assumption \eqref{eq:cond1n} in Theorem \ref{thm:uniqueness2} can be removed. We choose not to discuss the details about those extensions in this article. 
\end{rem}	


\section*{Acknowledgement}

The research of H Liu was supported by Hong Kong RGC GRF grants, No. 12302017, 12301218 and 12302919.

\section*{Appendix}

\begin{proof}[Proof of Lemma~\ref{lem:Tuu exp}] 
We first prove \eqref{eq:Tu1}. Recall that $\bmf{\nu }\big|_{\Gamma_h^+}$ is defined in \eqref{eq:nu}:
\begin{equation}\label{eq:normal der}
	\boldsymbol{\tau}=(-\cos\varphi_0,-\sin\varphi_0). 
\end{equation}
Substituting \eqref{eq:normal der} into \eqref{eq:Tu} yields
\begin{equation}\label{eq:Tu new exp}
\begin{aligned}
T_{\nu} \mathbf{u}\Big |_{\Gamma^+_h} & =2 \mu \left[\begin{array}{cc}
\partial_1 u_1 & \partial_2 u_1\\
\partial_1 u_2 & \partial_2 u_2\\
\end{array}\right]
\left[\begin{array}{c}
-\sin \varphi_0\\
\cos \varphi_0\\
\end{array}\right]
+\lambda \left[\begin{array}{c}
-\sin \varphi_0\\
\cos \varphi_0\\
\end{array}\right] \left(\partial_1 u_1+\partial_2 u_2\right)\\
& \quad + \mu
\left[\begin{array}{c}
-\cos \varphi_0\\
-\sin \varphi_0\\
\end{array}\right]\left(\partial_2 u_1-\partial_1 u_2\right):=\left[\begin{array}{c}{T_1 (\bmf{u}) }\big|_{\Gamma_h^+ }\\{T_2 (\bmf{u}) }\big|_{\Gamma_h^+ }\end{array}\right],
\end{aligned}
\end{equation}
where
\begin{equation*}
	\begin{split}
		T_1( \bmf{u})\big|_{\Gamma_h^+ }&=2\mu(-\sin \varphi_0 \partial_1 u_1+\cos \varphi_0  \partial_2 u_1 ) -\lambda \sin \varphi_0 (\partial_1 u_1+\partial_2 u_2 )-\mu \cos \varphi_0 (\partial_2 u_1-\partial_1 u_2  ) , \\
		T_2 (\bmf{u})\big|_{\Gamma_h^+ }&=2\mu(-\sin \varphi_0 \partial_1 u_2+\cos \varphi_0  \partial_2 u_2 ) -\lambda \cos  \varphi_0 (\partial_1 u_1+\partial_2 u_2 )-\mu \sin  \varphi_0 (\partial_2 u_1-\partial_1 u_2  ).
	\end{split}
\end{equation*}
Using \eqref{eq:u comp}, it is readily shown that
\begin{equation}\label{eq:u1 par}
\begin{aligned}
\frac{\partial u_1}{\partial r}
&= \sum_{m=0} ^\infty \left\{ \frac{k_p^2}{4}a_m \left\{\mathrm{e}^{\bsi \left(m-1\right) \varphi} J_{m-2}\left(k_p r\right) - \left[\mathrm{e}^{\bsi \left(m-1\right) \varphi}+\mathrm{e}^{\bsi \left(m+1\right) \varphi}\right] J_m \left(k_p r\right) \right\}\right.\\
&\quad  + \frac{\bsi k_s^2}{4}b_m \left\{\mathrm{e}^{\bsi \left(m-1\right) \varphi} J_{m-2}\left(k_s r\right) - \left[\mathrm{e}^{\bsi \left(m-1\right) \varphi}-\mathrm{e}^{\bsi \left(m+1\right) \varphi}\right] J_m \left(k_s r\right) \right\}\\
&\quad  + \frac{k_p^2}{4}a_m \mathrm{e}^{\bsi \left(m+1\right) \varphi} J_{m+2} \left(k_p r\right)-\frac{\bsi k_s^2}{4}b_m \mathrm{e}^{\bsi \left(m+1\right) \varphi} J_{m+2} \left(k_s r\right)
 \bigg\},\\
\frac{\partial u_1}{\partial \varphi}
 &= \sum_{m=0} ^\infty  \left\{\frac{\bsi\left(m-1\right)}{2} k_p \mathrm{e}^{\bsi \left(m-1\right) \varphi} J_{m-1} \left(k_p r\right) a_m - \frac{\bsi\left(m+1\right)}{2} k_p \mathrm{e}^{\bsi \left(m+1\right) \varphi} J_{m+1} \left(k_p r\right) a_m \right.\\
 & \quad - \frac{\left(m-1\right)}{2} k_s \mathrm{e}^{\bsi \left(m-1\right) \varphi} J_{m-1} \left(k_s r\right) b_m - \frac{\left(m+1\right)}{2} k_s \mathrm{e}^{\bsi \left(m+1\right) \varphi} J_{m+1} \left(k_s r\right) b_m\bigg\},
\end{aligned}
\end{equation}
and
\begin{equation}\label{eq:u2 par}
\begin{aligned}
\frac{\partial u_2}{\partial r}
& = \sum_{m=0} ^\infty \left\{\frac{\bsi k_p^2}{4}a_m \left\{\mathrm{e}^{\bsi \left(m-1\right) \varphi} J_{m-2}\left(k_p r\right) - \left[\mathrm{e}^{\bsi \left(m-1\right) \varphi}-\mathrm{e}^{\bsi \left(m+1\right) \varphi}\right] J_m \left(k_p r\right)\right\}\right.\\
&\quad  +  \frac{k_s^2}{4}b_m \left\{-\mathrm{e}^{\bsi \left(m-1\right) \varphi} J_{m-2}\left(k_s r\right) + \left[\mathrm{e}^{\bsi \left(m-1\right) \varphi}+\mathrm{e}^{\bsi \left(m+1\right) \varphi}\right] J_m \left(k_s r\right) \right\}\\
&\quad -  \frac{\bsi k_p^2}{4}a_m \mathrm{e}^{\bsi \left(m+1\right) \varphi} J_{m+2} \left(k_p r\right) -\frac{k_s^2}{4} b_m \mathrm{e}^{\bsi \left(m+1\right) \varphi} J_{m+2} \left(k_s r\right)
\bigg\},
\end{aligned}
\end{equation}
\begin{equation}
\begin{aligned}
 \frac{\partial u_2}{\partial \varphi}
&= \sum_{m=0} ^\infty  \left\{-\frac{\left(m-1\right)}{2} k_p \mathrm{e}^{\bsi \left(m-1\right) \varphi} J_{m-1} \left(k_p r\right) a_m - \frac{\left(m+1\right)}{2} k_p \mathrm{e}^{\bsi \left(m+1\right) \varphi} J_{m+1} \left(k_p r\right) a_m \right.\\
 & \quad - \frac{\bsi \left(m-1\right)}{2} k_s \mathrm{e}^{\bsi \left(m-1\right) \varphi} J_{m-1} \left(k_s r\right) b_m + \frac{\bsi \left(m+1\right)}{2} k_s \mathrm{e}^{\bsi \left(m+1\right) \varphi} J_{m+1} \left(k_s r\right) b_m\bigg\}.
\end{aligned}
\end{equation}
Using the fact that
\begin{equation}\label{eq:u3 par}
	\begin{split}
		\frac{\partial u_i}{\partial x_1}&=\cos\varphi \cdot \frac{\partial u_i}{\partial r}- \frac{\sin \varphi}{r} \cdot \frac{\partial u_i}{\partial \varphi},\quad
		\frac{\partial u_i}{\partial x_2}=\sin \varphi \cdot \frac{\partial u_i}{\partial r}+ \frac{\cos \varphi}{r} \cdot \frac{\partial u_i}{\partial \varphi},
	\end{split}
 i=1,2,
\end{equation}
as well as \eqref{eq:u1 par} and \eqref{eq:u2 par}, by tedious but straightforward calculations, one can obtain that
\begin{equation}\label{eq:u1 partial 127}
\begin{aligned}
& \partial_1 u_1 \cdot \left(-\sin \varphi_0\right)+\partial_2 u_1  \cdot \left(\cos \varphi_0\right)\\
& =  \sum_{m=0} ^\infty  \Bigg \{ \sin(\varphi -\varphi_0) \Big[\frac{k_p^2}{4} a_m \mathrm{e}^{\bsi \left(m-1\right) \varphi} J_{m-2}(k_p r)    - \frac{k_p^2}{4} a_m J_m(k_p r)\left(\mathrm{e}^{\bsi (m-1)\varphi}+\mathrm{e}^{\bsi (m+1)\varphi}\right)\\
&  + \frac{k_p^2}{4} a_m \mathrm{e}^{\bsi \left(m+1\right) \varphi} J_{m+2}(k_p r)  + \frac{\bsi k_s^2}{4} b_m \mathrm{e}^{\bsi \left(m-1\right) \varphi} J_{m-2}(k_s r)  \\
&  - \frac{\bsi k_s^2}{4} b_m J_m(k_s r)\left(\mathrm{e}^{\bsi (m-1)\varphi}-\mathrm{e}^{\bsi (m+1)\varphi}\right)  - \frac{\bsi k_s^2}{4} b_m \mathrm{e}^{\bsi \left(m+1\right) \varphi} J_{m+2}(k_s r) \Big ]  \\
& + \frac{ \cos(\varphi -\varphi) }{r} \Big [\frac{\bsi(m-1)k_p}{2} a_m \mathrm{e}^{\bsi (m-1)\varphi} J_{m-1}(k_p r) \\
&  - \frac{\bsi(m+1)k_p}{2} a_m \mathrm{e}^{\bsi (m+1)\varphi} J_{m+1}(k_p r) - \frac{(m-1)k_s}{2} b_m \mathrm{e}^{\bsi (m-1)\varphi} J_{m-1}(k_s r)\\
 & - \frac{(m+1)k_s}{2} b_m \mathrm{e}^{\bsi (m+1)\varphi} J_{m+1}(k_s r) \Big] \Bigg  \},\\
& \partial_1 u_2 \cdot \left(-\sin \varphi_0\right)+\partial_2 u_2 \cdot \left(\cos \varphi_0\right)\\
& =  \sum_{m=0} ^\infty \Bigg \{ \sin(\varphi -\varphi_0 ) \Big [ \frac{\bsi k_p^2}{4} a_m \mathrm{e}^{\bsi \left(m-1\right) \varphi} J_{m-2}(k_p r)   - \frac{\bsi k_p^2}{4} a_m J_m(k_p r)\left(\mathrm{e}^{\bsi (m-1)\varphi}-\mathrm{e}^{\bsi (m+1)\varphi}\right) \\
&  - \frac{\bsi k_p^2}{4} a_m \mathrm{e}^{\bsi \left(m+1\right) \varphi} J_{m+2}(k_p r)  - \frac{ k_s^2}{4} b_m \mathrm{e}^{\bsi \left(m-1\right) \varphi} J_{m-2}(k_s r) + \frac{ k_s^2}{4} b_m J_m(k_s r)\left(\mathrm{e}^{\bsi (m-1)\varphi}+\mathrm{e}^{\bsi (m+1)\varphi}\right)\\
& -\frac{ k_s^2}{4} b_m \mathrm{e}^{\bsi \left(m+1\right) \varphi} J_{m+2}(k_s r) \Big] + \frac{\cos(\varphi -\varphi_0)}{r} \Big[ \frac{-(m-1)k_p}{2} a_m \mathrm{e}^{\bsi (m-1)\varphi} J_{m-1}(k_p r) \\
& - \frac{(m+1)k_p}{2} a_m \mathrm{e}^{\bsi (m+1)\varphi} J_{m+1}(k_p r) - \frac{\bsi(m-1)k_s}{2} b_m \mathrm{e}^{\bsi (m-1)\varphi} J_{m-1}(k_s r)\\
 & + \frac{\bsi(m+1)k_s}{2} b_m \mathrm{e}^{\bsi (m+1)\varphi} J_{m+1}(k_s r) \Big] \Bigg \}.
\end{aligned}
\end{equation}
Similarly, from \eqref{eq:u1 par} and \eqref{eq:u2 par}, we have
\begin{equation}
\begin{aligned}
& \partial_1 u_1 +\partial_2 u_2  =  \sum_{m=0} ^\infty  \Bigg \{ \frac{k_p^2}{4} \mathrm{e}^{\bsi m \varphi} a_m \left( J_{m-2} \left(k_p r\right)-2 J_m \left(k_p r\right)
+ J_{m+2} \left(k_p r\right)\right ) \\
& +\frac{\bsi k_s^2}{4} \mathrm{e}^{\bsi m \varphi} b_m \left ( J_{m-2} \left(k_s r\right)-J_{m+2} \left(k_s r\right)\right) 
 + \frac{1}{r} \bigg [ -\frac{k_p}{2}\mathrm{e}^{\bsi m \varphi} a_m \big( \left(m-1\right)J_{m-1} \left(k_p r\right)\\
 &+\left(m+1\right)J_{m+1} \left(k_p r\right)\big)  - \frac{\bsi k_s}{2}\mathrm{e}^{\bsi m \varphi} b_m \left( \left(m-1\right)J_{m-1} \left(k_s r\right)-\left(m+1\right)J_{m+1} \left(k_s r\right)\right)   \bigg ]\Bigg\},
\end{aligned}
\end{equation}
and
\begin{equation}\label{eq:u1u2 129}
\begin{aligned}
& \partial_2 u_1-\partial_1 u_2
 =  \sum_{m=0} ^\infty  \Bigg \{\frac{\bsi k_p^2}{4} \mathrm{e}^{\bsi m \varphi} a_m \left( -J_{m-2} \left(k_p r\right)+J_{m+2} \left(k_p r\right)\right)  \\
 & +\frac{k_s^2}{4} \mathrm{e}^{\bsi m \varphi} b_m \left( J_{m-2} \left(k_s r\right)-2 J_m \left(k_s r\right)+J_{m+2} \left(k_s r\right)\right)  + \frac{1}{r} \bigg [ \frac{\bsi k_p}{2}\mathrm{e}^{\bsi m \varphi} a_m \bigg( \left(m-1\right)J_{m-1} \left(k_p r\right)\\
 & -\left(m+1\right)J_{m+1} \left(k_p r\right)\bigg) - \frac{k_s}{2}\mathrm{e}^{\bsi m \varphi} b_m \left( \left(m-1\right)J_{m-1} \left(k_s r\right)+\left(m+1\right)J_{m+1} \left(k_s r\right)\right)  \bigg] \Bigg\}.
\end{aligned}
\end{equation}
Plugging \eqref{eq:u1 partial 127}--\eqref{eq:u1u2 129} into \eqref{eq:Tu new exp}, after tedious but straightforward calculations, we have
\begin{equation}\label{eq:J1}
\begin{aligned}
&T_1(\bmf{u})\big |_{\Gamma_h^+  }  =  \sum_{m=0} ^\infty  \left\{\frac{k_p^2}{4} a_m \mathrm{e}^{\bsi \left(m-1\right) \varphi} J_{m-2}(k_p r) \left[2 \mu \sin(\varphi-\varphi_0)-\lambda \mathrm{e}^{\bsi \varphi} \sin \varphi_0 +\bsi \mu \mathrm{e}^{\bsi \varphi} \cos \varphi_0\right] \right.\\
& + \frac{k_p^2}{2} a_m \mathrm{e}^{\bsi m \varphi} J_m(k_p r)\left[\lambda \sin\varphi_0+2 \mu \cos\varphi \sin(\varphi_0-\varphi)\right]\\
& + \frac{k_p^2}{4} a_m \mathrm{e}^{\bsi m \varphi} J_{m+2}(k_p r) \left[2 \mu \mathrm{e}^{\bsi \varphi} \sin(\varphi-\varphi_0) - \lambda \sin\varphi_0 - \bsi \mu \cos \varphi_0 \right]\\
& + \frac{\bsi k_s^2}{4} b_m \mathrm{e}^{\bsi \left(m-1\right) \varphi} J_{m-2}(k_s r) \left[2 \mu \sin(\varphi-\varphi_0)-\lambda \mathrm{e}^{\bsi \varphi} \sin \varphi_0 +\bsi \mu \mathrm{e}^{\bsi \varphi} \cos \varphi_0\right]\\
& + \frac{ k_s^2}{2} b_m \mathrm{e}^{\bsi m\varphi} J_m(k_s r)\left[\mu \cos\varphi_0+2 \mu \sin\varphi \sin(\varphi_0-\varphi)\right] \\
& + \frac{\bsi k_s^2}{4} b_m \mathrm{e}^{\bsi m \varphi} J_{m+2}(k_s r) \left[2 \mu \mathrm{e}^{\bsi \varphi} \sin(\varphi_0-\varphi) + \lambda \sin\varphi_0 + \bsi \mu \cos \varphi_0\right]\\
& + \frac{1}{r} \left\{\frac{(m-1)k_p}{2} a_m \mathrm{e}^{\bsi (m-1)\varphi} J_{m-1}(k_p r)\left[2 \bsi \mu \cos(\varphi_0-\varphi) + \lambda \mathrm{e}^{\bsi \varphi} \sin\varphi_0 - \bsi \mu \mathrm{e}^{\bsi \varphi} \cos\varphi_0\right] \right.\\
& + \frac{(m+1)k_p}{2} a_m \mathrm{e}^{\bsi m \varphi} J_{m+1}(k_p r)\left[-2 \bsi \mu \mathrm{e}^{\bsi \varphi} \cos(\varphi_0-\varphi) + \lambda  \sin\varphi_0 + \bsi \mu  \cos\varphi_0\right]\\
& + \frac{(m-1)k_s}{2} b_m \mathrm{e}^{\bsi (m-1)\varphi} J_{m-1}(k_s r)\left[-2 \mu \cos(\varphi_0-\varphi) + \bsi \lambda \mathrm{e}^{\bsi \varphi} \sin\varphi_0 + \mu \mathrm{e}^{\bsi \varphi} \cos\varphi_0\right]\\
 & + \frac{(m+1)k_s}{2} b_m \mathrm{e}^{\bsi m\varphi} J_{m+1}(k_s r) \left[-2 \mu \mathrm{e}^{\bsi \varphi} \cos(\varphi_0-\varphi) - \bsi \lambda  \sin\varphi_0 + \mu  \cos\varphi_0\right] \bigg\}\bigg\}.
\end{aligned}
\end{equation}
Substituting \eqref{eq:J2} into \eqref{eq:J1}, we can obtain that
\begin{equation}\label{eq:T1 132}
\begin{aligned}
 & T_1( \bmf{u})\big |_{\Gamma_h^+ } =  \sum_{m=0} ^\infty  \bigg \{ \frac{\bsi k_p^2}{2} a_m \mathrm{e}^{\bsi (m-1) \varphi} \mathrm{e}^{-\bsi (\varphi-\varphi_0)} \mu
  J_{m-2} \left(k_p r\right)+ k_p^2 a_m \mathrm{e}^{\bsi m \varphi} \left(\lambda+\mu\right) \sin\varphi_0  J_m \left(k_p r\right) \\
& -\frac{\bsi k_p^2}{2} a_m \mathrm{e}^{\bsi (m+1) \varphi} \mathrm{e}^{\bsi (\varphi-\varphi_0)} \mu J_{m+2}
 -\frac{k_s^2}{2} b_m \mathrm{e}^{\bsi (m-1) \varphi} \mathrm{e}^{-\bsi (\varphi-\varphi_0)} \mu J_{m-2} \left(k_s r\right)\\
  &-\frac{k_s^2}{2} b_m \mathrm{e}^{\bsi (m+1) \varphi} \mathrm{e}^{\bsi (\varphi-\varphi_0)} \mu J_{m+2} \left(k_s r\right)\bigg\}.
\end{aligned}
\end{equation}

Similarly, substituting \eqref{eq:u1 partial 127} to \eqref{eq:u1u2 129} into \eqref{eq:Tu new exp}, after tedious but straightforward calculations, we have
\begin{equation}\label{eq:J3}
\begin{aligned}
& T_2(\bmf{u}  )\big|_{\Gamma_h^+ }=  \sum_{m=0} ^\infty  \left\{\frac{k_p^2}{4} a_m \mathrm{e}^{\bsi \left(m-1\right) \varphi} J_{m-2}(k_p r) \left[2 \bsi \mu \sin(\varphi-\varphi_0)+\lambda \mathrm{e}^{\bsi \varphi} \cos \varphi_0 +\bsi \mu \mathrm{e}^{\bsi \varphi} \sin \varphi_0\right] \right.\\
& + \frac{k_p^2}{2} a_m \mathrm{e}^{\bsi m \varphi} J_m(k_p r)\left[-\lambda \cos\varphi_0+2 \mu \sin\varphi \sin(\varphi_0-\varphi)\right]\\
& + \frac{k_p^2}{4} a_m \mathrm{e}^{\bsi m \varphi} J_{m+2}(k_p r) \left[2 \bsi \mu \mathrm{e}^{\bsi \varphi} \sin(\varphi_0-\varphi) + \lambda \cos\varphi_0 - \bsi \mu \sin \varphi_0 \right]\\
& + \frac{ k_s^2}{4} b_m \mathrm{e}^{\bsi \left(m-1\right) \varphi} J_{m-2}(k_s r) \left[2 \mu \sin(\varphi_0-\varphi)+ \bsi \lambda \mathrm{e}^{\bsi \varphi} \cos \varphi_0 - \mu \mathrm{e}^{\bsi \varphi} \sin \varphi_0\right]\\
& + \frac{ k_s^2}{2} b_m \mathrm{e}^{\bsi m\varphi} J_m(k_s r)\left[\mu \sin\varphi_0+2 \mu \cos\varphi \sin(\varphi-\varphi_0)\right] \\
& + \frac{ k_s^2}{4} b_m \mathrm{e}^{\bsi m \varphi} J_{m+2}(k_s r) \left[2 \mu \mathrm{e}^{\bsi \varphi} \sin(\varphi_0-\varphi) -\bsi \lambda \cos\varphi_0 - \mu \sin \varphi_0\right]\\
& + \frac{1}{r} \left\{-\frac{(m-1)k_p}{2} a_m \mathrm{e}^{\bsi (m-1)\varphi} J_{m-1}(k_p r)\left[2  \mu \cos(\varphi_0-\varphi) + \lambda \mathrm{e}^{\bsi \varphi} \cos\varphi_0 + \bsi \mu \mathrm{e}^{\bsi \varphi} \sin\varphi_0\right] \right.\\
\end{aligned}
\end{equation}
\[
\begin{aligned}
&\hspace*{-5mm} + \frac{(m+1)k_p}{2} a_m \mathrm{e}^{\bsi m \varphi} J_{m+1}(k_p r)\left[-2 \mu \mathrm{e}^{\bsi \varphi} \cos(\varphi_0-\varphi) - \lambda  \cos\varphi_0 + \bsi \mu  \sin\varphi_0\right]\\
&\hspace*{-5mm} + \frac{(m-1)k_s}{2} b_m \mathrm{e}^{\bsi (m-1)\varphi} J_{m-1}(k_s r)\left[-2 \bsi \mu \cos(\varphi_0-\varphi) - \bsi \lambda \mathrm{e}^{\bsi \varphi} \cos\varphi_0 + \mu \mathrm{e}^{\bsi \varphi} \sin\varphi_0\right]\\
 &\hspace*{-5mm}+ \frac{(m+1)k_s}{2} b_m \mathrm{e}^{\bsi m\varphi} J_{m+1}(k_s r) \left[2 \bsi \mu \mathrm{e}^{\bsi \varphi} \cos(\varphi_0-\varphi) + \bsi \lambda  \cos\varphi_0 + \mu  \sin\varphi_0\right] \bigg\}\bigg\},
\end{aligned}
\]
Plugging \eqref{eq:J2} into \eqref{eq:J3}, we can obtain that
\begin{equation}\label{eq:T2 134}
\begin{aligned}
&  T_2( \bmf{u})\big |_{\Gamma_h^+ } =  \sum_{m=0} ^\infty  \left\{- \frac{ k_p^2}{2} a_m \mathrm{e}^{\bsi (m-1) \varphi} \mathrm{e}^{-\bsi (\varphi-\varphi_0)} \mu
  J_{m-2} \left(k_p r\right)- k_p^2 a_m \mathrm{e}^{\bsi m \varphi} \left(\lambda+\mu\right) \cos\varphi_0  J_m \left(k_p r\right)\right.\\
& -\frac{ k_p^2}{2} a_m \mathrm{e}^{\bsi (m+1) \varphi} \mathrm{e}^{\bsi (\varphi-\varphi_0)} \mu J_{m+2}(k_p r)
 -\frac{\bsi k_s^2}{2} b_m \mathrm{e}^{\bsi (m-1) \varphi} \mathrm{e}^{-\bsi (\varphi-\varphi_0)} \mu J_{m-2} \left(k_s r\right)\\
  &+\frac{\bsi k_s^2}{2} b_m \mathrm{e}^{\bsi (m+1) \varphi} \mathrm{e}^{\bsi (\varphi-\varphi_0)} \mu J_{m+2} \left(k_s r\right)\bigg\}.
\end{aligned}
\end{equation}
Using the fact
 \[\left[\begin{array}{c}{\sin \varphi}\\{-\cos \varphi}\end{array}\right]=\frac{\bsi}{2}\left(\mathrm{e}^{-\bsi \varphi}\mathbf{e}_1-\mathrm{e}^{\bsi \varphi}\mathbf{e}_2\right),
 \]
substituting \eqref{eq:T1 132} and \eqref{eq:T2 134} into \eqref{eq:Tu new exp}, we can prove \eqref{eq:Tu1}.
The  proof of \eqref{eq:Tu2} is similar to \eqref{eq:Tu1}, which is omitted here.
\end{proof}


\begin{thebibliography}{99}

\bibitem{AR} G.~Alessandrini and L.~Rondi, {\it Determining a sound-soft polyhedral scatterer by a single far-field measurement}, Proc. Amer. Math. Soc., {\bf 35} (2005), 1685--1691.



\bibitem{Abr}
M.~Abramowitz and I. ~A.~Stegun, {\it Handbook of mathematical functions: with formulas, graphs, and mathematical tables}, vol. 55, Courier Corporation, 1964.
\bibitem{CY}
J.~Cheng and M.~Yamamoto, {\it Uniqueness in an inverse scattering problem within non-trapping polygonal obstacles with at most two incoming waves},
Inverse Problems, {\bf 19} (2003), 1361--1384.

\bibitem{EBL}
E. Bl{\aa}sten and Y.-H.  Lin, \emph{Radiating and non-radiating sources in elasticity}, Inverse Problems, 35, 1 (2019) 015005.



\bibitem{CDLZ}
X.~Cao, H.~Diao, H.~Liu and J.~Zou, \emph{On nodal and generalized singular structures of laplacian eigenfunctions and applications}, arXiv:1902.05798, 2019.

\bibitem{CDLZ2}
 X. Cao, H. Diao, H. Liu and J. Zou, {\it On nodal and generalized singular structures of Laplacian eigenfunctions and applications in $\mathbb{R}^3$}, arXiv:1909.10174, 2019. 



\bibitem{CK}
D.~Colton and R.~Kress, {\it Inverse Acoustic and Electromagnetic Scattering Theory}, 3rd edition, Springer-Verlag, Berlin, 2013.



\bibitem{CK18} D. Colton and R. Kress, {\it Looking back on inverse scattering theory}, SIAM Review, {\bf 60} (2018), no. 40, 779--807.

\bibitem{DR95}
G. Dassios and Z. Rigou, \emph{ Elastic Herglotz functions}, SIAM J. Appl. Math., {\bf 55}  (1995),  1345--1361.



\bibitem{ElschnerYama2010}

J. Elschner and M. Yamamoto, \emph{ Uniqueness in inverse elastic scattering with finitely many incident waves}, Inverse Problems, {\bf 26}  (2010), 045005.


\bibitem{krantz}
S.~G. Krantz and  H.~R. Parks, \emph{ A primer of real analytic functions}, 2nd edition, Birkh\"{a}user Boston, Inc., Boston, MA, 2002.


\bibitem{Kupradze}
V. D. Kupradze, \emph{Three-dimensional Problems of the Mathematical Theory of Elasticity and Thermoelasticity}, Amsterdam, North-Holland, 1979.





\bibitem{LPRX} H. Liu, M. Petrini, L. Rondi, Luca and J. Xiao, {Stable determination of sound-hard polyhedral scatterers by a minimal number of scattering measurements}, {\em J. Differential Equations}, {\bf 262} (2017), no. 3, 1631--1670.

\bibitem{LRX} {H. Liu, L. Rondi and J. Xiao}, {Mosco convergence for $H(curl)$ spaces, higher integrability for Maxwell's equations, and stability in direct and inverse EM scattering problems}, {\em J. Eur. Math. Soc. (JEMS)}, \textbf{21} (2019), 2945--2993.

\bibitem{LiuXiao} H. Liu and J. Xiao, Decoupling elastic waves and its applications, {\em J. Differential Equations}, \textbf{265} (2017), no. 8, 4442--4480.

\bibitem{Liu-Zou}
H.~Liu and J.~Zou,
{\it Uniqueness in an inverse acoustic obstacle scattering problem for both sound-hard and sound-soft polyhedral scatterers},
Inverse Problems, {\bf 22} (2006),
515--524.

\bibitem{Liu-Zou3}
H.~Liu and J.~Zou,
{\it On unique determination of partially coated polyhedral scatterers
with far field measurements},
Inverse Problems, {\bf 23} (2007), 297--308.



\bibitem{SP}
V. Sevroglou and G. Pelekanos,
\emph{Two-dimensional elastic Herglotz functions and their application in inverse scattering}, Journal of Elasticity, {\bf 68} (2002),  123--144.

\bibitem{TF} F. Treves, \emph{Introduction to Pseudodifferential and Fourier Integral Operators}, Vol. 1, Plenum Press, New York, 1980. 


\end{thebibliography}
\end{document}